\documentclass[11pt]{amsart}
\usepackage[utf8]{inputenc}
\usepackage{amsfonts,amssymb,amsmath,mathtools}
\usepackage{geometry}
\usepackage{indentfirst}
\usepackage{enumitem}
\usepackage{xcolor}
\usepackage{hyperref}
\usepackage{bookmark}
\setlist[itemize]{leftmargin=2em}
\numberwithin{equation}{section}
\newtheorem{theorem}{Theorem}

\newtheorem{lemma}{Lemma}[section]
\newtheorem{proposition}[lemma]{Proposition}
\theoremstyle{definition}

\theoremstyle{definition}
\newtheorem{remark}[lemma]{Remark}

\newcommand{\T}{\mathbb T}
\newcommand{\R}{\mathbb R}
\newcommand{\Z}{\mathbb Z}
\newcommand{\Id}{\mathrm{Id}}

\DeclareMathOperator{\Leb}{Leb}
\DeclareMathOperator{\supp}{supp}
\DeclareMathOperator{\dist}{dist}

\DeclareMathOperator*{\esssup}{ess\,sup}
\setlist[enumerate]{label=\textup{(\roman*)},leftmargin=*}
\hypersetup{pdftitle={Nondifferentiability on Irrational Invariant Graphs of Smooth Area-Preserving Twist Maps}}

\newcommand{\SinaiCutoffBW}{%
  \begingroup\normalfont\normalsize
  \ifdefined\color\color{black}\fi
  \setlength{\unitlength}{1bp}%
  \begin{picture}(420,161)
  \put(0,0){\hbox{\pdfliteral{q 0 G 0 g 1 j 0 J
0.55 w [] 0 d
37 34 m 409 34 l S
[] 0 d 412 34 m 408 35.36 l 408 32.64 l h f
0.55 w [] 0 d
37 34 m 37 121 l S
[] 0 d 37 124 m 35.64 120 l 38.36 120 l h f
0.42 w [1 2.7] 0 d
116 36 m 116 111 l S
0.6 w [] 0 d
116 34 m 116 31 l S
0.42 w [1 2.7] 0 d
208 36 m 208 111 l S
0.6 w [] 0 d
208 34 m 208 31 l S
0.42 w [1 2.7] 0 d
300 36 m 300 111 l S
0.6 w [] 0 d
300 34 m 300 31 l S
0.55 w [] 0 d
34 34 m 37 34 l S
0.55 w [] 0 d
34 73 m 37 73 l S
0.55 w [] 0 d
34 112 m 37 112 l S
0.95 w [] 0 d
37 112 m 300 112 l S
0.95 w [1.4 2.1] 0 d
116 112 m 121.09 111.98 l 123.28 111.7 l 124.73 111.1 l 125.46 110.61 l 126.91 109.18 l 127.64 108.24 l 128.37 107.15 l 129.82 104.51 l 132 99.5 l 134.19 93.41 l 137.82 82.21 l 139.28 78.55 l 140.01 77.17 l 140.74 76.1 l 142.92 73.77 l 145.83 70.97 l 149.47 67.76 l 154.56 63.72 l 159.65 60.16 l 164.74 57.03 l 170.56 53.91 l 176.38 51.21 l 182.2 48.88 l 188.75 46.63 l 195.3 44.72 l 201.84 43.1 l 209.12 41.58 l 216.39 40.32 l 224.4 39.17 l 233.13 38.16 l 241.86 37.34 l 251.31 36.64 l 262.23 36.01 l 285.51 35.12 l 314.61 34.54 l 352.44 34.21 l 407 34.05 l S
0.95 w [5.6 2.7] 0 d
208 112 m 213.47 111.96 l 215.46 111.65 l 216.46 111.25 l 217.45 110.61 l 218.45 109.71 l 219.44 108.51 l 221.43 105.27 l 222.43 103.24 l 223.92 99.72 l 225.91 94.24 l 229.89 82.02 l 231.38 78.34 l 232.88 75.93 l 235.86 72.84 l 239.84 69.16 l 243.82 65.82 l 248.3 62.45 l 252.78 59.43 l 257.25 56.74 l 262.23 54.08 l 267.2 51.73 l 272.18 49.65 l 277.15 47.82 l 282.62 46.05 l 288.6 44.38 l 294.56 42.94 l 301.03 41.6 l 307.5 40.47 l 314.47 39.43 l 321.93 38.51 l 329.89 37.69 l 338.35 36.99 l 357.25 35.86 l 379.64 35.06 l 407 34.54 l S
0.95 w [] 0 d
300 112 m 305.08 111.98 l 307.22 111.72 l 308.29 111.33 l 309.1 110.87 l 310.17 109.99 l 310.97 109.12 l 312.04 107.66 l 313.11 105.88 l 314.18 103.77 l 315.51 100.73 l 317.65 94.98 l 322.2 81.16 l 323.27 78.57 l 324.34 76.65 l 325.41 75.31 l 327.02 73.67 l 332.37 68.7 l 337.18 64.76 l 342.26 61.08 l 348.69 57.06 l 355.37 53.51 l 362.6 50.28 l 370.35 47.41 l 378.38 44.97 l 387.2 42.79 l 396.57 40.96 l 407 39.36 l S
0.55 w [] 0 d
40 69 m 113 69 l S
40 69 m 40 73 l S
113 69 m 113 73 l S
1 w [1.4 2.1] 0 d
40 150 m 62 150 l S
1 w [5.6 2.7] 0 d
174 150 m 196 150 l S
1 w [] 0 d
300 150 m 322 150 l S
Q}}}%
\put(39,130){\makebox(0,0)[l]{$Z(D(t)/R)$}}%
\put(116,23){\makebox(0,0){$R_{n+1}$}}%
\put(208,23){\makebox(0,0){$R_n$}}%
\put(300,23){\makebox(0,0){$R_{n-1}$}}%
\put(220,8){\makebox(0,0){$D(t)$ (logarithmic scale)}}%
\put(30,34){\makebox(0,0)[r]{$0$}}%
\put(30,73.0){\makebox(0,0)[r]{$1/2$}}%
\put(30,112){\makebox(0,0)[r]{$1$}}%
\put(77,80){\makebox(0,0){$E_j=0$}}%
\put(68,150){\makebox(0,0)[l]{$R=R_{n+1}$}}%
\put(202,150){\makebox(0,0)[l]{$R=R_n$}}%
\put(328,150){\makebox(0,0)[l]{$R=R_{n-1}$}}%
\end{picture}\endgroup}

\newcommand{\SinaiNodesBW}{%
  \begingroup\normalfont\normalsize
  \ifdefined\color\color{black}\fi
  \setlength{\unitlength}{1bp}%
  \begin{picture}(420,141)
  \put(0,0){\hbox{\pdfliteral{q 0 G 0 g 1 j 0 J
0.4 w [1 3] 0 d
117 109 m 117 118 l S
0.4 w [1 3] 0 d
199 109 m 199 118 l S
199 55 m 199 69 l S
199 86 m 199 92 l S
0.4 w [1 3] 0 d
282 109 m 282 118 l S
282 55 m 282 69 l S
282 86 m 282 92 l S
0.4 w [1 3] 0 d
365 109 m 365 118 l S
0.45 w [] 0 d
103 98 m 199 98 l S
365 98 m 396 98 l S
103 48 m 117 48 l S
282 48 m 396 48 l S
1.25 w [] 0 d
199 98 m 365 98 l S
117 48 m 282 48 l S
201.85 98 m 201.85 99.57 200.57 100.85 199 100.85 c 197.43 100.85 196.15 99.57 196.15 98 c 196.15 96.43 197.43 95.15 199 95.15 c 200.57 95.15 201.85 96.43 201.85 98 c h f
367.85 98 m 367.85 99.57 366.57 100.85 365 100.85 c 363.43 100.85 362.15 99.57 362.15 98 c 362.15 96.43 363.43 95.15 365 95.15 c 366.57 95.15 367.85 96.43 367.85 98 c h f
119.85 48 m 119.85 49.57 118.57 50.85 117 50.85 c 115.43 50.85 114.15 49.57 114.15 48 c 114.15 46.43 115.43 45.15 117 45.15 c 118.57 45.15 119.85 46.43 119.85 48 c h f
284.85 48 m 284.85 49.57 283.57 50.85 282 50.85 c 280.43 50.85 279.15 49.57 279.15 48 c 279.15 46.43 280.43 45.15 282 45.15 c 283.57 45.15 284.85 46.43 284.85 48 c h f
0.85 w [] 0 d
q 1 g 282 101.55 m 285.55 98 l 282 94.45 l 278.45 98 l h B Q
q 1 g 199 51.55 m 202.55 48 l 199 44.45 l 195.45 48 l h B Q
0.8 w [] 0 d
96.8 9 m 96.8 10.55 95.55 11.8 94 11.8 c 92.45 11.8 91.2 10.55 91.2 9 c 91.2 7.45 92.45 6.2 94 6.2 c 95.55 6.2 96.8 7.45 96.8 9 c h f
q 1 g 178 12.55 m 181.55 9 l 178 5.45 l 174.45 9 l h B Q
Q}}}%
\put(117,129){\makebox(0,0){$q_{2n}$}}%
\put(199,129){\makebox(0,0){$k_n$}}%
\put(282,129){\makebox(0,0){$q_{2n+2}$}}%
\put(365,129){\makebox(0,0){$k_{n+1}$}}%
\put(409,129){\makebox(0,0){$j$}}%
\put(9,98){\makebox(0,0)[l]{First orbit}}%
\put(9,48){\makebox(0,0)[l]{Second orbit}}%
\put(282,77){\makebox(0,0){$U_{q_{2n+2}}\asymp1$}}%
\put(199,29){\makebox(0,0){$B_{k_n}=e^{-M_n}$}}%
\put(103,9){\makebox(0,0)[l]{cutoff index}}%
\put(188,9){\makebox(0,0)[l]{term used in the sum}}%
\end{picture}\endgroup}

\newcommand{\SinaiChordsBW}{%
  \begingroup\normalfont\normalsize
  \ifdefined\color\color{black}\fi
  \setlength{\unitlength}{1bp}%
  \begin{picture}(420,149)
  \put(0,0){\hbox{\pdfliteral{q 0 G 0 g 1 j 0 J
0.7 w [] 0 d
129 72 m 110.66 112.36 l S
[] 0 d 109 116 m 109.52 110.74 l 112.62 112.15 l h f
0.7 w [4 2.7] 0 d
215 72 m 233.34 112.36 l S
[] 0 d 235 116 m 231.38 112.15 l 234.48 110.74 l h f
0.7 w [] 0 d
129 66 m 101.19 23.35 l S
[] 0 d 99 20 m 103.16 23.26 l 100.31 25.12 l h f
0.7 w [4 2.7] 0 d
215 66 m 242.81 23.35 l S
[] 0 d 245 20 m 243.69 25.12 l 240.84 23.26 l h f
0.5 w [] 0 d
63 120 m 290 120 l S
63 117 m 63 123 l S
290 117 m 290 123 l S
1.55 w [] 0 d
109 120 m 235 120 l S
0.82 w [] 0 d
111.75 120 m 111.75 121.52 110.52 122.75 109 122.75 c 107.48 122.75 106.25 121.52 106.25 120 c 106.25 118.48 107.48 117.25 109 117.25 c 110.52 117.25 111.75 118.48 111.75 120 c h f
q 1 g 232.3 117.3 5.4 5.4 re B Q
0.5 w [] 0 d
67 69 m 301 69 l S
67 66 m 67 72 l S
301 66 m 301 72 l S
1.55 w [] 0 d
129 69 m 215 69 l S
0.82 w [] 0 d
131.75 69 m 131.75 70.52 130.52 71.75 129 71.75 c 127.48 71.75 126.25 70.52 126.25 69 c 126.25 67.48 127.48 66.25 129 66.25 c 130.52 66.25 131.75 67.48 131.75 69 c h f
q 1 g 212.3 66.3 5.4 5.4 re B Q
0.5 w [] 0 d
55 16 m 282 16 l S
55 13 m 55 19 l S
282 13 m 282 19 l S
1.55 w [] 0 d
99 16 m 245 16 l S
0.82 w [] 0 d
101.75 16 m 101.75 17.52 100.52 18.75 99 18.75 c 97.48 18.75 96.25 17.52 96.25 16 c 96.25 14.48 97.48 13.25 99 13.25 c 100.52 13.25 101.75 14.48 101.75 16 c h f
q 1 g 242.3 13.3 5.4 5.4 re B Q
Q}}}%
\put(34,120){\makebox(0,0){$I_{j+1}$}}%
\put(355,120){\makebox(0,0){$\Delta_hX_{j+1}$}}%
\put(34,69){\makebox(0,0){$I_j$}}%
\put(355,69){\makebox(0,0){$\Delta_hX_j$}}%
\put(34,16){\makebox(0,0){$I_{j-1}$}}%
\put(355,16){\makebox(0,0){$\Delta_hX_{j-1}$}}%
\put(109,137){\makebox(0,0){$t$}}%
\put(235,137){\makebox(0,0){$t+h$}}%
\put(355,138){\makebox(0,0){Chord length}}%
\put(258,94){\makebox(0,0){$g$}}%
\put(269,43){\makebox(0,0){$g^{-1}$}}%
\end{picture}\endgroup}

\AtBeginDocument{\hypersetup{hidelinks}}
\newcommand{\SinaiOverviewCirclesBW}{%
\begingroup\normalfont\normalsize\color{black}%
\setlength{\unitlength}{1bp}%
\begin{picture}(430,158)
\put(0,0){\hbox{\pdfliteral{q 0 G 0 g 1 J 1 j
1 0 0 1 0 -422 cm
0.55 w [] 0 d 126.000 510.000 m 126.000 528.629 114.005 545.139 96.288 550.895 c 78.571 556.652 59.162 550.346 48.212 535.275 c 37.263 520.204 37.263 499.796 48.212 484.725 c 59.162 469.654 78.571 463.348 96.288 469.105 c 114.005 474.861 126.000 491.371 126.000 510.000 c S
0 g 125.269 526.108 m 125.269 527.434 124.194 528.508 122.869 528.508 c 121.543 528.508 120.469 527.434 120.469 526.108 c 120.469 524.783 121.543 523.708 122.869 523.708 c 124.194 523.708 125.269 524.783 125.269 526.108 c h f 0 g 0.6 w 125.269 526.108 m 125.269 527.434 124.194 528.508 122.869 528.508 c 121.543 528.508 120.469 527.434 120.469 526.108 c 120.469 524.783 121.543 523.708 122.869 523.708 c 124.194 523.708 125.269 524.783 125.269 526.108 c h S
0 g 51.516 536.473 m 51.516 537.799 50.441 538.873 49.116 538.873 c 47.790 538.873 46.716 537.799 46.716 536.473 c 46.716 535.148 47.790 534.073 49.116 534.073 c 50.441 534.073 51.516 535.148 51.516 536.473 c h f 0 g 0.6 w 51.516 536.473 m 51.516 537.799 50.441 538.873 49.116 538.873 c 47.790 538.873 46.716 537.799 46.716 536.473 c 46.716 535.148 47.790 534.073 49.116 534.073 c 50.441 534.073 51.516 535.148 51.516 536.473 c h S
0 g 79.416 467.418 m 79.416 468.744 78.341 469.818 77.016 469.818 c 75.690 469.818 74.616 468.744 74.616 467.418 c 74.616 466.093 75.690 465.018 77.016 465.018 c 78.341 465.018 79.416 466.093 79.416 467.418 c h f 0 g 0.6 w 79.416 467.418 m 79.416 468.744 78.341 469.818 77.016 469.818 c 75.690 469.818 74.616 468.744 74.616 467.418 c 74.616 466.093 75.690 465.018 77.016 465.018 c 78.341 465.018 79.416 466.093 79.416 467.418 c h S
1 g 91.384 552.582 m 91.384 553.907 90.310 554.982 88.984 554.982 c 87.659 554.982 86.584 553.907 86.584 552.582 c 86.584 551.256 87.659 550.182 88.984 550.182 c 90.310 550.182 91.384 551.256 91.384 552.582 c h f 0 g 0.6 w 91.384 552.582 m 91.384 553.907 90.310 554.982 88.984 554.982 c 87.659 554.982 86.584 553.907 86.584 552.582 c 86.584 551.256 87.659 550.182 88.984 550.182 c 90.310 550.182 91.384 551.256 91.384 552.582 c h S
1 g 45.531 493.892 m 45.531 495.217 44.457 496.292 43.131 496.292 c 41.806 496.292 40.731 495.217 40.731 493.892 c 40.731 492.566 41.806 491.492 43.131 491.492 c 44.457 491.492 45.531 492.566 45.531 493.892 c h f 0 g 0.6 w 45.531 493.892 m 45.531 495.217 44.457 496.292 43.131 496.292 c 41.806 496.292 40.731 495.217 40.731 493.892 c 40.731 492.566 41.806 491.492 43.131 491.492 c 44.457 491.492 45.531 492.566 45.531 493.892 c h S
1 g 119.284 483.527 m 119.284 484.852 118.210 485.927 116.884 485.927 c 115.559 485.927 114.484 484.852 114.484 483.527 c 114.484 482.201 115.559 481.127 116.884 481.127 c 118.210 481.127 119.284 482.201 119.284 483.527 c h f 0 g 0.6 w 119.284 483.527 m 119.284 484.852 118.210 485.927 116.884 485.927 c 115.559 485.927 114.484 484.852 114.484 483.527 c 114.484 482.201 115.559 481.127 116.884 481.127 c 118.210 481.127 119.284 482.201 119.284 483.527 c h S
0.5 w [] 0 d 100.443 557.924 m 87.443 562.656 73.078 561.903 60.643 555.838 c S
0.5 w [] 0 d 62.890 556.934 m 60.643 555.838 l S
0 g 60.643 555.838 m 63.603 558.895 l 64.874 556.289 l h f
0.45 w [] 0 d 386.000 510.000 m 386.000 528.629 374.005 545.139 356.288 550.895 c 338.571 556.652 319.162 550.346 308.212 535.275 c 297.263 520.204 297.263 499.796 308.212 484.725 c 319.162 469.654 338.571 463.348 356.288 469.105 c 374.005 474.861 386.000 491.371 386.000 510.000 c S
2 w [] 0 d 385.679 515.240 m 384.729 522.985 381.688 530.325 376.884 536.473 c S
2 w [] 0 d 317.122 544.341 m 311.500 540.105 307.004 534.553 304.029 528.173 c S
2 w [] 0 d 326.199 470.418 m 332.679 467.668 339.735 466.550 346.748 467.164 c S
1.6 w [3 2] 0 d 359.801 549.582 m 353.321 552.332 346.265 553.450 339.252 552.836 c S
1.6 w [3 2] 0 d 300.321 504.760 m 301.178 497.772 303.739 491.103 307.776 485.336 c S
1.6 w [3 2] 0 d 368.878 475.659 m 374.500 479.895 378.996 485.447 381.971 491.827 c S
0.55 w [] 0 d 286.000 510.000 m 140.000 510.000 l S
0 g 140.000 510.000 m 144.000 511.450 l 144.000 508.550 l h f
0.4 w [] 0 d 383.000 527.000 m 398.000 536.000 l S
0.4 w [] 0 d 378.000 480.000 m 395.000 473.000 l S
0.4 w [] 0 d 331.000 502.000 m 307.776 485.337 l S
0 g 308.776 485.337 m 308.776 485.889 308.328 486.337 307.776 486.337 c 307.224 486.337 306.776 485.889 306.776 485.337 c 306.776 484.785 307.224 484.337 307.776 484.337 c 308.328 484.337 308.776 484.785 308.776 485.337 c h f 0 g 0.6 w 308.776 485.337 m 308.776 485.889 308.328 486.337 307.776 486.337 c 307.224 486.337 306.776 485.889 306.776 485.337 c 306.776 484.785 307.224 484.337 307.776 484.337 c 308.328 484.337 308.776 484.785 308.776 485.337 c h S
Q}}}
\put(83,150){\makebox(0,0)[c]{Rotation circle $\T$}}
\put(343,150){\makebox(0,0)[c]{Circle of $g$: $\T$}}
\put(83,88){\makebox(0,0)[c]{$\theta$}}
\put(28,126){\makebox(0,0)[c]{$R_\alpha$}}
\put(213,105){\makebox(0,0)[c]{$h_g$}}
\put(213,72){\makebox(0,0)[c]{collapse the fibres}}
\put(213,52){\makebox(0,0)[c]{$h_g\circ g=R_\alpha\circ h_g$}}
\put(343,92){\makebox(0,0)[c]{$x$}}
\put(327,78){\makebox(0,0)[c]{$K$}}
\put(400,118){\makebox(0,0)[l]{$I_j$}}
\put(397,48){\makebox(0,0)[l]{$\widehat I_j$}}
\put(83,29){\makebox(0,0)[c]{$\bullet\ j\alpha\qquad\circ\ (1/2+j\alpha)$}}
\put(343,29){\makebox(0,0)[c]{Only selected fibres are shown}}
\put(215,13){\makebox(0,0)[c]{$h_g(\overline I_j)=j\alpha,\qquad h_g(\overline{\widehat I_j})=1/2+j\alpha$}}
\end{picture}%
\endgroup}

\newcommand{\SinaiOverviewDensitiesBW}{%
\begingroup\normalfont\normalsize\color{black}%
\setlength{\unitlength}{1bp}%
\begin{picture}(430,223.04)
\put(0,0){\hbox{\pdfliteral{q 1 0 0 0.82 0 0 cm 0 G 0 g 1 J 1 j
1 0 0 1 0 -150 cm
0.6 w [] 0 d 22.000 322.000 m 195.000 322.000 l S
0.5 w [] 0 d 22.000 319.000 m 22.000 325.000 l S
0.5 w [] 0 d 195.000 319.000 m 195.000 325.000 l S
0.6 w [] 0 d 237.000 322.000 m 410.000 322.000 l S
0.5 w [] 0 d 237.000 319.000 m 237.000 325.000 l S
0.5 w [] 0 d 410.000 319.000 m 410.000 325.000 l S
2.1 w [] 0 d 78.000 322.000 m 80.259 322.000 l S
2.1 w [] 0 d 82.519 322.000 m 84.778 322.000 l S
2.1 w [] 0 d 91.556 322.000 m 93.815 322.000 l S
2.1 w [] 0 d 96.074 322.000 m 98.333 322.000 l S
2.1 w [] 0 d 118.667 322.000 m 120.926 322.000 l S
2.1 w [] 0 d 123.185 322.000 m 125.444 322.000 l S
2.1 w [] 0 d 132.222 322.000 m 134.481 322.000 l S
2.1 w [] 0 d 136.741 322.000 m 139.000 322.000 l S
0 g 103.000 322.000 m 103.000 323.105 102.105 324.000 101.000 324.000 c 99.895 324.000 99.000 323.105 99.000 322.000 c 99.000 320.895 99.895 320.000 101.000 320.000 c 102.105 320.000 103.000 320.895 103.000 322.000 c h f 0 g 0.6 w 103.000 322.000 m 103.000 323.105 102.105 324.000 101.000 324.000 c 99.895 324.000 99.000 323.105 99.000 322.000 c 99.000 320.895 99.895 320.000 101.000 320.000 c 102.105 320.000 103.000 320.895 103.000 322.000 c h S
1 g 324.200 322.000 m 324.200 323.215 323.215 324.200 322.000 324.200 c 320.785 324.200 319.800 323.215 319.800 322.000 c 319.800 320.785 320.785 319.800 322.000 319.800 c 323.215 319.800 324.200 320.785 324.200 322.000 c h f 0 g 0.6 w 324.200 322.000 m 324.200 323.215 323.215 324.200 322.000 324.200 c 320.785 324.200 319.800 323.215 319.800 322.000 c 319.800 320.785 320.785 319.800 322.000 319.800 c 323.215 319.800 324.200 320.785 324.200 322.000 c h S
0.6 w [] 0 d 27.000 242.000 m 193.000 242.000 l S
0.5 w [] 0 d 27.000 239.000 m 27.000 245.000 l S
0.5 w [] 0 d 193.000 239.000 m 193.000 245.000 l S
0.6 w [] 0 d 18.000 183.000 m 188.000 183.000 l S
0.5 w [] 0 d 18.000 180.000 m 18.000 186.000 l S
0.5 w [] 0 d 188.000 180.000 m 188.000 186.000 l S
2.1 w [] 0 d 85.000 242.000 m 87.148 242.000 l S
2.1 w [] 0 d 89.296 242.000 m 91.444 242.000 l S
2.1 w [] 0 d 97.889 242.000 m 100.037 242.000 l S
2.1 w [] 0 d 102.185 242.000 m 104.333 242.000 l S
2.1 w [] 0 d 123.667 242.000 m 125.815 242.000 l S
2.1 w [] 0 d 127.963 242.000 m 130.111 242.000 l S
2.1 w [] 0 d 136.556 242.000 m 138.704 242.000 l S
2.1 w [] 0 d 140.852 242.000 m 143.000 242.000 l S
2.1 w [] 0 d 72.000 183.000 m 74.259 183.000 l S
2.1 w [] 0 d 76.519 183.000 m 78.778 183.000 l S
2.1 w [] 0 d 85.556 183.000 m 87.815 183.000 l S
2.1 w [] 0 d 90.074 183.000 m 92.333 183.000 l S
2.1 w [] 0 d 112.667 183.000 m 114.926 183.000 l S
2.1 w [] 0 d 117.185 183.000 m 119.444 183.000 l S
2.1 w [] 0 d 126.222 183.000 m 128.481 183.000 l S
2.1 w [] 0 d 130.741 183.000 m 133.000 183.000 l S
0 g 110.200 242.000 m 110.200 243.215 109.215 244.200 108.000 244.200 c 106.785 244.200 105.800 243.215 105.800 242.000 c 105.800 240.785 106.785 239.800 108.000 239.800 c 109.215 239.800 110.200 240.785 110.200 242.000 c h f 0 g 0.6 w 110.200 242.000 m 110.200 243.215 109.215 244.200 108.000 244.200 c 106.785 244.200 105.800 243.215 105.800 242.000 c 105.800 240.785 106.785 239.800 108.000 239.800 c 109.215 239.800 110.200 240.785 110.200 242.000 c h S
0 g 101.200 183.000 m 101.200 184.215 100.215 185.200 99.000 185.200 c 97.785 185.200 96.800 184.215 96.800 183.000 c 96.800 181.785 97.785 180.800 99.000 180.800 c 100.215 180.800 101.200 181.785 101.200 183.000 c h f 0 g 0.6 w 101.200 183.000 m 101.200 184.215 100.215 185.200 99.000 185.200 c 97.785 185.200 96.800 184.215 96.800 183.000 c 96.800 181.785 97.785 180.800 99.000 180.800 c 100.215 180.800 101.200 181.785 101.200 183.000 c h S
0.7 w [] 0 d 108.000 237.000 m 99.000 188.000 l S
0 g 99.000 188.000 m 98.296 192.196 l 101.149 191.672 l h f
0.9 w [3 2] 0 d 244.000 242.000 m 402.000 242.000 l S
0.5 w [] 0 d 244.000 239.000 m 244.000 245.000 l S
0.5 w [] 0 d 402.000 239.000 m 402.000 245.000 l S
0.9 w [3 2] 0 d 235.000 183.000 m 410.000 183.000 l S
0.5 w [] 0 d 235.000 180.000 m 235.000 186.000 l S
0.5 w [] 0 d 410.000 180.000 m 410.000 186.000 l S
1 g 322.200 242.000 m 322.200 243.215 321.215 244.200 320.000 244.200 c 318.785 244.200 317.800 243.215 317.800 242.000 c 317.800 240.785 318.785 239.800 320.000 239.800 c 321.215 239.800 322.200 240.785 322.200 242.000 c h f 0 g 0.6 w 322.200 242.000 m 322.200 243.215 321.215 244.200 320.000 244.200 c 318.785 244.200 317.800 243.215 317.800 242.000 c 317.800 240.785 318.785 239.800 320.000 239.800 c 321.215 239.800 322.200 240.785 322.200 242.000 c h S
1 g 339.200 183.000 m 339.200 184.215 338.215 185.200 337.000 185.200 c 335.785 185.200 334.800 184.215 334.800 183.000 c 334.800 181.785 335.785 180.800 337.000 180.800 c 338.215 180.800 339.200 181.785 339.200 183.000 c h f 0 g 0.6 w 339.200 183.000 m 339.200 184.215 338.215 185.200 337.000 185.200 c 335.785 185.200 334.800 184.215 334.800 183.000 c 334.800 181.785 335.785 180.800 337.000 180.800 c 338.215 180.800 339.200 181.785 339.200 183.000 c h S
0.7 w [] 0 d 320.000 237.000 m 337.000 188.000 l S
0 g 337.000 188.000 m 334.319 191.304 l 337.059 192.254 l h f
Q}}}
\put(107,206.64){\makebox(0,0)[c]{First family: singular parameters}}
\put(323,206.64){\makebox(0,0)[c]{Second family: smooth densities}}
\put(107,189.42){\makebox(0,0)[c]{$\rho_j=P_j Z_{R_n}$}}
\put(323,189.42){\makebox(0,0)[c]{$\widehat\rho_j=\widehat P_j Z(\widehat D/\widehat R_n)$}}
\put(107,174.66){\makebox(0,0)[c]{$P_j=A_j(1+T_j\Theta/4)$ near $E$}}
\put(323,174.66){\makebox(0,0)[c]{$\widehat P_j=B_j\widehat H_j$}}
\put(215,157.44){\makebox(0,0)[c]{Prescribed at cutoff indices; recurrence between them}}
\put(22,128.74){\makebox(0,0)[c]{$0$}}
\put(195,128.74){\makebox(0,0)[c]{$1$}}
\put(237,128.74){\makebox(0,0)[c]{$0$}}
\put(410,128.74){\makebox(0,0)[c]{$1$}}
\put(107,128.74){\makebox(0,0)[c]{$t\in[0,1],\quad E$}}
\put(323,128.74){\makebox(0,0)[c]{$t\in[0,1]$}}
\put(107,110.7){\makebox(0,0)[c]{$x_j(t)=\int_0^t\rho_j(v)\,dv$}}
\put(323,110.7){\makebox(0,0)[c]{$\widehat x_j(t)=\int_0^t\widehat\rho_j(v)\,dv$}}
\put(107,94.3){\makebox(0,0)[c]{$X_j(t)=a(j\alpha)+x_j(t)/S$}}
\put(323,94.3){\makebox(0,0)[c]{$\widehat X_j(t)=a(1/2+j\alpha)+\widehat x_j(t)/S$}}
\put(28,83.64){\makebox(0,0)[l]{$I_j$}}
\put(18,35.26){\makebox(0,0)[l]{$I_{j+1}$}}
\put(244,83.64){\makebox(0,0)[l]{$\widehat I_j$}}
\put(235,35.26){\makebox(0,0)[l]{$\widehat I_{j+1}$}}
\put(125,52.48){\makebox(0,0)[l]{$g$}}
\put(83,52.48){\makebox(0,0)[r]{same $t$}}
\put(345,52.48){\makebox(0,0)[l]{$g$}}
\put(305,52.48){\makebox(0,0)[r]{same $t$}}
\put(107,9.84){\makebox(0,0)[c]{$X_j(E)$ and $X_{j+1}(E)$}}
\put(323,9.84){\makebox(0,0)[c]{Smooth interval maps}}
\end{picture}%
\endgroup}

\newcommand{\SinaiOverviewRegularityBW}{%
\begingroup\normalfont\normalsize\color{black}%
\setlength{\unitlength}{1bp}%
\begin{picture}(430,129)
\put(215,111.8){\makebox(0,0)[c]{$N(u)=\bigcup_{j\in\Z}X_j(E),\qquad |I_j|=w_j/S,\quad |\widehat I_j|=\widehat w_j/S$}}
\put(215,90.3){\makebox(0,0)[c]{$Q_j=\dfrac{\rho_{j+1}+\rho_{j-1}}{\rho_j}-2,\qquad \mathrm D\phi(X_j(t))=Q_j(t)$}}
\put(215,68.8){\makebox(0,0)[c]{$u=g-\mathrm{Id},\qquad \phi=g+g^{-1}-2\mathrm{Id}$}}
\put(215,48.16){\makebox(0,0)[c]{Near $X_j(E)$: exact affine cancellation}}
\put(215,30.1){\makebox(0,0)[c]{Across $K$: $C^{r-1}$ estimates yield $\phi\in C^r$}}
\put(215,12.04){\makebox(0,0)[c]{Cutoff transitions: $D(t)\asymp R_n$ and $\widehat D(t)\asymp\widehat R_n$}}
\end{picture}%
\endgroup}

\makeatletter
\newenvironment{SinaiAppendixFigure}{%
  \par\smallskip\noindent\begin{minipage}{\linewidth}%
  \def\@captype{figure}\centering
}{\end{minipage}\par\smallskip}
\makeatother

\begin{document}
\title[Nondifferentiability on Invariant Graphs]{Nondifferentiability on Irrational Invariant Graphs\\of Smooth Area-Preserving Twist Maps}
\author{Lin Wang}
\address{School of Mathematics and Statistics, Beijing Institute of Technology, Beijing 100081, China}
	\email{lwang@bit.edu.cn}

	

\begin{abstract}
For a fixed irrational rotation number and every $1\le r<\infty$, we construct an exact area-preserving $C^r$ twist map preserving a Lipschitz graph. The set of nondifferentiability points has Hausdorff dimension $d=2/(r+1)$ and zero $d$-dimensional Hausdorff measure. At the same rotation number, we construct a $C^\infty$ twist map whose invariant graph is nondifferentiable on exactly one two-sided orbit. The dynamics on these graphs has wandering intervals.

\end{abstract}

\maketitle
{\footnotesize\setlength{\baselineskip}{10pt}\vspace{-\baselineskip}\tableofcontents}
\newpage

\section{Introduction and main results}\label{sec:intro}

The existence of an invariant circle, its regularity, and the dynamics
on it are different questions for an area-preserving twist map.
Birkhoff's graph theorem gives Lipschitz regularity for an invariant
essential circle; see  \cite{Ar09} and \cite[Chap.~I]{Her}.
KAM theory gives more: under suitable smallness, regularity and
arithmetic assumptions, smooth invariant circles persist and carry
dynamics conjugate to irrational rotations
\cite{Arn63,Her,Her86,Kol54,Mos62}.
Aubry--Mather theory provides ordered minimizing invariant sets for
rotation numbers in the admissible range, without requiring a small
perturbation of an integrable map \cite{ALD83,Ban88,Mat82}.
For an irrational rotation number, the minimal set can be a circle or
a Cantor set. Mather's celebrated results \cite{Mat82,Mat84} show that a
complete invariant circle need not survive even though ordered
minimizing sets remain.

We study an intermediate problem between these two descriptions:
the complete invariant graph is retained, but its geometry is allowed
to be nondifferentiable and its irrational dynamics is allowed to have
wandering intervals. The distinction is geometric and dynamical; we do
not assume that the constructed maps belong to a specified parameter
family near a KAM breakdown threshold. Our aim is to prescribe the
Hausdorff dimension of the nondifferentiability set while increasing
the regularity of the twist map.

Arnaud constructed nondifferentiable invariant curves with Denjoy
dynamics for $C^1$ and $C^2$ symplectic twist maps \cite{Ar11,Ar14}.
Avila and Fayad \cite{AF} obtained a $C^1$ twist map with a
nondifferentiable invariant curve carrying minimal dynamics at any
prescribed irrational rotation number. Guo, Liu and Wang
\cite[Theorem~4]{GLW} constructed invariant graphs outside every
$C^{2,\beta}$ class, $\beta>0$, at a fixed constant-type rotation
number, using arbitrarily small $C^{3-\varepsilon}$ perturbations.
These results
separate the regularity of the map from that of its invariant graph.
The present construction treats every finite regularity and the
$C^\infty$ class, with wandering intervals, and determines the entire
set of nondifferentiability points.

The dissipative setting provides a related comparison. Hall
\cite{Hall83} constructed a smooth annulus diffeomorphism with an
attracting Lipschitz invariant circle carrying Denjoy dynamics,
without imposing a twist condition. Le Calvez \cite{LeC88} studied
the geometry and rotation properties of Birkhoff attractors for
dissipative twist maps. The dissipative counterpart of the
Herman--Mather formula already appears in an equivalent graph form
in \cite[\S11.5, p.~285]{LeC88}. Sorrentino and Wang
\cite[Proposition~2.1]{SW26} formulate this relation as a
characterization of invariant graphs and extend it to conformally
symplectic maps; they also destroy all continuous invariant Lagrangian
graphs by perturbations arbitrarily small in $C^{1-\varepsilon}$,
$0<\varepsilon<1$.
Huang, Li and Wang \cite{HLW26} study sharp perturbation thresholds
for graph regularity and normal hyperbolicity. Kong, Massetti and
Wang \cite{KMW26} give quantitative $C^r$ persistence results and,
in their one-dimensional dissipative setting with irrational rotation
number, upgrade an existing $C^1$ graph to the $C^k$ regularity of
the map, for $k=2,3,\ldots,\infty$. They also give analytic
dissipative twist maps with Lipschitz non-$C^1$ invariant graphs;
the regularity theorem therefore does not extend to arbitrary
Lipschitz graphs. Li and Wang \cite{LW26}
study persistence and internal rotation numbers, and construct
Denjoy graphs at a critical Lipschitz bound with Lipschitz
perturbations that are smooth outside a null Cantor set. These works use the relation
between the graph equation and the induced circle map. In the
area-preserving problem considered here, the transverse contraction
used in dissipative persistence arguments is unavailable.

The construction has to satisfy the following requirements together.
\begin{itemize}[leftmargin=*,itemsep=.5em]
\item \emph{Constructing a complete invariant graph.}
An ordered Aubry--Mather Cantor set does not provide invariant graph
segments in its complementary intervals. Their circular order,
lengths and two-sided dynamics must be compatible.
There is also a geometric restriction: Arnaud and Le Calvez
\cite[Theorem~1.2 and Corollary~1.3]{ALC17} exclude $C^2$ Denjoy
sub-systems. Here the $C^2$ requirement concerns the embedding of
the Cantor dynamics in a curve, not just the ambient map.
Thus a smooth ambient map does not permit us to place the
exceptional dynamics on an arbitrarily chosen smooth curve.

\item \emph{Keeping the map smooth while the graph is rough.}
For a standard twist map, the Herman--Mather formula expresses
the perturbation through the sum of an increasing circle lift and
its inverse. The two maps may fail to be differentiable, but their
symmetric sum must have all the derivatives required of the
perturbation. The local cancellations in \cite{Ar11,Ar14,AF}
give a starting point. To prescribe a larger singular set or obtain
one $C^\infty$ example, these cancellations must also be compatible
with summable interval lengths and with derivatives in the circle
coordinate.

\item \emph{Controlling points outside the support of the invariant measure.}
For $C^{1+\beta}$ symplectic twist maps, $\beta>0$, Arnaud and
Berger \cite[Theorem~1 and Corollary~1]{AB16} prove that the
invariant probability on an irrational invariant circle has zero
Lyapunov exponents and that the circle is $C^1$-regular almost
everywhere for this measure. Arnaud \cite{ArAnn11} relates zero
Lyapunov exponents of irrational Mather measures to almost-everywhere
$C^1$-regularity of their supports, using Green bundles and
paratingent cones. This regularity means that limiting secant
directions at the point lie in one line. Every wandering interval has zero invariant
measure. These results therefore do not control the secant slopes
at the points where our graph will be nondifferentiable.
Moreover, regularity of the perturbation inside each interval is
insufficient: all derivative estimates must remain compatible
as the intervals accumulate on the minimal Cantor set.
\end{itemize}
\vspace{1ex}\vspace{1ex}\paragraph{\bf \bf Main results.}

Let $\T=\R/\Z$. For a periodic real-valued function $v$, let
$N(v)\subset\T$ denote the points at which its periodic lift has no
 derivative. Given a $1$-periodic function $\phi:\R\to\R$, define
\begin{equation}\label{in:map}
 F_\phi:\R^2\longrightarrow\R^2,
 \qquad F_\phi(x,y):=(x+y,\;y+\phi(x+y)).
\end{equation}
The identity $F_\phi(x+1,y)=F_\phi(x,y)+(1,0)$ induces a map on
$\T\times\R$, denoted by the same symbol. For $u:\T\to\R$, also denoting
its periodic lift by $u$, define the graph embedding and projection by
\begin{equation}\label{pre:maps}
 \begin{aligned}
 \iota_u:\T&\longrightarrow\T\times\R,& \iota_u(x)&:=(x,u(x)),\\
 \pi_1:\T\times\R&\longrightarrow\T,& \pi_1(x,y)&:=x.
 \end{aligned}
\end{equation}
The graph is $\Gamma:=\iota_u(\T)$. If it is invariant under $F_\phi$,
the induced circle map has real lift
\[
 g:\R\longrightarrow\R,\qquad g(x):=x+u(x).
\]
Thus $\pi_1\circ F_\phi\circ\iota_u=g$ on the circle. We use the same
symbol for a degree-one lift and its induced circle map; formulas with
real differences are always interpreted on the lift.
Throughout the paper $C^k$, for an integer $k$, means continuous derivatives through order $k$. For $r=k+{\beta}$, $0<{\beta}<1$, $C^r$ means $C^{k,{\beta}}$. We use $C^0$ for continuity.

In both theorems, we fix the irrational rotation number
\begin{equation}\label{in:frequency}
 {\alpha}=[0;a_1,a_2,\ldots],\qquad
 a_{2n+1}=2(n+2)^4,\quad a_{2n+2}=2\qquad(n\ge0).
\end{equation}

\begin{theorem}[Finite regularity]\label{thm:finite}
For every $1\le r<\infty$ there is a periodic $\phi\in C^r$ of mean zero such that $F_\phi$ preserves the graph of a Lipschitz function $u:\T\to\R$ with induced rotation number ${\alpha}$. Its set of nondifferentiability points satisfies
\[
 \dim_H N(u)=d:=\frac{2}{r+1},\qquad \mathcal H^d(N(u))=0.
\]
\end{theorem}

The map in Theorem~\ref{thm:finite} may depend on $r$.
Consequently, letting $r$ tend to infinity does not give one smooth
example. For the next result, we take the singular parameter set to
be a singleton and choose one family of cutoff radii satisfying the
estimates at every derivative order.

\begin{theorem}[Smooth perturbation]\label{thm:smooth}
There is a periodic $\phi\in C^\infty$ of mean zero such that $F_\phi$ preserves the graph of a Lipschitz function $u:\T\to\R$ with induced rotation number ${\alpha}$ and $N(u)$ is exactly the two-sided orbit of one point in a wandering interval. In particular the invariant graph is not $C^1$.
\end{theorem}

At $r=1$, Theorem~\ref{thm:finite} gives $\dim_H N(u)=1$ and
zero Lebesgue measure.
Arnaud's theorem \cite[Theorem~1]{Ar09} gives differentiability on a
dense $G_\delta$ set of full Lebesgue measure. These conclusions do
not imply a Hausdorff dimension bound below one for $N(u)$.
Theorem~\ref{thm:smooth} gives one smooth map and identifies its
set of nondifferentiability points.

\begin{remark}[Sharpness in the locally affine class]
The exponent $2/(r+1)$ is sharp within the class of standard twist
maps whose perturbation is locally affine near the orbit of every
nondifferentiability point in the wandering interval.  Proposition~\ref{sharp:upper}
proves the upper bound, including zero Hausdorff measure at the critical
dimension, for every finite $r\ge1$. By Lemma~\ref{gl:local-affine}, the examples of
Theorem~\ref{thm:finite} belong to this class and attain the upper
bound at the fixed frequency~\eqref{in:frequency}.
Proposition~\ref{sharp:optimal} also shows that the perturbation constructed
for a given $r$ does not belong to $C^s$ for any $s>r$.
\end{remark}

\begin{remark}[Admissible rotation numbers]\label{rem:admissible-rotations}
Both theorems hold at every frequency in the set $\mathcal S$
defined below. The maps may depend on the frequency.
Write ${p}_i/q_i=[0;a_1,\ldots,a_i]$, with $q_{-1}=0$, $q_0=1$.
Fix integers $B,k_0\ge2$. At the positions $m=2k^2$, $k\ge k_0$, prescribe
\[
 a_m=\begin{cases}2,&q_{m-2}\text{ odd},\\1,&q_{m-2}\text{ even},\end{cases}
 \qquad a_{m+1}=2(k+2)^4,\qquad a_{m+2}=2.
\]
Choose all remaining partial quotients freely in $\{1,\ldots,B\}$
and take the union over $B,k_0$. This defines $\mathcal S$.
It has Lebesgue measure zero, but
\[
 \dim_H(\mathcal S\cap J)=1
 \qquad\text{in every nonempty open interval }J\subset(0,1).
\]
Every member has irrationality exponent two (Roth type in the
Diophantine sense) and unbounded partial quotients.  There are also admissible Liouville frequencies outside
this set.
%

\end{remark}

\vspace{1ex}\paragraph{\bf Method of construction.}
By the Herman--Mather formula (Proposition~\ref{in:herman}), it
suffices to construct an increasing bi-Lipschitz circle lift $g$
with the required dynamics and singular set such that
\[
 \phi=g+g^{-1}-2\Id
\]
has the prescribed regularity. The graph function is $u=g-\Id$,
so $N(u)=N(g)$.
\begin{itemize}[leftmargin=*,itemsep=.5em]
\item \emph{Ordered intervals and the full circle map.}
Following Herman's interval construction
\cite[Chap.~X, \S\S3.8--3.11]{Her79}, we replace two orbits of
the rigid rotation by intervals in their circular order.
A positive density parametrizes each interval, and $g$ sends a
point to the point with the same parameter in the next interval.
The other fibres are singletons. Summability of the lengths gives
a circle homeomorphism and a collapsing semiconjugacy to the
rotation. This step constructs the complete invariant graph,
including the wandering intervals, rather than only its minimal set.

\item \emph{Cancellation through a second-order recurrence.}
Near the chosen singular parameters, the densities combine two
linearly independent solutions of the same
\emph{second-order linear recurrence}~\eqref{pf:linear-recurrence}.
The combination gives different limits for forward secant slopes,
while the recurrence fixes the sum of the forward and backward
slopes. This separates nondifferentiability of $g$ from regularity
of $\phi$. The use of a smooth symmetric sum is shared with
\cite{Ar11,Ar14,AF}; the construction here prescribes a singular
parameter set and carries this cancellation through all the
intervals. Remark~\ref{rem:construction-sources} gives the precise
relation to the interval constructions of Herman and Arnaud.

\item \emph{Two rotation orbits and global derivative estimates.}
The arithmetic part concerns only the rigid rotation.
A step function takes opposite values on the orbits of zero and
one half, and its exponentiated sums give reciprocal sequences.
Selected returns provide the indices at which we truncate and
interpolate the densities. We use different indices on the two
orbits so that both families have summable lengths while the
derivatives of the perturbation remain controlled in the circle
coordinate. This is the reason for the two-orbit choice in
Figure~\ref{fig:nodes}. After constructing $g$, we first extend
the density quotients to a global coefficient $Q$ and establish
its regularity. Absolute continuity then identifies $Q$ with
$\mathrm D\phi$.  One choice of radii
valid for all derivative orders gives Theorem~\ref{thm:smooth}.
\end{itemize}


\vspace{1ex}\paragraph{\bf Questions of Sinai and Mather and remaining cases.}
Sinai \cite[\S I, item~4]{Sinai} asked for examples of nonsmooth
invariant curves of twist maps in his discussion of regions with
different dynamics. Theorem~\ref{thm:smooth} gives such a curve
for a smooth twist map and identifies its nondifferentiability set.
Mather's questions appear in
\cite[\S3.1, Problems~3.1.1 and~3.1.3, pp.~181--182]{MY}.
Problem~3.1.1 asks for a rotational invariant curve that is not
$C^1$, for each regularity of the  area-preserving twist
map; here rotational means homotopically nontrivial.
Problem~3.1.3 asks for an irrational rotational curve on which the
dynamics is not topologically transitive, with the map of class
$C^3$, of higher regularity, or real analytic.
Theorem~\ref{thm:smooth} gives both properties in the $C^\infty$
class. Indeed, an orbit visits a given wandering interval at most
once and hence cannot be dense in the circle.

Two cases remain outside this construction:
\begin{itemize}[leftmargin=*,itemsep=.3em]
\item a $C^\infty$ area-preserving twist map with a non-$C^1$
invariant graph of irrational rotation number and minimal
restricted dynamics;
\item a real-analytic area-preserving twist map with a non-$C^1$
invariant graph of irrational rotation number, with either minimal
or nonminimal restricted dynamics.
\end{itemize}
Fayad and Krikorian \cite[\S4.4.1, Question~26]{FK} ask for a
non-$C^1$ invariant curve with minimal restricted dynamics for
twist maps of class $C^r$, $r\ge2$, and in the real-analytic
category. A $C^\infty$ example with minimal dynamics would answer
their finite-regularity question simultaneously for every $r$.
The interval insertion used here produces nonminimal dynamics,
so it does not address that case. The analytic problem requires
a different cancellation mechanism: our perturbation is locally
affine near the prescribed singular orbits, whereas a real-analytic
periodic function that is affine on a nonempty interval must be
constant. The analytic dissipative examples in \cite{KMW26}
therefore do not settle the area-preserving question.

\vspace{1ex}\paragraph{\bf Organization of this work.}
Section~\ref{sec:arith} recalls the Herman--Mather formula and
establishes the arithmetic estimates used in the construction.
Section~\ref{sec:profiles} constructs densities on the orbit of zero
through an interpolation formula for a second-order recurrence and
estimates the resulting perturbation. Section~\ref{sec:reciprocal}
constructs smooth densities on the orbit of one half and proves the
corresponding estimates. Section~\ref{sec:scales} chooses the singular
set and cutoff radii to obtain summability and the required regularity,
including the endpoint $r=1$ and the smooth case.
Section~\ref{sec:realization} realizes these densities as a circle
homeomorphism and establishes the global regularity of the perturbation.
Section~\ref{sec:singular} identifies the  nondifferentiability
set, proves strict differentiability on the minimal set, and completes
the proofs of both theorems. Appendix~\ref{app:maps} illustrates the
maps and coordinates. Appendix~\ref{app:frequencies} extends the
construction to further rotation numbers and studies their arithmetic
properties. Finally, Appendix~\ref{app:sharp} proves the dimension
bound under the assumption that the perturbation is locally affine
near every iterate of each nondifferentiability point in a wandering
interval, establishes sharpness within this class, and determines the
exact regularity of the perturbations in Theorem~\ref{thm:finite}.

\vspace{1ex}\paragraph{\bf Acknowledgments.}
This work was partially supported by the National Natural Science Foundation of China (Grant Nos.~12671217 and 12122109).

\section{Preliminaries and arithmetic estimates}\label{sec:arith}

\subsection{Maps, regularity and the Herman--Mather formula}
Following Herman \cite[Chaps.~I--II]{Her79}, write
\[
 \Id:\R\to\R,\quad \Id(x):=x,
 \qquad R_\alpha:\R\to\R,\quad R_\alpha(x):=x+\alpha.
\]
The latter also denotes the induced rotation on $\T$.
An increasing degree-one lift is an increasing homeomorphism
$g:\R\to\R$ satisfying $g(x+1)=g(x)+1$; we write $g\in D^0(\T)$.
Its real rotation number is
\[
 \rho(g):=\lim_{n\to\infty}\frac{g^n(x)-x}{n},
\]
which is independent of $x$. The rotation number of the circle map
is this number modulo one.

A semiconjugacy in this paper is a continuous nondecreasing map
$h_g:\R\to\R$ such that
\[
 h_g(x+1)=h_g(x)+1,
 \qquad h_g(g(x))=h_g(x)+\alpha.
\]
It induces a surjection of circles. Section~\ref{sec:realization}
constructs it by specifying every fibre. The subscript distinguishes
it from a positive increment $h$ in estimates. We write $\mathrm D$
for differentiation in the circle coordinate and use primes on
interval parameterizations for differentiation in their parameter.
A wandering interval has pairwise disjoint iterates.
For irrational rotation number, if there is a wandering interval, the
unique minimal set is a Cantor set; its complementary intervals are
wandering \cite[Chap.~II, \S\S7.2--7.3]{Her79}.

Strict differentiability of $u$ at $x$ means that $(u(z)-u(y))/(z-y)$ has a limit as $y,z\to x$, $y\ne z$. Ordinary differentiability only requires the quotients based at $x$ to have a limit.

A function $v$ is ${\beta}$-H\"older on an interval $J$ if
\[
 [v]_{{\beta},J}:=\sup_{\substack{x,y\in J\\x\ne y}}
 \frac{|v(x)-v(y)|}{|x-y|^{\beta}}<\infty.
\]
The $C^{m,{\beta}}$ norm used below is the sum of the supremum norms of derivatives through order $m$ and $[v^{(m)}]_{{\beta},J}$. For integer order we omit the seminorm. {On $\T$, the H\"older seminorm is defined
using the shorter arc distance $\dist(x,y)$; interval seminorms
use the real length coordinate.} A Lipschitz function on an interval is absolutely continuous:
the sum of its absolute increments over disjoint subintervals is
bounded by its Lipschitz constant times their total length. For absolutely continuous functions, the derivative exists almost everywhere, is integrable, and its integral recovers every increment.

A null set means a set of Lebesgue measure zero. We always specify the
exponent when referring to vanishing Hausdorff measure.
For a bounded set $B\subset\R$, $\mathcal H^s(B)$ denotes the limit as $\eta\downarrow0$ of the infimum of $\sum_i |J_i|^s$ over countable interval covers with $|J_i|\le\eta$ (we use interval lengths without an additional constant). The Hausdorff dimension is $\inf\{s:\mathcal H^s(B)=0\}$. If a probability measure supported on $B$ assigns each sufficiently short interval $J$ mass at most $C|J|^t$, every such cover gives $1\le C\sum_i|J_i|^t$, hence $\dim_H B\ge t$. An $L$-Lipschitz map multiplies covering lengths by at most $L$; consequently bi-Lipschitz changes preserve dimension and the vanishing of a Hausdorff measure. The dimension of a countable union is the supremum of the dimensions of its members, by countable subadditivity of the covering measures.

We use $A\asymp B$ for positive quantities if $cB\le A\le CB$ with positive constants independent of the varying index or scale. In an estimate $O_m(\cdot)$ the constant may depend on the indicated, fixed derivative order. Constants may also depend on fixed functions chosen in the proof.

Appendix~\ref{app:maps} lists the maps used in the construction
and illustrates the circles and their fibres
(Figure~\ref{appfig:circles}), the density coordinates and interval maps
(Figure~\ref{appfig:densities}), and the singular set and perturbation
(Figure~\ref{appfig:regularity}).

The Herman--Mather formula originates in the work of Herman (1983) and Mather (1984). Herman \cite[Chap.~II, \S2]{Her} constructs a twist map from an increasing circle lift. Mather \cite[\S\S5--6, p.~306, (3)--(4)]{Mat84} obtains the corresponding orbit identity and its almost-everywhere differentiated form from an invariant circle of the sine standard map. We use the normalization of \eqref{in:map}; see also \cite[p.~493, (1)--(2)]{AF}.

\begin{proposition}[Herman--Mather formula]\label{in:herman}
Let $g$ be an increasing bi-Lipschitz lift and put
\[
 u=g-\Id,\qquad \phi=g+g^{-1}-2\Id .
\]
Then $u$ and $\phi$ are periodic, $\int_0^1\phi(x)\,dx=0$, and
$F_\phi(x,u(x))=(g(x),u(g(x)))$. Conversely, every invariant graph for \eqref{in:map} with induced lift $g=\Id+u$ satisfies this identity for $\phi$.
\end{proposition}

When $\phi\in C^1$, as in both theorems, the twist derivative is $\partial_y(x+y)=1$.
Define $\Phi:\R\to\R$ by $\Phi(x):=\int_0^x\phi(s)\,ds$.
The mean-zero condition makes this primitive periodic, and
\[
 F_\phi^*(y\,dx)-y\,dx
 = d\bigl(y^2/2+\Phi(x+y)\bigr).
\]
In particular $F_\phi$ is exact and area-preserving; it is the composition
of the shears $(x,y)\mapsto(x+y,y)$ and $(x,y)\mapsto(x,y+\phi(x))$.

\subsection{Continued fractions and orbit products}

We first establish two estimates for the rigid rotation
$R_\alpha$. At convergent denominators, the Denjoy--Koksma inequality
bounds the multiplier products above and below, uniformly in the
starting point. Between selected consecutive returns, the partial
quotients in \eqref{in:frequency} give a large product along the orbit
of zero. These estimates will be used to construct $g$.

Let $p_i/q_i=[0;a_1,\ldots,a_i]$ denote the continued-fraction convergents of $\alpha$, with
\[
\begin{gathered}
(p_{-1},q_{-1})=(1,0),\qquad (p_0,q_0)=(0,1),\\
p_{i+1}=a_{i+1}p_i+p_{i-1},\qquad
q_{i+1}=a_{i+1}q_i+q_{i-1}.
\end{gathered}
\]
Here the subscripted $q_i$ is a denominator; the function $q$ introduced below is a recurrence coefficient.

\begin{lemma}[Continued fractions and denominator growth]\label{ar:denominators}
For every irrational $\alpha\in(0,1)$, consecutive denominators are coprime, and $p_i,q_i$ are coprime. For $\tau=[a_{i+1};a_{i+2},\ldots]$,
\begin{equation}\label{ar:convergents}
0<(-1)^i(q_i\alpha-p_i)
=(\tau q_i+q_{i-1})^{-1}<q_{i+1}^{-1}.
\end{equation}
For the frequency in \eqref{in:frequency}, $q_{2n}$ is odd, $q_{2n+1}$ is even, and
\begin{equation}\label{ar:denominator-growth}
\log q_{2n+2}=O(n\log(n+2))\qquad(n\ge1).
\end{equation}
The implicit constant in \eqref{ar:denominator-growth} is independent of $n$.
\end{lemma}
\begin{proof}
The coprimality and approximation identities are the standard continued-fraction formulas; see \cite[Chap.~V, \S7.2]{Her79}.
For \eqref{in:frequency}, reduction of the denominator recurrence modulo two gives $q_{i+1}\equiv q_{i-1}\pmod2$, proving the parity assertion from $q_{-1}=0,q_0=1$.
Since $q_{i+1}\le(a_{i+1}+1)q_i$, induction yields
$q_{2n+2}\le\prod_{i=1}^{2n+2}(a_i+1)$.
Taking logarithms and using $a_{2l+1}=2(l+2)^4$ and $a_{2l+2}=2$ proves \eqref{ar:denominator-growth}.
\end{proof}

Define the right-continuous step function $\sigma:\T\to\{-1,1\}$ by
\[
\sigma(\theta):=
\begin{cases}
1,&0\le\theta<1/2,\\
-1,&1/2\le\theta<1,
\end{cases}
\]
using the representative of $\theta$ in $[0,1)$. Define the multiplier $\lambda:\T\to(0,\infty)$ and the recurrence coefficient $q:\T\to\R$ by
\begin{equation}\label{ar:fields}
\begin{aligned}
 \lambda(\theta)&:=e^{\sigma(\theta)},&
 q(\theta)&:=\lambda(\theta)+\lambda(\theta-\alpha)^{-1}-2,\\
 U_0&:=1,& U_{j+1}&:=\lambda(j\alpha)U_j,\qquad S_j:=\log U_j
 \quad (j\in\Z).
\end{aligned}
\end{equation}
The recurrence defines $U_j$ in both time directions. Its logarithm is the signed Birkhoff sum
\begin{equation}\label{ar:signed-sums}
S_j=
\begin{cases}
\displaystyle\sum_{i=0}^{j-1}\sigma(i\alpha),&j>0,\\[2pt]
0,&j=0,\\[2pt]
\displaystyle-\sum_{i=j}^{-1}\sigma(i\alpha),&j<0.
\end{cases}
\end{equation}
In particular, $S_{j+1}-S_j=\sigma(j\alpha)$. The symbols $S_j$ and $S$ have different meanings: the later $S$ is a total interval length, not a limit or a sum of the $S_j$.
The two sequences of rotation points are
\[
\theta_j:=R_\alpha^j(0)=j\alpha\pmod1,\qquad
\widehat\theta_j:=R_\alpha^j(1/2)=1/2+j\alpha\pmod1.
\]
For a step function $v$, write $v(\theta^\pm)$ for its right and left limits, respectively; equivalently $v^\pm(\theta)=v(\theta^\pm)$. For our right-continuous functions, $v^+(\theta)=v(\theta)$.

\begin{lemma}[Multiplier products and the two rotation orbits]\label{ar:products}
For every irrational $\alpha$, the orbits of $0$ and $1/2$ are disjoint,
\[
\lambda(\theta+1/2)=\lambda(\theta)^{-1},\qquad
 e^{-1}\le\frac{U_{j+1}}{U_j}\le e,
\qquad 2e^{-1}\le 2+q(\theta)\le2e.
\]
The multiplier product based at $1/2$, normalized to one at time zero, is $U_j^{-1}$ at time $j\in\Z$. The function $q$ is locally constant off
\[
\{0,\alpha,1/2,1/2+\alpha\}\pmod1.
\]
These assertions also hold with a fixed choice of one-sided values wherever applicable.
\end{lemma}
\begin{proof}
An intersection of the two orbits would give $(i-j)\alpha=1/2\pmod1$, contradicting irrationality.
The identity $\sigma(\theta+1/2)=-\sigma(\theta)$ gives reciprocity of $\lambda$ and hence of the products in both time directions. The bounds follow from $\lambda\in\{e^{-1},e\}$. The only jumps of $\lambda$ are at $0$ and $1/2$, and translation by $\alpha$ gives the displayed possible jump set of $q$.
\end{proof}

\begin{lemma}[Nonnegative sums at the origin]\label{ar:lower}
{If the odd-indexed partial quotients of $\alpha$ are all even, then
$S_j\ge0$ and $U_j\ge1$ for every $j\in\Z$. This applies in
particular to \eqref{in:frequency}.}
\end{lemma}
\begin{proof}
{By the half-interval criterion of Boshernitzan and Ralston
\cite[Theorem~1, with $m=2$]{BR},
$\sum_{k=1}^n\sigma(k\alpha)\ge0$ for every $n\ge1$.
For $j>0$, the convention $\sigma(0)=1$ and the identity
$\sigma(-k\alpha)=-\sigma(k\alpha)$ for $k\ge1$ give
\[
 S_j=1+\sum_{k=1}^{j-1}\sigma(k\alpha),\qquad
 S_{-j}=\sum_{k=1}^j\sigma(k\alpha)=S_{j+1}-1 .
\]
These identities prove the assertion for both positive and
negative indices.}
\end{proof}

\begin{lemma}[Denjoy--Koksma at convergent denominators]\label{ar:grid}
For every convergent denominator $q_i$, every $z\in\T$, and either
one-sided representative $\sigma^\pm$ of $\sigma$,
\[
 \left|\sum_{l=0}^{q_i-1}\sigma^\pm(z+l{\alpha})\right|\le4 .
\]
Consequently, if $P_0=1$ and
$P_{j+1}/P_j=\lambda((z+j{\alpha})^\pm)$ with a fixed side, then
$e^{-4}\le P_{\pm q_i}\le e^4$. Both assertions hold for every
irrational ${\alpha}$, independently of \eqref{in:frequency}.
\end{lemma}
\begin{proof}
{Apply the Denjoy--Koksma inequality
\cite[Chap.~VI, Theorem~3.1, p.~73]{Her79} to the rotation $R_\alpha$.
Both one-sided representatives of $\sigma$ have zero integral and total
variation four. This gives the bound for every starting point. The
logarithm of $P_{q_i}$ is the displayed sum; the logarithm of
$P_{-q_i}$ is its negative with starting point $z-q_i\alpha$.}
\end{proof}

\begin{lemma}[Maximal sums between returns]\label{ar:peaks}
Assume that $\alpha$ is given by \eqref{in:frequency}. Put $M_n:=\max_{0\le j<q_{2n+2}}S_j$. For constants $c,C>0$
independent of $n$ and all sufficiently large $n$,
\begin{equation}\label{ar:increments}
 cn^4\le M_n-M_{n-1}\le Cn^4,\qquad cn^5\le M_n\le Cn^5 .
\end{equation}
Choose an integer $k_n\in[0,q_{2n+2})$ with $S_{k_n}=M_n$. Every such choice satisfies $q_{2n}<k_n<q_{2n+2}$ for all sufficiently large $n$.
In particular $k_n<q_{2n+2}<k_{n+1}$, and
\begin{equation}\label{ar:signedpeaks}
 \log {U}_{k_n}=M_n,\qquad |\log {U}_{-k_n}-M_n|\le2 .
\end{equation}
Moreover,
\begin{equation}\label{ar:separation}
\min\{k_n-q_{2n},\ q_{2n+2}-k_n\}\ge M_n-4.
\end{equation}
Thus $k_n$ and its immediate neighbours lie strictly between consecutive selected returns when $n$ is sufficiently large. The reflected assertion holds at $-k_n$.
\end{lemma}
\begin{proof}
To prove the lower bound, we increase the previous maximum by
order $a_{2n+1}$. We do this by placing a prefix that attains the
previous maximum after consecutive blocks whose sums are all $+1$.
Write $h=q_{2n}$, $H=q_{2n+1}$, $a=a_{2n+1}$ and
$\varepsilon=h{\alpha}-{p}_{2n}$. Then
\[
 H=ah+q_{2n-1},\qquad q_{2n+2}=2H+h,\qquad
 0<\varepsilon<1/H .
\]
At time $lh+i$ ($0\le i<h$) the rotation point is, modulo one,
\[
 \frac{i{p}_{2n}}h+(l+i/h)\varepsilon .
\]
For $0\le l\le L:=\lfloor a/8\rfloor$ the displacement from the
grid is less than $(L+1)/H<1/(2h)$ for large $n$. The grid point
zero moves to its right and keeps its sign. Every other grid point
is at least $1/h$ from zero, and $h$ is odd, so every grid point
is at least $1/(2h)$ from $1/2$. All these blocks have the same
ordered sequence of signs, with $(h+1)/2$ plus and $(h-1)/2$
minus signs. If $S_i=M_{n-1}$ for $0\le i<h$, it follows that
$S_{Lh+i}=L+S_i$; this index is less than $q_{2n+2}$. This proves
the required lower bound.

For the upper bound write $j=lh+i<2H+h$ with $0\le i<h$. Then
$l\le2a+2$. Lemma~\ref{ar:grid} bounds the sum on each complete
block by $4$, regardless of its starting point. The remaining prefix
of length $i$ differs from the length-$i$ prefix starting at zero
by a translation $l\varepsilon\le3/h$. A sign can change only when
the corresponding point of the original $h$-point grid is within
$4/h$ of zero or $1/2$; the extra $1/h$ allows for the displacement
at times $0,\ldots,h-1$. There are at most $18$ such grid points.
Thus $S_j\le4(2a+2)+M_{n-1}+36$. As $a=2(n+2)^4$, summing the
bounds on $M_n-M_{n-1}$ proves \eqref{ar:increments}.

No index less than $h$ now attains $M_n$, and $|S_h|\le4$ by
Lemma~\ref{ar:grid}. Thus $k_n>h$ for large $n$. The identity
$S_{-j}=S_{j+1}-1$ and $|S_{j+1}-S_j|=1$ give
\eqref{ar:signedpeaks}. Finally, \eqref{ar:separation} follows because $|S_{q_{2n}}|,|S_{q_{2n+2}}|\le4$ and each step changes the sum $S_j$
by only one. Thus each maximizing time and its immediate neighbours lie in
the same interval between returns, for both positive and negative
indices.
\end{proof}

\begin{lemma}[Uniform bounds on the interpolation blocks]\label{ar:blocks}
Assume \eqref{in:frequency}, and write
\[
J_n^+:=[k_n,k_{n+1}]\cap\Z,
\qquad J_n^-:=-J_n^+.
\]
There is a constant $C$ independent of $n$ such that, for all sufficiently large $n$,
\[
\begin{gathered}
\log\bigl(2+\max_{j\in J_n^+\cup J_n^-}|j|\bigr)
+\log\bigl(\#J_n^++\#J_n^-\bigr)
\le Cn\log(n+2),\\
1\le U_j\le e^{M_{n+1}+2}\qquad(j\in J_n^+\cup J_n^-).
\end{gathered}
\]
The return $q_{2n+2}$ and its immediate neighbours lie in $J_n^+$; their negatives lie in $J_n^-$. At these six indices,
\[
e^{-5}\le U_j\le e^5.
\]
At $j=\pm k_n$ one has $e^{M_n-2}\le U_j\le e^{M_n+2}$.
\end{lemma}
\begin{proof}
The number and absolute value of the indices in either block are bounded by $2q_{2n+4}+1$, so Lemma~\ref{ar:denominators} gives the first estimate.
For $j\ge0$ in the positive block, $S_j\le M_{n+1}$ by definition. For the reflected block use $|S_{-j}-S_j|\le2$, proved in Lemma~\ref{ar:peaks}; Lemma~\ref{ar:lower} supplies the lower bound.
The separation estimate \eqref{ar:separation} at $k_n$ and $k_{n+1}$ places the return and its neighbours in the indicated blocks. Lemma~\ref{ar:grid} bounds the products at the return by $e^{\pm4}$, and one multiplication or division by $\lambda$ gives $e^{\pm5}$ at its neighbours. The last assertion is \eqref{ar:signedpeaks}.
\end{proof}

The selected return times and maximizing times alternate. At the return times
$\pm q_{2n}$ the two reciprocal products are both bounded above and
below; at $\pm k_n$ one is large and the other is small.
Sections~\ref{sec:profiles} and~\ref{sec:reciprocal} use the latter
times to prescribe densities on the first orbit and the former times
on the second. The uniformity in the starting point in
Lemma~\ref{ar:grid} will also rule out additional singular points
on the minimal set in Section~\ref{sec:singular}.

\section{Densities on the orbit of zero}\label{sec:profiles}

To obtain a smooth perturbation with a nonsmooth induced map, we
construct densities for which the quotient involving three consecutive
densities is smooth, although
the quotient of two consecutive densities need not have a limit.
We first construct the densities on the orbit of zero, using the
sequence ${U}_j$ and the integers $k_n$ from Section~\ref{sec:arith}.

The integer $j$ is the iterate number: an interval will replace the
point $j{\alpha}$ on the rotation circle. A parameter $t\in[0,1]$
specifies a point in that interval. Under one iterate $j$ increases
by one and $t$ stays fixed. Choose a nonempty compact set
$E\subset[1/3,2/3]$ of Lebesgue measure zero; its elements are
the parameters at which differentiability will fail.
For each $j$ we will construct a Borel density
$\rho_j:[0,1]\to(0,\infty)$, bounded above and away from zero for that
fixed $j$. The integral of the density determines interval lengths. Define
\begin{equation}\label{pf:length-coordinate}
 x_j:[0,1]\longrightarrow[0,w_j],\qquad
 x_j(t):=\int_0^t\rho_j(v)\,dv,
 \qquad w_j:=\int_0^1\rho_j(v)\,dv.
\end{equation}
This is an unnormalized length coordinate measured from the left
endpoint. The common parameter $t$ identifies points along an orbit;
it is neither the iterate number nor, in general, the fraction of an
interval's length.
Thus two parameters $t_1<t_2$ determine a subinterval of length
$x_j(t_2)-x_j(t_1)$. There is no probability normalization.
When the interval maps send $t$ to the same $t$ at the next time,
their secant slopes are
\[
 \frac{x_{j+1}(t_2)-x_{j+1}(t_1)}{x_j(t_2)-x_j(t_1)}
 =\frac{\int_{t_1}^{t_2}\rho_{j+1}(v)\,dv}
        {\int_{t_1}^{t_2}\rho_j(v)\,dv}.
\]
To insert the intervals into a circle we must have
$\sum_j w_j<\infty$. To prescribe the regularity of the perturbation
we must also estimate
\begin{equation}\label{pf:Q}
 Q_j:[0,1]\to\R,\qquad Q_j(t):=\frac{\rho_{j+1}(t)+\rho_{j-1}(t)}{\rho_j(t)}-2 .
\end{equation}
Indeed, Section~\ref{sec:realization} identifies $Q_j$ with the
derivative of the perturbation in the Herman--Mather formula.
Although $Q_j$ is written as a function of $t$, its higher derivatives
must be estimated with respect to the length coordinate $x_j$.
The integral coordinates and their role in defining $g$ are shown
in Figure~\ref{appfig:densities} of Appendix~\ref{app:maps}.

\subsection{An interpolation formula for a second-order recurrence}

For an arbitrary positive sequence $A_j$, put
\[
 a_j^A:=\frac{A_{j+1}+A_{j-1}}{A_j}-2,
 \qquad d_j:=\frac1{A_jA_{j+1}}.
\]
We use the \emph{second-order linear recurrence equation}
\begin{equation}\label{pf:linear-recurrence}
 z_{j+1}+z_{j-1}-(2+a_j^A)z_j=0,
 \qquad j\in\mathbb Z.
\end{equation}
The coefficients $a_j^A$ are prescribed, and the sequence $z_j$ is the
unknown. Two consecutive values determine a unique solution on all of
$\mathbb Z$: the equation gives
$z_{j+1}=(2+a_j^A)z_j-z_{j-1}$ in the forward direction and
$z_{j-1}=(2+a_j^A)z_j-z_{j+1}$ in the backward direction.

For $z_j=A_j y_j$, multiplication and subtraction give
\begin{equation}\label{pf:flux}
 A_j\bigl(z_{j+1}+z_{j-1}-(2+a_j^A)z_j\bigr)
 =\frac{y_{j+1}-y_j}{d_j}-\frac{y_j-y_{j-1}}{d_{j-1}}.
\end{equation}
By \eqref{pf:flux}, equation~\eqref{pf:linear-recurrence} implies that
$(y_{i+1}-y_i)/d_i$ has the same value for all $l\le i<h$.
Summing these increments from $l$ to $h$ shows that this common
value is $(y_h-y_l)/(\sum_{l\le i<h}d_i)$. In particular, given the
two endpoint values, the unique solution inside the block is
\begin{equation}\label{pf:interpolate}
 y_j=(1-\tau)y_l+\tau y_h,\qquad
 \tau=\frac{\sum_{l\le i<j}d_i}{\sum_{l\le i<h}d_i}.
\end{equation}
Since $0\le\tau\le1$, the interpolated value is positive if both prescribed values are positive.

Define the signed sum $L_j$ by $L_0=0$,
$L_{j+1}-L_j=({U}_j{U}_{j+1})^{-1}$ for all $j\in\mathbb Z$, and put
\begin{equation}\label{pf:AT}
 A_j={U}_j(1+|L_j|),\qquad
 W=\sum_{i\in\mathbb Z}\frac1{A_iA_{i+1}},\qquad
 T_j=\frac1W\sum_{i<j}\frac1{A_iA_{i+1}} .
\end{equation}
Applying \eqref{pf:flux} first with reference sequence ${U}$ shows
that $a_j^A=q(j{\alpha})$ for $j\ne0$:
$1+|L_j|$ is affine in $L_j$ on each of the two half-lines.
At $j=0$ the slopes of $1+|L_j|$ are $-1$ to the left and $1$ to the
right. Since $A_0={U}_0=1$, \eqref{pf:flux} gives
$a_0^A=q(0)+2$. With this change at zero, the recurrence has two linearly independent
positive solutions, as shown below.
Also $0<W\le2$. Indeed, for $i\ge0$ the summand of $W$ equals
$(1+L_i)^{-1}-(1+L_{i+1})^{-1}$; for $i<0$ it equals
$(1-L_{i+1})^{-1}-(1-L_i)^{-1}$. Since every summand is positive and both tails contain
infinitely many terms, $0<T_j<T_{j+1}<1$. Moreover,
\begin{equation}\label{pf:T-increment}
 T_{j+1}-T_j=\frac{d_j}{W}
              =\frac{1}{W A_jA_{j+1}}>0.
\end{equation}
For $z_j=A_j$, the corresponding sequence in \eqref{pf:flux} is
$y_j=1$, so both weighted increments are zero. For $z_j=A_jT_j$,
that sequence is $y_j=T_j$, and both weighted increments equal
$1/W$ by \eqref{pf:T-increment}. Thus both sequences solve
\eqref{pf:linear-recurrence} at every index, including zero.

Their values at two consecutive indices satisfy
\begin{equation}\label{pf:wronskian}
 \det\begin{pmatrix}
 A_j&A_jT_j\\
 A_{j+1}&A_{j+1}T_{j+1}
 \end{pmatrix}
 =A_jA_{j+1}(T_{j+1}-T_j)=\frac1W>0.
\end{equation}
 Hence the solutions
are linearly independent. For any solution $z$, the nonzero determinant gives unique
$c_0,c_1$ matching $z_j,z_{j+1}$.
The uniqueness from two consecutive values then gives
$z_i=c_0A_i+c_1A_iT_i$ for every $i\in\mathbb Z$.
In particular, $A_j(1+cT_j)$ solves the recurrence for every fixed
real number $c$. The bound ${U}_i\ge1$ gives $|L_j|\le |j|$.
We next explain why the coefficient is changed at zero.
At the indices $i=\pm q_l$, Lemma~\ref{ar:grid} bounds ${U}_i$
from above independently of $l$, and ${U}_{i+1}/{U}_i\le e$.
Thus $({U}_i {U}_{i+1})^{-1}$ is bounded below by a positive constant
on each of these two infinite subsequences. It follows that
$L_j\to+\infty$ as $j\to+\infty$ and $L_j\to-\infty$ as $j\to-\infty$.
For the unmodified coefficient $q(j{\alpha})$, \eqref{pf:flux} would
give $z_j/{U}_j=a+bL_j$ on the whole integer line, so an everywhere
positive solution would have to be a multiple of ${U}$. For the modified
coefficient the telescoping sums above give instead
\[
 W=2,\qquad
 T_j=\tfrac12\!\left(1+\frac{L_j}{1+|L_j|}\right).
\]
The bounded, strictly increasing $T_j$ therefore permits different
positive solutions along the same orbit.

The ratio $T_j$ of the second solution to the first is bounded. Since $0<T_j<1$, for every $|c|\le1/4$
their combination satisfies
\[
 \frac34A_j\le A_j(1+cT_j)\le\frac54A_j.
\]
We will allow $c$ to depend irregularly on $t$, using the same value
at every iterate. For each fixed $t$, the densities still satisfy the recurrence,
while their successive ratios may vary with $t$. This property depends
on the defining increments of $T_j$.

\subsection{The two endpoint values and the interior density}

For $\theta\in\mathbb T$, superscripts $-$ and $+$ denote the limits from the left and right in an increasing local lift to $\mathbb R$.
Define ${U}_0^-=1$, ${U}_{j+1}^-=\lambda((j{\alpha})^-){U}_j^-$;
the right-hand sequence is ${U}_j$. Their successive ratios differ
only at step $j=0$, because $j{\alpha}$ never equals $1/2$ and equals
$0$ only for $j=0$. At this step the right multiplier is $e$ and the
left multiplier is $e^{-1}$. It follows that ${U}_j^-/{U}_j=1$ for
$j\le0$ and ${U}_j^-/{U}_j=e^{-2}$ for $j\ge1$. The two sequences have
respective recurrence coefficients $q((j{\alpha})^-)$ and $q(j{\alpha})$;
these can differ only at $j=0,1$, since a coefficient involves the
multipliers at times $j$ and $j-1$.

Put $E_0:=E\cup\{0,1\}$ and define $D:[0,1]\to[0,\infty)$ by
\[
 D(t):=\begin{cases}
 (t-a)(b-t)/(b-a),&t\in(a,b),\ (a,b)\text{ a component of }[0,1]\setminus E_0,\\
 0,&t\in E_0.
 \end{cases}
\]
Then $\tfrac12\dist(t,E_0)\le D(t)\le\dist(t,E_0)$, $|D'|\le1$,
$D|D''|\le1$, and $D^{(i)}=0$ for $i\ge3$ on that component.
For a fixed $C_*>e^2$, define the bounded Borel function
$\Theta:[0,1]\to[-1,1]$ by
\[
 \Theta(t):=\begin{cases}
 \sin(\log\log(C_*/D(t))),&D(t)>0,\\
 0,&t\in E_0.
 \end{cases}
\]
Its values on the null set $E_0$ do not change any integral defining
an interval coordinate. Differentiating this formula gives
\begin{equation}\label{pf:theta-weight}
 D^i|\Theta^{(i)}|\le C_i\quad (i\ge0,\ D>0).
\end{equation}
The estimate holds at every derivative order, as follows. Write $L=\log(C_*/v)$.
Each $k$th derivative of $\sin(\log L)$ with respect to $L$
is a linear combination of $L^{-k}\sin(\log L)$ and $L^{-k}\cos(\log L)$.
Since $v\partial_v=-\partial_L$ and
$v^k\partial_v^k=(v\partial_v)(v\partial_v-1)\cdots(v\partial_v-k+1)$,
the $k$th derivative with respect to $v$ is bounded by $C_kv^{-k}$
for $0<v\le1$. Substitute $v=D(t)$ in the chain rule. A term with $k$
factors $D^{(i_1)},\ldots,D^{(i_k)}$, where $i_1+\cdots+i_k=i$,
is multiplied by at most $C D^{-k}$. The listed bounds on $D$
give $|D^{(i_l)}|\le C_{i_l} D^{1-i_l}$, proving \eqref{pf:theta-weight}.

Choose fixed nonnegative smooth functions
$\eta_0,\eta_1,\eta_c$ on $[0,1]$ with sum one,
$\eta_0=1$ near $0$, $\eta_1=1$ near $1$, and
$\eta_c=1$ on an open neighborhood of $[1/3,2/3]$.
For each iterate $j\in\Z$, define the uncut density
$P_j:[0,1]\to(0,\infty)$ and its quotient
$H_j:[0,1]\to(0,\infty)$ as functions of the parameter $t$ by
\begin{equation}\label{pf:raw}
 P_j(t)=\eta_0(t){U}_j^-+\eta_1(t){U}_j+
 \eta_c(t) A_j\bigl(1+T_j\Theta(t)/4\bigr),\qquad H_j=P_j/A_j .
\end{equation}
The factor multiplying $A_j$ lies between $3/4$ and $5/4$ wherever
$\eta_c=1$. Uniformly in $j$,
\begin{equation}\label{pf:H}
 c(1+|L_j|)^{-1}\le H_j(t)\le C,\qquad
 D^i|\partial_t^i H_j(t)|\le C_i\quad(D>0).
\end{equation}
In fact $H_j$ is a convex combination of ${U}_j^-/A_j$, ${U}_j/A_j$
and $1+T_j\Theta/4$. The first two numbers lie between
$e^{-2}/(1+|L_j|)$ and $1$, and the last lies in $[3/4,5/4]$.
Differentiating the fixed smooth weights and applying
\eqref{pf:theta-weight} gives the derivative bound.
For fixed $t$ the weights in \eqref{pf:raw} are independent of $j$;
the same holds for the coefficient $c=\Theta(t)/4$ of $T_j$.
Linearity of the recurrence shows that $P_j(t)$ solves it with
coefficient $a_j^A$ away from $0,1$. Thus $H_j$ satisfies
\eqref{pf:interpolate} inside every block whose interior avoids
these two indices.

\subsection{Cutoff times and estimates for the perturbation}

The uncut densities retain the oscillation near $E$ but their total
length is infinite: on the middle interval
$P_j\ge 3A_j/4$, and $A_{k_n}\to\infty$. Multiplying every $P_j$
by the same positive function of $t$ would preserve the ratios in
\eqref{pf:Q}, but the total length would remain infinite. We instead
reduce the density outside successively smaller neighborhoods
of $E\cup\{0,1\}$. We prescribe the reduced densities
at selected times and use interpolation to preserve the recurrence
between those times. Formula~\eqref{pf:nodeidentity} below shows that the resulting error
at time $j$ is multiplied by $A_j^{-2}$. We therefore choose times
when ${U}_j$ is large, so that this error is small.

Fix a nonincreasing smooth $Z:[0,\infty)\to(0,1]$ with $Z(v)=1$
for $v\le1$ and $Z(v)=1/v$ for $v\ge2$. Such a function is obtained
by using a smooth nondecreasing cutoff to pass from $1$ to $1/v$
on $[1,2]$. Its values are comparable to $\min(1,1/v)$ and
$v^i|Z^{(i)}(v)|\le C_i Z(v)$ for all $i$. For $R>0$ define $Z_R:[0,1]\to(0,1]$ by
$Z_R(t):=Z(D(t)/R)$. Let $R_n\downarrow0$ be positive parameter scales. The function
$Z_{R_n}:[0,1]\to(0,1]$ equals one where $D\le R_n$ and equals
$R_n/D$ where $D\ge2R_n$. Thus $R_n$ is a scale in $t$, not an
interval length in the circle coordinate. Decreasing $R_n$ reduces the
total length and preserves both the oscillation near $E$ and the
endpoint values. The choices for the three regularity regimes are
given in Section~\ref{sec:scales}. Their rate of decrease enters the estimates through
\begin{equation}\label{pf:radiusassumption}
 \Lambda_n:=\max\!\left\{\log\frac{R_{n-1}}{R_n},
                              \log\frac{R_n}{R_{n+1}}\right\}.
\end{equation}
We call $j=\pm k_n$, $n\ge n_0$, the outer cutoff indices and take
$n_0$ large enough for Lemma~\ref{ar:peaks}. Include every integer
strictly between $-k_{n_0}$ and $k_{n_0}$ as a central cutoff index.
In particular, the exceptional indices $0,1$ are prescribed indices;
neither lies in the interior of an interpolation block.
At $j=\pm k_n$ prescribe $\rho_j=P_j Z_{R_n}$; at each central
index set $\rho_j=P_j Z_1=P_j$. Between consecutive prescribed
indices define $y_j=\rho_j/A_j$ by \eqref{pf:interpolate}.
For each fixed $j$, the density is bounded and has a positive lower
bound, since it is obtained from two fixed positive radii and
\eqref{pf:H}. These bounds may depend on $j$. In the interior of
each block \eqref{pf:flux} gives $Q_j=a_j^A=q(j{\alpha})$.

On the block $k_n\le j\le k_{n+1}$, and on its reflection,
Lemma~\ref{ar:blocks} and $|L_j|\le|j|$ give
\[
 \log(1+|j|)=O(n\log(n+2)),\qquad
 A_j\le \exp(M_{n+1}+O(n\log(n+2))).
\]
At its two endpoints one has
$\log A_{\pm k_i}=M_i+O(i\log(i+2))$, $i=n,n+1$;
for negative times use $|\log {U}_{-k_i}-M_i|\le2$.
The interval of indices contains a return $\pm q_{2n+2}$ where
${U}_i,{U}_{i+1}$ are bounded. For all sufficiently large $n$ this
return lies strictly between its endpoints, since the values at the
two maximizing times tend to infinity. At that return $A_i,A_{i+1}$
are at most $\exp(O(n\log(n+2)))$. Hence for either block $[l,h]$
the number
\begin{equation}\label{pf:conductance}
 c(l,h):=\left(\sum_{l\le i<h}\frac1{A_iA_{i+1}}\right)^{-1}
 \quad\text{satisfies}\quad c(l,h)\le\exp(O(n\log(n+2))).
\end{equation}
We use these estimates to bound the total length of the intervals
in either block $J$. Indeed $Z_{R_{n+1}}\le Z_{R_n}$ since $Z$ is nonincreasing;
the prescribed value of $\rho_i/A_i$ is therefore bounded by $C Z_{R_n}$
at either end.
Interpolation is a convex combination at each $t$, so
$\rho_j(t)\le C A_j Z_{R_n}(t)$ throughout the block. Multiply
the bound for $\max_{j\in J}A_j$ by the number of indices, at most
$1+q_{2n+4}$. Integration gives
\begin{equation}\label{pf:mass}
 \sum_{j\in J}\int_0^1\rho_j(t)\,dt
 \le \exp(M_{n+1}+O(n\log(n+2)))\int_0^1 Z_{R_n}(t)\,dt .
\end{equation}
The exponential bound includes the number of indices in $J$.

\begin{figure}[htbp]
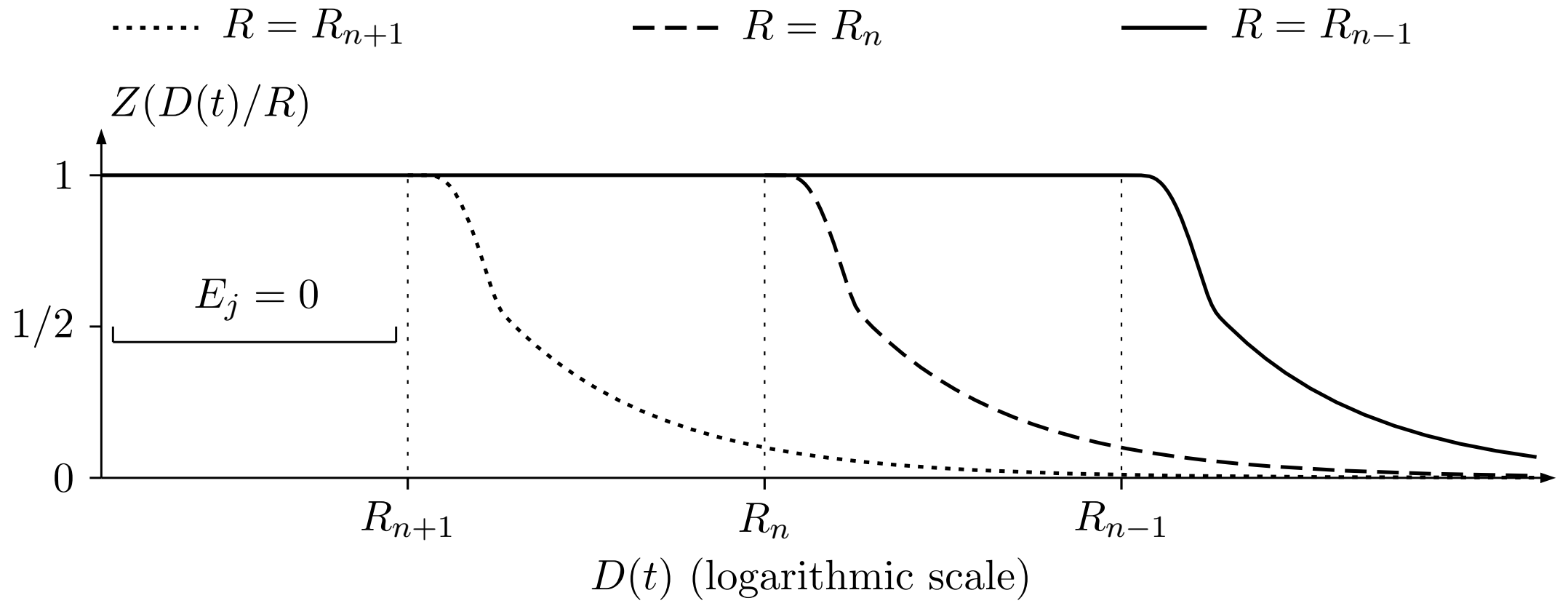

\centering
\begingroup\color{black}\normalsize\SinaiCutoffBW\endgroup
\caption{The three cutoffs in \eqref{pf:nodeidentity}. If $D(t)\le R_{n+1}$, all three equal one and the difference $Q_j-q(j{\alpha})$ vanishes. }\label{fig:cutoffs}
\end{figure}

We give the derivative estimate in detail. Fix $j=k_n$ or $-k_n$,
$n>n_0$, and let $l<j<h$ be the preceding and following cutoff indices. For each
of the three indices $i$ write $Z_i^*$ for its assigned cutoff.
Put $c_l=c(l,j)$ and $c_h=c(j,h)$.
We express the two differences adjacent to $j$ in terms of the
endpoint values. Indeed,
\eqref{pf:interpolate} gives
\[
 \frac{y_{j+1}-y_j}{d_j}=c_h(y_h-y_j),\qquad
 \frac{y_j-y_{j-1}}{d_{j-1}}=c_l(y_j-y_l).
\]
In \eqref{pf:flux}, setting $z_i=\rho_i$ then makes the left side
$A_j\rho_j(Q_j-a_j^A)=A_j^2y_j(Q_j-a_j^A)$. Solving for $Q_j-a_j^A$ gives the factor $A_j^{-2}$.
Subtracting the second slope from the first gives
\[
 Q_j-a_j^A=\frac{c_l(y_l-y_j)+c_h(y_h-y_j)}{A_j^2 y_j}.
\]
The same numerator with $H$ in place of $y$ vanishes by the recurrence
for $P$. Subtract it multiplied by $Z_j^*$ and use $y_i=H_i Z_i^*$.
One obtains the exact identity
\begin{equation}\label{pf:nodeidentity}
 E_j:=Q_j-q(j{\alpha})=
 A_j^{-2}\sum_{i\in\{l,h\}}c_i\frac{H_i}{H_j}
 \left(\frac{Z_i^*}{Z_j^*}-1\right).
\end{equation}
In Figure~\ref{fig:cutoffs}, all three cutoffs equal one up to the smallest radius $R_{n+1}$. Thus \eqref{pf:nodeidentity} gives $E_j=0$ there pointwise, independently of the oscillation of $\Theta$.

We estimate all terms in \eqref{pf:nodeidentity} on a component of
$[0,1]\setminus E_0$. For a smooth function $f$ there,
set $[f]_{m,D}(t)=\sum_{i=0}^mD(t)^i|f^{(i)}(t)|$.
The product rule gives $[fg]_{m,D}\le C_m[f]_{m,D}[g]_{m,D}$.
For $f=Z$ and $f=1/Z$ one has $v^k|f^{(k)}(v)|\le C_k f(v)$:
this follows from the formulas outside $[1,2]$ and from positivity
on the compact interval $1\le v\le2$. For $f(D/R)$ the same
chain rule as in \eqref{pf:theta-weight} gives
$[f(D/R)]_{m,D}\le C_m f(D/R)$, with $C_m$ independent of $R$.
For the reciprocal of $H_i$ differentiate $H_i(1/H_i)=1$:
\[
 \partial_t^a(1/H_i)=-H_i^{-1}\sum_{b=1}^a
    \binom{a}{b}(\partial_t^b H_i)\,\partial_t^{a-b}(1/H_i)
 \qquad(a\ge1).
\]
In view of \eqref{pf:H}, induction yields
$[1/H_i]_{m,D}\le C_m(1+|L_i|)^{m+1}$.
For the choice \eqref{in:frequency}, $|L_i|\le |i|$;
at the three indices this bound is $\exp(O_m(n\log(n+2)))$.
Finally $Z(D/R')/Z(D/R)\le C\max(1,R'/R)$ by the comparison with
$\min(1,R/D)$. Together these inequalities bound every derivative
of each ratio in \eqref{pf:nodeidentity}, and the product
$1/\rho_j=A_j^{-1}H_j^{-1}/Z_j^*$. The ratio of the larger radius to the smaller radius at adjacent
cutoff indices is at most $e^{\Lambda_n}$. Thus
\eqref{pf:nodeidentity} gives
\begin{equation}\label{pf:weighted}
 \begin{aligned}
 [E_j]_{m,D}&\le
    \exp\bigl(-2M_n+\Lambda_n+O_m(n\log(n+2))\bigr),\\
 [1/\rho_j]_{m,D}&\le
    \exp\bigl(O_m(n\log(n+2))\bigr)/\rho_j .
 \end{aligned}
\end{equation}
On $D\ge R_{n+1}$, use $D Z_{R_n}\ge c\min(D,R_n)$ to get
\begin{equation}\label{pf:scale}
 \rho_jD\ge\exp(M_n-O(n\log(n+2)))R_{n+1}.
\end{equation}

Define the interval error as the function
$e_j:[0,w_j]\to\R$ given by $e_j:=E_j\circ x_j^{-1}$; equivalently,
$e_j(x_j(t))=E_j(t)$.
On each component where $D>0$, $x_j'=\rho_j>0$ and
$\partial_x=\rho_j^{-1}\partial_t$. For instance, writing $b=1/\rho_j$,
one has $\partial_x^2 e_j=b b'E_j'+b^2E_j''$. More generally, every
term in the $m$th derivative is a numerical coefficient times
\[
 (\partial_t^{a_1} b)\cdots(\partial_t^{a_m} b)
                 (\partial_t^{a_0}E_j),\qquad
 a_0+\cdots+a_m=m,\quad a_i\ge0.
\]
This follows by induction from $\partial_x=b\partial_t$ and the
product rule. Multiplying the bounds
in \eqref{pf:weighted} for these factors gives an upper bound
$\exp(-2M_n+\Lambda_n+O_m(n\log(n+2)))(\rho_jD)^{-m}$ for the absolute value of the
sum. Applying \eqref{pf:scale} on the support of the error gives
\begin{equation}\label{pf:physical}
 \|e_{\pm k_n}^{(m)}\|_\infty
 \le \exp\bigl(-2M_n+\Lambda_n+O_m(n\log(n+2))\bigr)\,
             (e^{M_n}R_{n+1})^{-m}
 \qquad (m=0,1,2,\ldots).
\end{equation}
In particular, if $\Lambda_n=O(n^4)$, the exponent preceding
the last factor is $-2M_n+O_m(n^4)$. At order $m=0$ there is no
inverse length factor: $\Lambda_n=o(M_n)$ already suffices to make
$\|e_{\pm k_n}\|_\infty\to0$.
{For $j=\pm k_n$ and all sufficiently large $n$, the error
vanishes where $D\le R_{n+1}$. Consequently,
\[
 \supp e_j\subset x_j\bigl(\{D\ge R_{n+1}\}\bigr)
 \Subset (0,w_j)\setminus x_j(E).
\]
On a neighborhood of this compact set the density is smooth and
positive, so $x_j^{-1}$ is smooth there. Thus $e_j$, extended by
zero off its interval, is smooth. Estimate~\eqref{pf:physical}
does not require differentiating the inverse parameterization at
points of $E_0$.} A common rescaling of all interval coordinates
by a fixed positive constant changes only the constants at each order.
The remaining exceptional indices lie in
$\{0,1,\pm k_{n_0},\pm(k_{n_0}-1)\}$; the last pair consists of the
central indices adjacent to the innermost outer cutoffs.
For each fixed $j$, there is $\epsilon_j>0$ such that whenever
$D(t)<\epsilon_j$, the densities at $j-1,j,j+1$ equal
$P_{j-1},P_j,P_{j+1}$, respectively. Indeed only finitely many
cutoff radii enter those three interpolations, and interpolation of $H$
recovers $H$ itself. In an open neighborhood of $E$ the three
uncut densities then have the form
$A_i(1+T_i\Theta/4)$; their quotient \eqref{pf:Q} is the constant
$a_j^A$, even at a central index. Near $t=0$ and $t=1$ it equals
$q((j{\alpha})^-)$ and $q((j{\alpha})^+)$, respectively. Elsewhere
in a fixed interval the density is smooth and positive. It follows
that $Q_j$ as a function of $x_j$ is smooth for each fixed $j$
and has the stated one-sided endpoint values, with all positive-order
derivatives zero near each endpoint.

\section{Densities on the orbit of one half}\label{sec:reciprocal}

The multiplier $\lambda$ has jumps at both zero and $1/2$.
We must also match the left and right values on the orbit
of $1/2$. No singular point is prescribed inside its intervals.
For the second orbit $\widehat\theta_j=1/2+j\alpha$, we construct
positive smooth functions $\widehat\rho_j:[0,1]\to(0,\infty)$.
Define their length coordinates and lengths by
\begin{equation}\label{pf:hat-length-coordinate}
 \widehat x_j:[0,1]\to[0,\widehat w_j],\qquad
 \widehat x_j(t):=\int_0^t\widehat\rho_j(v)\,dv,
 \qquad \widehat w_j:=\int_0^1\widehat\rho_j(v)\,dv.
\end{equation}
The coefficient $\widehat Q_j:[0,1]\to\R$ is defined by
$\widehat Q_j:=(\widehat\rho_{j+1}+\widehat\rho_{j-1})/\widehat\rho_j-2$. The corresponding right-hand
sequence is $B_j={U}_j^{-1}$,
since $\lambda(\theta+1/2)=1/\lambda(\theta)$, also with either
choice of one-sided value. Its left-hand sequence is
$B_j^-=1/{U}_j^-$. Choose a smooth function $\chi:[0,1]\to[0,1]$ equal to
one near $0$ and zero near $1$. Define
$\widehat P_j,\widehat H_j:[0,1]\to(0,\infty)$ by
\[
 \widehat P_j=\chi B_j^-+(1-\chi)B_j,\qquad
 \widehat H_j=\widehat P_j/B_j.
\]
Here $B_j\le1$ by Lemma~\ref{ar:lower}. The calculation for the two one-sided
sequences on the first orbit gives $B_j^-/B_j=1$ if $j\le0$ and
$B_j^-/B_j=e^2$ if $j\ge1$. Hence $\widehat H_j$ is either $1$ or
$1+(e^2-1)\chi(t)$. This proves uniform upper and positive lower
bounds and uniform bounds on all derivatives of any fixed order,
also for $1/\widehat H_j$.
The recurrence coefficients of $B$ and $B^-$ agree except at $0,1$.

\begin{figure}[htbp]
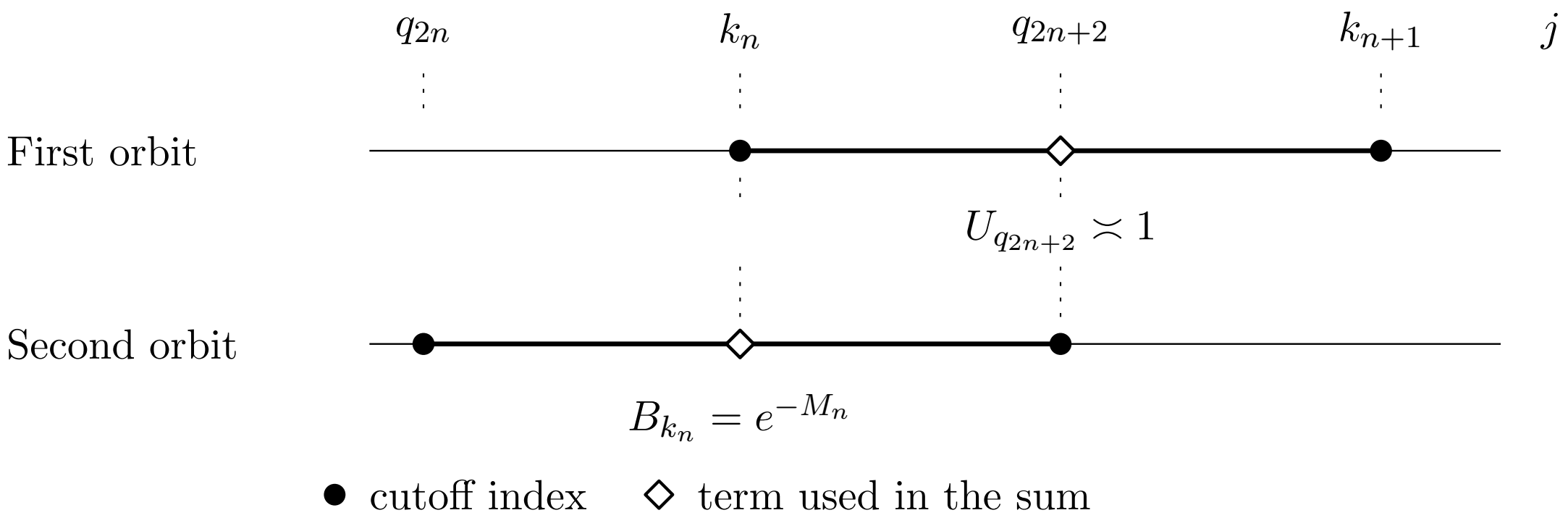

\centering
\begingroup\color{black}\normalsize\SinaiNodesBW\endgroup
\caption{The interpolation blocks on the two orbits, for positive indices. }\label{fig:nodes}
\end{figure}

On this orbit, we obtain a small error in \eqref{pf:nodeidentity}
by choosing small interpolation weights $c_l,c_h$. Since $B_j=1/{U}_j$ is very small
at the maximizing times, each block should contain one such time
in its interior. We accordingly prescribe the densities at the
return times $\pm q_{2n}$, $n\ge n_0$, and at every integer
strictly between the two innermost return times.
Define $\widehat D:[0,1]\to[0,\infty)$ by $\widehat D(t):=t(1-t)$.
Choose the radii and prescribed densities by
\[
\qquad \widehat R_n=e^{-\sqrt{M_n}},\qquad
 \widehat\rho_{\pm q_{2n}}=\widehat P_{\pm q_{2n}}
                                  Z(\widehat D/\widehat R_n).
\]
At each central index use radius one.
Interpolate $\widehat\rho_j/B_j$ by \eqref{pf:interpolate} with
$d_i=(B_iB_{i+1})^{-1}$. The resulting
$\widehat Q_j=(\widehat\rho_{j+1}+\widehat\rho_{j-1})/\widehat\rho_j-2$
equals $q(\widehat\theta_j)$ between successive prescribed indices.
Figure~\ref{fig:nodes} shows the order of the return times and the maximizing times on the two orbits.
For a block $[l,h]$ between return times of orders $n,n+1$,
define the reciprocal block sum
\[
 \widehat c(l,h):=
 \left(\sum_{l\le i<h}\frac1{B_iB_{i+1}}\right)^{-1}.
\]
The block contains $i=k_n$ or $i=-k_n$, together with its next
index. Since $U_i,U_{i+1}\ge c e^{M_n}$,
\[
 0<\widehat c(l,h)\le B_iB_{i+1}\le C e^{-2M_n}.
\]
We use this bound in the second-orbit version of
\eqref{pf:nodeidentity}. The coefficient $\widehat c(l,h)$ differs from
the convex interpolation parameter $\tau$ in \eqref{pf:interpolate}.

For $0<R<1/4$,
\[
 \int_0^1 Z(\widehat D(t)/R)\,dt\le C R(1+|\log R|).
\]
Indeed, on either half of the interval, compare $\widehat D$ with the
distance to its endpoint, split at distance $R$, and integrate $R/t$.
The preceding interpolation and $B_j\le C$ give, on a block $J$
of order $n$,
\begin{equation}\label{pf:recipmass}
 \sum_{j\in J}\int_0^1\widehat\rho_j(t)\,dt
 \le C\exp(O(n\log(n+2)))\widehat R_n(1+|\log\widehat R_n|).
\end{equation}
The logarithm of the bound in \eqref{pf:recipmass} is at most
\[-\sqrt{M_n}+O(n\log(n+2))+\log(1+\sqrt{M_n}).\]
Since $\sqrt{M_n}\asymp n^{5/2}$, this is at most $-c n^{5/2}$
for large $n$. The same estimate holds for the reflected block.
Counting a common endpoint twice only increases the bound on the
total length.

Now fix a return time $j=\pm q_{2n}$, $n>n_0$.
For the adjacent prescribed indices $l<j<h$,
\[
 \widehat c(l,j),\widehat c(j,h)\le C e^{-2M_{n-1}},
 \qquad c\le B_j\le C.
\]
The last bounds follow from Lemma~\ref{ar:grid}. Apply \eqref{pf:nodeidentity}
on the second orbit, using $B$ in place of $A$. Ratios of the cutoffs and
all their weighted derivatives of fixed order are at most
$C_m e^{\sqrt{M_{n+1}}}$; the error is zero on
$\widehat D\le\widehat R_{n+1}$. On the complementary set
$\widehat\rho_j\widehat D\ge c\widehat R_{n+1}$.
Define the parameter error and its expression in the length coordinate by
\[
 \widehat E_j(t):=\widehat Q_j(t)-q(\widehat\theta_j),\qquad
 \widehat e_j:=\widehat E_j\circ\widehat x_j^{-1}
       :[0,\widehat w_j]\longrightarrow\R.
\]
Then
\begin{equation}\label{pf:recipbounds}
 \|\widehat e_{\pm q_{2n}}^{(m)}\|_\infty
 \le \exp(-2M_{n-1}+O_m(1+\sqrt{M_{n+1}}))\longrightarrow0
 \quad(m\ge0).
\end{equation}
To verify this estimate, in the quantity
$[\widehat E_j]_{m,\widehat D}$ the right side of the
identity \eqref{pf:nodeidentity} has a bounded reference factor $B_j^{-2}$, two weights
bounded by $C e^{-2M_{n-1}}$, and weighted ratios bounded by
$C_m e^{\sqrt{M_{n+1}}}$. Also
$[1/\widehat\rho_j]_{m,\widehat D}\le C_m/\widehat\rho_j$ because
$\widehat H_j$ is uniformly bounded above and below with bounded
derivatives. The product expansion used for \eqref{pf:physical}
introduces a factor bounded by $C_m(\widehat\rho_j\widehat D)^{-m}\le
C_m e^{m\sqrt{M_{n+1}}}$ on the support. This proves
\eqref{pf:recipbounds}. Every fixed interval of this second family
has a smooth positive density as a function of $t$.
At their two ends $\widehat Q_j$ is identically
$q(\widehat\theta_j^-)$ and $q(\widehat\theta_j^+)$, respectively:
near each end the three neighboring cutoffs equal one and the
interpolation equals the corresponding one-sided uncut sequence.

\section{Choice of the singular set and the cutoff radii}\label{sec:scales}

We choose a closed set $E$ of interval parameters whose image
will be the nondifferentiability set in each first-orbit interval.
The set $E$ and the radii $R_n$ must satisfy three conditions: the
lengths in \eqref{pf:mass} are summable; the recurrence error,
expressed in the circle coordinate, tends to zero in $C^{r-1}$;
and the averages of $\Theta$ have two distinct limits at each point
of $E$. The first two conditions give the required regularity of the
perturbation. The third gives nondifferentiability on the image of
$E$. Only the first-orbit construction depends on the order of
regularity. We use $|J|$ for the length of an interval and $|B|$ for
the Lebesgue measure of a measurable subset of the line.

For reference, the parameter radii used below are
\[
\begin{array}{c|c|c}
\text{regularity of }\phi&\text{parameter set }E&R_n\text{ for large }n\\ \hline
1<r<\infty&\dim_HE=2/(r+1)&
 \exp[-(M_n-n^{9/2})/(1-d)]\\[2pt]
r=1&\dim_HE=1,\ |E|=0&\exp(-n^{11/2})\\[2pt]
C^\infty&E=\{1/2\}&\exp(-M_n-n^{9/2})
\end{array}
\]
In the first row $d=2/(r+1)$. The second family always uses
$\widehat R_n=e^{-\sqrt{M_n}}$. Here $j$ indexes iterates, whereas
$n$ indexes selected cutoff times. We verify summability of the lengths and the derivative
estimates below.

\subsection{A Cantor set at the prescribed dimension}

For fixed $1\le r<\infty$, write
\begin{equation}\label{sc:exponents}
 s=r-1,\qquad d=\frac{2}{r+1},\qquad
 b=1-d=\frac{s}{s+2},\qquad \tau=\frac{19}{20}.
\end{equation}
Start from $J_0=[1/3,2/3]$. At step $m+1$ retain the two closed
subintervals at the ends of each step-$m$ interval, each of length
$\ell_{m+1}$, where
\begin{equation}\label{sc:lengths}
 \ell_m=|J_0|\,2^{-m/d}\exp(-m^\tau),\qquad m\ge0 .
\end{equation}
Let $E$ be the intersection of the resulting nested unions.
Since $0<(m+1)^\tau-m^\tau\le1$,
\[
 0<\eta_-:=e^{-1}2^{-1/d}
 \le \ell_{m+1}/\ell_m<2^{-1/d}\le 1/2.
\]
At every finite step a nonempty middle interval is removed, even when
$d=1$. The lengths tend to zero and every interval has two descendants,
so $E$ is compact, perfect and contains no interval. We will use the
following elementary count also when $d=1$: an interval of length $h$
meets at most $2+h/\ell_m$ level-$m$ intervals. Apart from the first
and last that it meets, they lie inside it and have disjoint interiors.

Let $E_0=E\cup\{0,1\}$ and use the function $D$ of
Section~\ref{sec:profiles}. The inequalities
$\dist(t,E_0)/2\le D(t)\le\dist(t,E_0)$ let us estimate the
measure of $\{D<v\}$ by covering  neighborhoods of $E_0$.
All constants in the following estimates may depend on $d$ but not
on the intervals, their level, or the point $t_*$.

\begin{lemma}\label{sc:tube}
$\dim_H E=d$, $\mathcal H^d(E)=0$, and in particular $|E|=0$.
For some $C,c>0$,
\begin{equation}\label{sc:global-tube}
 |\{D<v\}|\le C v^{1-d}
                \exp\bigl[-c(\log(1/v))^\tau\bigr],
 \qquad 0<v<1/4 .
\end{equation}
If $d<1$, then for $t_*\in E$, $0<v\le h$ and $0<h<1/24$,
\begin{equation}\label{sc:local-tube}
 |\{D<v\}\cap[t_*,t_*+h]|
        \le C h(v/h)^{1-d}.
\end{equation}
When $d=1$, put $L=\log(1/h)$ and $v_h=h\exp(-\sqrt L)$. Uniformly
for $t_*\in E$ and small $h>0$,
\begin{equation}\label{sc:endpoint-local}
 \frac{1}{h}\big|\{D<v_h\}\cap[t_*,t_*+h]\big|
 \le C\exp[-cL^{\tau-1/2}] .
\end{equation}
\end{lemma}
\begin{proof}
The level-$m$ cover satisfies
$2^m\ell_m^d=|J_0|^d e^{-d m^\tau}\to0$, proving the Hausdorff
measure assertion. For the lower bound on dimension, give each of the two
descendants half the mass of its parent. The resulting probability
measure on $E$ assigns mass $2^{-m}$ to each level-$m$ interval. If
$\ell_m\le |J|<\ell_{m-1}$, then
$|J|/\ell_m<1/\eta_-$ and $J$ meets at most $2+1/\eta_-$ such
intervals by the count above.
For each $0<t<d$,
\[
 \frac{2^{-m}}{\ell_m^t}
 =|J_0|^{-t}\exp\bigl[-(1-t/d)m\log 2+t m^\tau\bigr]\le C_t,
\]
since $m^\tau=o(m)$. Thus the mass of every sufficiently short $J$ is
at most $C'_t|J|^t$. {The interval-cover argument
in Section~\ref{sec:arith} yields
$\dim_H E\ge t$ for every $t<d$.}

Choose $m$ so that $\ell_m\le v<\ell_{m-1}$.
Enlarging the $2^m$ level-$m$ intervals by $2v$ on either side covers
the $2v$-neighborhood of $E$. The endpoints $0,1$ contribute at most
$4v$ more in $[0,1]$. Since $\eta_-v<\ell_m\le v$ and
$\log(1/\ell_m)=\log(1/|J_0|)+(m/d)\log 2+m^\tau$,
$m\asymp\log(1/v)$ at small scales. Thus
\[
 2^m v\le C 2^m\ell_m
  = C |J_0|^d\ell_m^{1-d}e^{-d m^\tau}
  \le C' v^{1-d}e^{-c(\log(1/v))^\tau}.
\]
The extra term $4v$ satisfies the same bound since
$d\log(1/v)-c(\log(1/v))^\tau\to+\infty$. This proves
\eqref{sc:global-tube}.
We also obtain a lower bound with the same power $1-d$. Both endpoints of every level-$m$ interval belong to $E$.
Every point in such an interval therefore has
$D\le\ell_m/2<v$. Hence
\[
 |\{D<v\}|\ge 2^m\ell_m
  = |J_0|^d\ell_m^{1-d}e^{-d m^\tau}
  \ge c v^{1-d}\exp[-C(\log(1/v))^\tau].
\]
In particular $\log|\{D<v\}|/\log v\to1-d$. This determines the leading power in the choice of radii
for the Cantor set defined by \eqref{sc:lengths}.

For the local estimate, select levels $m,N$ by
$\ell_m\le h<\ell_{m-1}$ and $\ell_N\le v<\ell_{N-1}$, so $N\ge m$.
A point in $[t_*,t_*+h]$ with $D<v$ lies within $2v$ of $E$;
the points $0,1$ are too far away to contribute. Only the level-$m$
intervals intersecting $[t_*-2v,t_*+h+2v]$ need to be considered.
There are at most $2+5/\eta_-$ of them: the enlarged interval
has length at most $5h$ and $\ell_m>\eta_-h$. Each has $2^{N-m}$ level-$N$ descendants.
Now $\ell_N/\ell_m\le 2^{-(N-m)/d}$ and $\ell_N>\eta_-v$, so
\[
 2^{N-m}\le(\ell_m/\ell_N)^d\le \eta_-^{-d}(h/v)^d .
\]
An enlarged descendant has length at most $\ell_N+4v\le5v$.
Summing the lengths gives \eqref{sc:local-tube} for $d<1$, including
at endpoints of construction intervals.

For $d=1$ take $v=v_h$. The same count, with the factor from the
lengths in \eqref{sc:lengths} retained, gives
\[
 \frac{|\{D<v\}\cap[t_*,t_*+h]|}{h}
 \le C 2^{N-m}\frac{v}{h}
 \le C' 2^{N-m}\frac{\ell_N}{\ell_m}
 = C'\exp[-(N^\tau-m^\tau)].
\]
Indeed $v<\ell_N/\eta_-$ and $h\ge\ell_m$. Since
$\log(\ell_m/\ell_N)=\sqrt L+O(1)$ and
$\log(\ell_j/\ell_{j+1})\le\log 2+1$, we have
$N-m\ge c\sqrt L$. Also $N=O(L)$, because
$\ell_N>\eta_-v$ and $\ell_N\le |J_0|2^{-N}$. The mean value
theorem gives $N^\tau-m^\tau\ge\tau N^{\tau-1}(N-m)
\ge c'L^{\tau-1/2}$. This proves \eqref{sc:endpoint-local} without
assuming a uniform lower bound on relative gap lengths.
\end{proof}

\begin{lemma}\label{sc:average}
Let $\Theta$ be the function from Section~\ref{sec:profiles}.
Uniformly for $t_*\in E$,
\begin{equation}\label{sc:mean}
 \frac1h\int_{t_*}^{t_*+h}\Theta(t)\,dt
  =\sin(\log\log(C_*/h))+o(1)\qquad(h\downarrow0).
\end{equation}
The same assertion holds when $E=\{1/2\}$ and $D,\Theta$ are
constructed with that set instead.
\end{lemma}
\begin{proof}
Let $L=\log(1/h)$ and $v=h e^{-\sqrt L}$. By
\eqref{sc:local-tube}, when $d<1$ the fraction of $[t_*,t_*+h]$
on which $D<v$ is at most $C e^{-(1-d)\sqrt L}$. At $d=1$ this
fraction tends to zero by \eqref{sc:endpoint-local}. For the
singleton it is at most $4v/h$. On the rest of the interval, $v\le D(t)\le h$ since
$t_*\in E$. Hence
\[
 0\le\log\log(C_*/D(t))-\log\log(C_*/h)
 \le \log\left(1+\frac{\sqrt L}{\log(C_*/h)}\right)\longrightarrow0.
\]
On this set, use $|\sin u-\sin v|\le |u-v|$; on the set
where $D<v$, use the bound $2$. These estimates give,
with constants independent of $t_*\in E$,
\[
 \left|\frac1h\int_{t_*}^{t_*+h}\Theta(t)\,dt
       -\sin\bigl(\log\log(C_*/h)\bigr)\right|
 \le C L^{-1/2}+
 \begin{cases}
  C e^{-(1-d)\sqrt L},&d<1,\\
  C e^{-cL^{\tau-1/2}},&d=1,\\
  C e^{-\sqrt L},&E=\{1/2\}.
 \end{cases}
\]
Here $L=\log(1/h)$, and the last case is the singleton construction.
This proves \eqref{sc:mean}. In particular, along
$h=C_*\exp[-\exp(\pi/2+2\pi j)]$ the averages tend to $1$,
and along $h=C_*\exp[-\exp(3\pi/2+2\pi j)]$ they tend to $-1$.
Both sequences are independent of $t_*\in E$. Since
$E\subset[1/3,2/3]$, every interval $[t_*,t_*+h]$ lies in $(0,1)$
for small $h$, including when $t_*$ is an endpoint of $E$.
\end{proof}

\subsection{Summability and regularity for $1<r<\infty$}

Here $s=r-1>0$ and $0<b<1$. The choice of $d$ follows from
the mass and derivative estimates. To choose the radii, we first
compare the leading exponential factors with $R_n=e^{-\gamma M_n}$. A neighborhood of thickness $R_n$
of a set of dimension $d$ in this construction has length of order
$R_n^{1-d}$ on the logarithmic scale, by Lemma~\ref{sc:tube}. Accordingly the leading exponents
in the mass bound \eqref{pf:mass} and the order-$s$ derivative bound
\eqref{pf:physical} are $1-\gamma b$ and $-2-s+s\gamma$,
respectively, in units of $M_n$; here $M_{n+1}/M_n\to1$.
Both exponents vanish when $\gamma=1/b=(s+2)/s$. To obtain decay in both estimates when these exponents
vanish, we use the factor $e^{-m^\tau}$ in \eqref{sc:lengths}
and add a lower-order correction to $\log R_n$. This calculation
applies to the present construction.

First, using $Z(w)\le C\min(1,w^{-1})$, we obtain
\begin{equation}\label{sc:z-int}
 \int_0^1 Z(D(t)/R)\,dt
 \le C R^b\exp[-c'(\log(1/R))^\tau],
 \qquad 0<R<1/4 .
\end{equation}
Indeed $D<R$ is covered by \eqref{sc:global-tube}. On
$2^iR\le D<2^{i+1}R$ the integral is at most
$C 2^{-i}|\{D<2^{i+1}R\}|$. For $2^{i+1}R\le \sqrt R$
and $R$ sufficiently small, this is at most
\[
 C R^b 2^{-i(1-b)}
      \exp[-c\,2^{-\tau}(\log(1/R))^\tau].
\]
For the remaining indices use $|\{D<w\}|\le C w^b$, valid for
all $w>0$ after increasing $C$. The first such index $i$ satisfies
$2^i>\tfrac12 R^{-1/2}$, so their total contribution is at most
$C R^bR^{(1-b)/2}$. Both geometric sums satisfy \eqref{sc:z-int}.

For all sufficiently large $n$ set
\begin{equation}\label{sc:choice}
 h_n=n^{9/2},\qquad R_n=\exp[-(M_n-h_n)/b].
\end{equation}
By \eqref{ar:increments}, $h_n=o(M_n)$ and
$h_{n+1}-h_n=O(n^{7/2})=o(M_{n+1}-M_n)$. Thus $R_n$ decreases
to zero, and $\log(R_{n-1}/R_n)=O_r(n^4)$. We may take the
first cutoff index sufficiently large, since all central intervals are
treated individually. On either side of the origin, let $J$ be the
block between the cutoff indices of orders $n$ and $n+1$. By \eqref{pf:mass} and \eqref{sc:z-int}, the logarithm of
its mass is at most
\[
 M_{n+1}-M_n+h_n-c' M_n^\tau+O(n\log(n+2)).
\]
Indeed $b\log R_n=-M_n+h_n$ and $\log(1/R_n)\asymp M_n$.
The positive terms are $O(n^{9/2})$, whereas
$M_n^\tau\ge c''n^{19/4}$. Consequently the mass of $J$ is at
most $\exp(-c'''n^{19/4})$ for all large $n$. These bounds are
summable on both sides.

For the errors put $\zeta_n=e^{M_n}R_{n+1}$. Its logarithm
divided by $M_n$ tends to $1-1/b<0$. Thus $\zeta_n<1$ for large
$n$, and for integer $s$ all derivative bounds through order $s$ in
\eqref{pf:physical} are bounded by
$\exp(-2M_n+O_r(n^4))\zeta_n^{-s}$. If $s=m+{\beta}$ is
noninteger, $m=\lfloor s\rfloor$ and $0<{\beta}<1$, apply to
$v=\partial_x^m e_{\pm k_n}$ the inequality
\[
 [v]_{{\beta}}\le (2\|v\|_\infty)^{1-{\beta}}\|v'\|_\infty^{\beta}.
\]
Here $v$ is extended by zero to the real line. The inequality
follows by bounding $|v(x)-v(y)|$ by the smaller of
$2\|v\|_\infty$ and $\|v'\|_\infty|x-y|$; its constant is
independent of the interval length. This is the first-order
interpolation estimate; see also
\cite[Chap.~VIII, Proposition~3.1.1]{Her79}.

In \eqref{pf:physical} the powers of
$\zeta_n^{-1}$ then combine to
$m(1-{\beta})+(m+1){\beta}=s$.
In both cases we obtain, since $s/b=s+2$,
\begin{align}
 \log\|e_{\pm k_n}\|_{C^s}
 &\le -(s+2)M_n+(s/b)(M_{n+1}-h_{n+1})+O_r(n^4)\notag\\
 &=(s+2)(M_{n+1}-M_n-h_{n+1})+O_r(n^4)
   \longrightarrow-\infty .\label{sc:norm}
\end{align}
If the error is identically zero the assertion is immediate.
For large $n$ the last upper bound is at most $-c_r n^{9/2}$.
Thus the errors tend to zero in $C^{r-1}$ for every $r>1$.
In particular $1<r<2$ uses $m=0$ in the same interpolation
inequality; an integer $r\ge2$ uses only the derivative bounds
through order $r-1$.
The two inequalities use the correction $h_n$ with opposite signs:
it must dominate $M_{n+1}-M_n$ in \eqref{sc:norm}, but must be
smaller than $M_n^\tau$ in the length estimate. Explicitly,
\[
 \frac{M_{n+1}-M_n}{h_n}=O(n^{-1/2}),\qquad
 \frac{h_n}{M_n^\tau}=O(n^{-1/4}).
\]
The constants in these estimates and in the $C^{r-1}$ bound need not
be uniform as $r$ varies.

\subsection{The endpoint $r=1$}\label{sc:endpoint-subsection}

Now $d=1$ and $b=0$, so \eqref{sc:choice} is undefined. We keep
the length sequence \eqref{sc:lengths} with $d=1$. By
Lemma~\ref{sc:tube} it gives a Cantor set of Hausdorff dimension one
and Lebesgue measure zero. To obtain a $C^1$ perturbation it suffices
to control the recurrence errors in the supremum norm. We can therefore choose smaller radii than for $r>1$:
\begin{equation}\label{sc:endpoint-radius}
 R_n=\exp(-n^{11/2}).
\end{equation}
First split the parameter interval at $D=\sqrt R$. On $D<\sqrt R$
use $Z\le1$ and \eqref{sc:global-tube}; on the complement use
$Z(D/R)\le C R/D\le C\sqrt R$. Thus for small $R$,
\[
 \int_0^1 Z(D/R)\,dt
 \le C\exp[-c(\log(1/R))^\tau]+C\sqrt R .
\]
As $(11/2)\tau=209/40>5$ and $M_{n+1}=O(n^5)$, inequality
\eqref{pf:mass} now bounds the mass of either block of order $n$
by $\exp(-c' n^{209/40})$ for large $n$. Hence all lengths are
summable.

The radii \eqref{sc:endpoint-radius} satisfy
$\Lambda_n=O(n^{9/2})$ in \eqref{pf:radiusassumption}. Only order
zero of \eqref{pf:physical} is needed here. To verify the error estimate for this choice of radii,
we use \eqref{pf:nodeidentity}.
For either neighbour $i$ of $j=\pm k_n$,
\[
 \sup_t\frac{Z_i^*(t)}{Z_j^*(t)}
 \le C\max(1,R_{n-1}/R_n)\le \exp(Cn^{9/2}).
\]
The bounds on the two interpolation coefficients, on $H_i/H_j$,
and on $A_j$ in Section~\ref{sec:profiles} are independent of the
radii. Taking the supremum in \eqref{pf:nodeidentity} gives
\begin{equation}\label{sc:endpoint-error}
 \|e_{\pm k_n}\|_\infty
 \le \exp[-2M_n+O(n^{9/2})]\longrightarrow0.
\end{equation}
For ${\alpha}\in\mathcal S$ in Remark~\ref{rem:admissible-rotations}, Appendix~\ref{app:frequencies} bounds
the additional arithmetic terms in the logarithms by $O_B(n^2)$; these terms are included in $O(n^{9/2})$, so the same
choices apply.
The extension in Section~\ref{sec:realization} then gives
a continuous derivative. By Lemma~\ref{sc:average}, the local
estimate \eqref{sc:endpoint-local} gives the two limiting averages
at every $t_*\in E$; the condition $|E|=0$ alone would not suffice.
This verifies the three required conditions for $r=1$.

\subsection{One choice of radii for all derivative orders}\label{sc:smooth-subsection}

For the smooth example let $E=\{1/2\}$, so that
$E_0=\{0,1/2,1\}$. Keep $h_n=n^{9/2}$ and now choose
\begin{equation}\label{sc:smooth-radii}
 R_n=\exp[-M_n-h_n].
\end{equation}
The set $\{D<v\}$ has length at most $Cv$ because it lies within
$2v$ of three points. Decomposing as above, the term with $D<R$ and
each nonempty shell $2^iR\le D<2^{i+1}R$ contributes at most $CR$.
Since $D\le1$, there are at most $C(1+\log(1/R))$ such shells
for $0<R<1/4$. Thus $\int_0^1 Z(D/R)\,dt\le CR(1+\log(1/R))$.
Now $\log(1/R_n)=O(n^5)$; taking logarithms in \eqref{pf:mass}
and inserting this integral gives the upper bound
\[
 M_{n+1}-M_n-h_n+O(n\log(n+2))\le-\tfrac12 n^{9/2}
\]
for the logarithm of the mass when $n$ is large. Also
$\log(R_{n-1}/R_n)=O(n^4)$. For every fixed integer $m\ge0$,
\eqref{pf:physical} gives
\[
 \log\|\partial_x^m e_{\pm k_n}\|_\infty
 \le -2M_n+m(M_{n+1}-M_n+h_{n+1})+O_m(n^4)
 =-2M_n+O_m(n^{9/2})\longrightarrow-\infty.
\]
The choice $h_n=n^{9/2}$ dominates
$M_{n+1}-M_n=O(n^4)$ and satisfies $h_n=o(M_n)$, since
$M_n\asymp n^5$. The first relation gives summability of the block
masses. The second gives convergence to zero in every fixed
derivative order: each differentiation contributes
$O(n^{9/2})$ to the logarithmic bound, whereas $M_n\ge c n^5$.
Thus the radii are independent of $m$, as required in
Theorem~\ref{thm:smooth}.


\section{Realization on the circle}\label{sec:realization}

To apply Proposition~\ref{in:herman}, we construct a bi-Lipschitz
circle homeomorphism $g$ from the densities and prove that
$g+g^{-1}-2\Id$ is $C^r$ (or $C^\infty$). For the second assertion,
we find a continuous function $Q$ of the required order that agrees
almost everywhere with the derivative of this sum. The sum is
absolutely continuous, so the conclusion follows from the fundamental
theorem of calculus.

\subsection{Positions of the intervals and the induced homeomorphism}

{We use Herman's construction with the two disjoint orbits of $0$ and
$1/2$: the cumulative mass of the preceding rotation points gives the
left endpoint of each interval \cite[Chap.~X, \S\S3.8--3.11]{Her79}.
Herman allows more than one orbit of intervals
\cite[Chap.~II, Remark~7.2.4]{Her79}.}

The quantities $w_j$ and $\widehat w_j$ were defined in
\eqref{pf:length-coordinate} and \eqref{pf:hat-length-coordinate}.
The preceding estimates give
\begin{equation}\label{gl:total-length}
 S:=\sum_{j\in\Z}(w_j+\widehat w_j)\in(0,\infty).
\end{equation}
This normalization constant is distinct from the orbit sums
$S_j=\log U_j$. It is the total length determined by the densities, whereas
$S_j$ is a Birkhoff sum.

On the rotation circle, define the set of points to be replaced
and the atomic measure that determines the interval lengths by
\begin{equation}\label{gl:placement-measure}
 \begin{aligned}
 \mathcal O&:=\{j\alpha:j\in\Z\}
       \cup\{1/2+j\alpha:j\in\Z\}\pmod1,\\
 \nu&:=\frac1S\sum_{j\in\Z}
       \bigl(w_j\delta_{j\alpha}+\widehat w_j\delta_{1/2+j\alpha}\bigr).
 \end{aligned}
\end{equation}
The two orbits are disjoint and dense, and all weights are positive.
Thus $\nu$ is a probability that gives positive mass to every open arc.
The measure specifies the lengths and circular order of the
intervals.

Take $\eta_*=1/4$, which is not an atom. Define the functions
$a,b:\R\to\R$ first on $[\eta_*,\eta_*+1)$ by
\begin{equation}\label{gl:endpoints}
 a(\eta):=\nu([\eta_*,\eta)),\qquad
 b(\eta):=a(\eta)+\nu(\{\eta\}),
\end{equation}
and then by $a(\eta+m):=a(\eta)+m$, $b(\eta+m):=b(\eta)+m$ for
$m\in\Z$. Thus $a(\eta)$ is the total normalized length preceding
$\eta$; the half-open arc excludes the length assigned to $\eta$ itself.
Explicitly, if $\theta_i^*,\widehat\theta_i^*\in[\eta_*,\eta_*+1)$ are
the representatives of $i\alpha$ and $1/2+i\alpha$, respectively, then
\[
 a(\eta)=\frac1S\left(
 \sum_{\theta_i^*<\eta}w_i+
 \sum_{\widehat\theta_i^*<\eta}\widehat w_i\right)
 \quad(\eta_*\le\eta<\eta_*+1).
\]
For a real lift $j\alpha=\theta_j^*+m_j$ this means
$a(j\alpha)=a(\theta_j^*)+m_j$.

If $\eta<\xi<\eta+1$, then $b(\eta)<a(\xi)$ because
$\nu((\eta,\xi))>0$. The ordered closed intervals $[a(\eta),b(\eta)]$
cover the line and define the collapsing map $h_g:\R\to\R$ by
\begin{equation}\label{gl:collapse}
 h_g(x)=\eta\quad\Longleftrightarrow\quad
                 a(\eta)\le x\le b(\eta).
\end{equation}
This is the continuous nondecreasing degree-one map of Herman's
construction \cite[Chap.~X, \S3.8]{Her79}. A fibre is a nondegenerate
interval precisely when its angle belongs to $\mathcal O$; otherwise
it is the singleton $\{a(\eta)\}$.

Write
\[
 I_j:=(a(j\alpha),b(j\alpha)),\qquad
 \widehat I_j:=(a(1/2+j\alpha),b(1/2+j\alpha))
\]
for the interiors of the fibres on chosen lifts, using the same symbols for their
projections to the circle. The circle decomposes as
\[
 \T=K\sqcup\bigsqcup_{j\in\Z}I_j
          \sqcup\bigsqcup_{j\in\Z}\widehat I_j,
 \qquad
 K:=\T\setminus\bigcup_{j\in\Z}(I_j\cup\widehat I_j).
\]
The interval lengths are $|I_j|=w_j/S$ and
$|\widehat I_j|=\widehat w_j/S$, and their sum is one. Hence $K$
has Lebesgue measure zero. It contains all singleton fibres
and all endpoints of the inserted intervals.
Figure~\ref{appfig:circles} in Appendix~\ref{app:maps} shows
the two circles, the two orbit families, and the collapsing map $h_g$.

Define the interval parameterizations by
\begin{equation}\label{gl:charts}
 \begin{aligned}
 X_j:[0,1]&\longrightarrow\overline I_j,&
 X_j(t)&:=a(j\alpha)+S^{-1}x_j(t),\\
 \widehat X_j:[0,1]&\longrightarrow\overline{\widehat I_j},&
 \widehat X_j(t)&:=a(1/2+j\alpha)+S^{-1}\widehat x_j(t).
 \end{aligned}
\end{equation}
Thus $x_j$ measures unnormalized length from a left endpoint, whereas
$X_j$ is the position in the circle lift. The parameter $t$ labels a
point inside an interval; it is constant along the point's orbit.
In general $t$ is not the relative length, since
\[
 \frac{X_j(t)-X_j(0)}{|I_j|}
 =\frac{\int_0^t\rho_j(v)\,dv}{\int_0^1\rho_j(v)\,dv}.
\]

\begin{lemma}[Regularity of the interval coordinates]\label{gl:chart-regularity}
For each fixed $j\in\Z$, the maps $x_j$ and $X_j$ are increasing
bi-Lipschitz homeomorphisms onto their respective intervals, and
\begin{equation}\label{gl:coordinate-relation}
 X_j(t)=a(j\alpha)+\frac{x_j(t)}S,
 \qquad X_j'(t)=\frac{\rho_j(t)}S\quad\text{for almost every }t.
\end{equation}
They are smooth locally on $[0,1]\setminus E$ and fail to have a
 derivative at every point of $E$. Endpoint smoothness is
understood by extension to a neighborhood of the endpoint. The maps
$\widehat x_j$ and $\widehat X_j$ are smooth interval diffeomorphisms
on all of $[0,1]$.
\end{lemma}
\begin{proof}
The coordinate relation follows from their definitions. For each
fixed $j$, positivity of the prescribed densities and the finite
interpolation formula give constants $0<c_j\le C_j<\infty$ with
$c_j\le\rho_j\le C_j$. Thus, for $s<t$,
\[
 \frac{c_j}S(t-s)\le X_j(t)-X_j(s)
 =\frac1S\int_s^t\rho_j(v)\,dv
 \le\frac{C_j}S(t-s).
\]
This proves bi-Lipschitz continuity and absolute continuity. The
fundamental theorem for the Lebesgue integral gives the derivative
formula almost everywhere. These constants need not be uniform in $j$.

Away from $E$, the density is locally smooth and positive. Near each
endpoint all the finitely many relevant cutoffs are one, and the
endpoint weights in \eqref{pf:raw} are constant. Hence $X_j$ is affine
near $0$ and $1$. Fix $t_*\in E$. On a sufficiently small neighborhood
of $t_*$, interpolation reproduces the uncut density, so that
\[
 \rho_j(t)=A_j\bigl(1+T_j\Theta(t)/4\bigr).
\]
Consequently, for sufficiently small $h>0$,
\[
 \frac{X_j(t_*+h)-X_j(t_*)}{h}
 =\frac{A_j}S\left(1+\frac{T_j}{4h}
                  \int_{t_*}^{t_*+h}\Theta(t)\,dt\right).
\]
By Lemma~\ref{sc:average}, along two sequences $h\downarrow0$ these
quotients tend to
\[
 \frac{A_j}S\left(1+\frac{T_j}4\right)
 \quad\text{and}\quad
 \frac{A_j}S\left(1-\frac{T_j}4\right).
\]
Since $T_j>0$, the two limits differ, so the right derivative
does not exist at $t_*$. The affine coordinate relation gives the
same assertions for $x_j$.

On the second orbit every fixed density $\widehat\rho_j$ is smooth
and strictly positive on $[0,1]$, by Section~\ref{sec:reciprocal}.
Its integral and the affine rescaling of that integral are smooth
with positive derivative, as are their inverse maps.
\end{proof}

Define the maps between closed fibres by
\begin{equation}\label{gl:interval-maps}
 \begin{aligned}
 g_j:\overline I_j&\longrightarrow\overline I_{j+1},&
 g_j&:=X_{j+1}\circ X_j^{-1},\\
 \widehat g_j:\overline{\widehat I_j}&\longrightarrow
             \overline{\widehat I_{j+1}},&
 \widehat g_j&:=\widehat X_{j+1}\circ\widehat X_j^{-1}.
 \end{aligned}
\end{equation}
The global lift $g:\R\to\R$ is specified on these fibres by $g_j$
and $\widehat g_j$, with integer translates of the same formulas.
On a singleton fibre it is defined by
\[
 g(a(\eta)):=a(\eta+\alpha)\qquad(\eta\notin\mathcal O\pmod1).
\]
This defines $g$ on every fibre: the complement of $\mathcal O$
is invariant under the rotation and its inverse, and
\eqref{gl:interval-maps} includes the endpoint values. The maps are increasing
and preserve circular order. Herman's ordered extension
\cite[Chap.~X, \S3.11]{Her79} gives an increasing degree-one
homeomorphism satisfying
\[
 g(X_j(t))=X_{j+1}(t),\qquad
 g(\widehat X_j(t))=\widehat X_{j+1}(t),\qquad
 h_g(g(x))=h_g(x)+\alpha.
\]
It follows that $\rho(g)=\alpha$
\cite[Chap.~II, Proposition~2.10]{Her79}, that the interiors of the fibres are
wandering, and that $K$ is the unique minimal Cantor set
\cite[Chap.~II, Proposition~7.2.3]{Her79}.

The choice of two orbits is specific to the multiplier used here.
Its jumps at $0$ and $1/2$, and the four possible jumps of $q$ in
Lemma~\ref{ar:products}, all belong to these two orbits. Their intervals
allow the required left and right data to be connected. The first
family also contains the prescribed singularities; the second has
smooth internal maps. Other rotation orbits persist as singleton
fibres.

In Herman's notation the interval maps are chosen directly. The
common parameter is an equivalent description: after choosing $X_0$,
one can write $X_j=g^j\circ X_0$. Here the densities are chosen first
and \eqref{gl:interval-maps} then defines the maps. The common parameter expresses the three-term recurrence at
corresponding points in neighbouring intervals. In particular, $h_g(X_j(t))=j\alpha$ is independent of $t$.

Figure~\ref{appfig:densities} in Appendix~\ref{app:maps}
shows how the two density families define the interval coordinates and
the maps between successive fibres.

\vspace{1ex}\paragraph{\bf {Interval densities and the invariant probability.}}
{Let $m_x$ and $m_\theta$ denote normalized Lebesgue measure on the
circle of $g$ and on the rotation circle, respectively. For $0\le s<t\le1$,
\[
 m_x\bigl([X_j(s),X_j(t)]\bigr)=\frac{1}{S}\int_s^t\rho_j(v)\,dv .
\]
Thus $\rho_j/S$ is the density of  interval length in the
parameter $t$. Applying $g$ to the segment replaces $\rho_j$ by
$\rho_{j+1}$, so its length need not be preserved. Collapsing the fibres
gives $(h_g)_*m_x=\nu$: the mass at each rotation point is the length
of its fibre. This is an atomic measure, whereas the invariant
probability of the irrational rotation is $m_\theta$.

The unique $g$-invariant probability is instead the Stieltjes measure
$\mu=dh_g$ \cite[Chap.~II, \S\S7.1, 8.5]{Her79}. More explicitly,
on an arc with lifted endpoints $x<y<x+1$ one has
\[
 \mu([x,y])=h_g(y)-h_g(x),\qquad (h_g)_*\mu=m_\theta .
\]
The increment is unchanged when $x,y$ are replaced by $g(x),g(y)$,
by the lift identity above. Its pushforward under $h_g$ is therefore an
invariant probability for the irrational rotation, hence is $m_\theta$.
Continuity of $h_g$ makes $\mu$ nonatomic.
It is zero on every inserted interval, including the endpoints,
and its support is $K$: every open arc meeting $K$ contains points
in distinct fibres, on which $h_g$ takes different values. Hence
$\mu(K)=1$ whereas $m_x(K)=0$, so $\mu$ is singular with respect to
 length. On the invariant graph the invariant probability is
the pushforward $(\iota_u)_*\mu$ by the graph embedding in
\eqref{pre:maps}. In particular
$\mu(N(u))=0$ in both examples, since Section~\ref{sec:singular}
places every point of $N(u)$ inside a wandering interval.}

\subsection{Local coefficients and the global perturbation}

Define $\mathcal Q_j:I_j\to\R$ by
$\mathcal Q_j:=Q_j\circ X_j^{-1}$; equivalently,
\begin{equation}\label{gl:Q}
 \mathcal Q_j(X_j(t))=Q_j(t)=
 \frac{\rho_{j+1}(t)+\rho_{j-1}(t)}{\rho_j(t)}-2
 \qquad(0<t<1).
\end{equation}
Define $\widehat{\mathcal Q}_j:\widehat I_j\to\R$ by
$\widehat{\mathcal Q}_j:=\widehat Q_j\circ\widehat X_j^{-1}$. Both quotients are smooth functions
of $x$ on their respective open intervals by
Sections~\ref{sec:profiles} and~\ref{sec:reciprocal}. The common normalization
by $S^{-1}$ cancels in \eqref{gl:Q}. Passing from the unnormalized
length coordinate in the derivative estimates to $X_j$ only multiplies
the derivative of order $m$ by $S^m$. For $s=m+\beta$,
$0<\beta<1$, the H\"older seminorm of that derivative is multiplied
by $S^{m+\beta}$: the distance in the unnormalized coordinate is $S$ times the
distance in the normalized coordinate. The same
constant applies to every interval.
When an error $e_j$ originally defined on $[0,w_j]$ is regarded
as a function on $\overline I_j$, it means
$x\mapsto e_j(S(x-a(j\alpha)))$. We retain the symbol $e_j$ for this
affine change of coordinate; the analogous convention applies to
$\widehat e_j$.

The coefficients on both families of intervals
are uniformly bounded. On the first family use \eqref{sc:norm} for $r>1$,
\eqref{sc:endpoint-error} at $r=1$, or the estimates following
\eqref{sc:smooth-radii} for the smooth example; on the second use
\eqref{pf:recipbounds}. Only finitely many exceptional intervals remain,
and each has a bounded coefficient. Positivity in \eqref{gl:Q} gives
$\rho_{j+1}/\rho_j\le C$ for a constant independent of $j$, and at $j+1$
the same identity gives $\rho_j/\rho_{j+1}\le C$. Consequently
\begin{equation}\label{gl:ratios}
 C^{-1}\rho_j(t)\le\rho_{j+1}(t)\le C\rho_j(t),\qquad
 C^{-1}\widehat\rho_j(t)\le\widehat\rho_{j+1}(t)\le C\widehat\rho_j(t)
\end{equation}
for almost every $t$. Since there are only countably many indices,
we may use a common set of full Lebesgue measure for all $j$.

Integrating \eqref{gl:ratios} between any two parameters proves the
same bi-Lipschitz bounds for the restriction of $g$ to a fibre.
The fibrewise definition and its inverse give $g(K)=K$. To
obtain the estimate on an arbitrary arc $J$, split $J\setminus K$
into its countably many open interval pieces. Their images are disjoint
intervals, and $g(J\cap K)\subset K$ has zero length. Therefore
$|g(J)|\le C|J|$. Applying the same length estimate to $g^{-1}$ gives
the reverse bound. In particular the periodic function
\begin{equation}\label{gl:actual}
       \phi(x):=g(x)+g^{-1}(x)-2x
\end{equation}
is Lipschitz and hence absolutely continuous on every compact interval
of the line.

\begin{figure}[htbp]
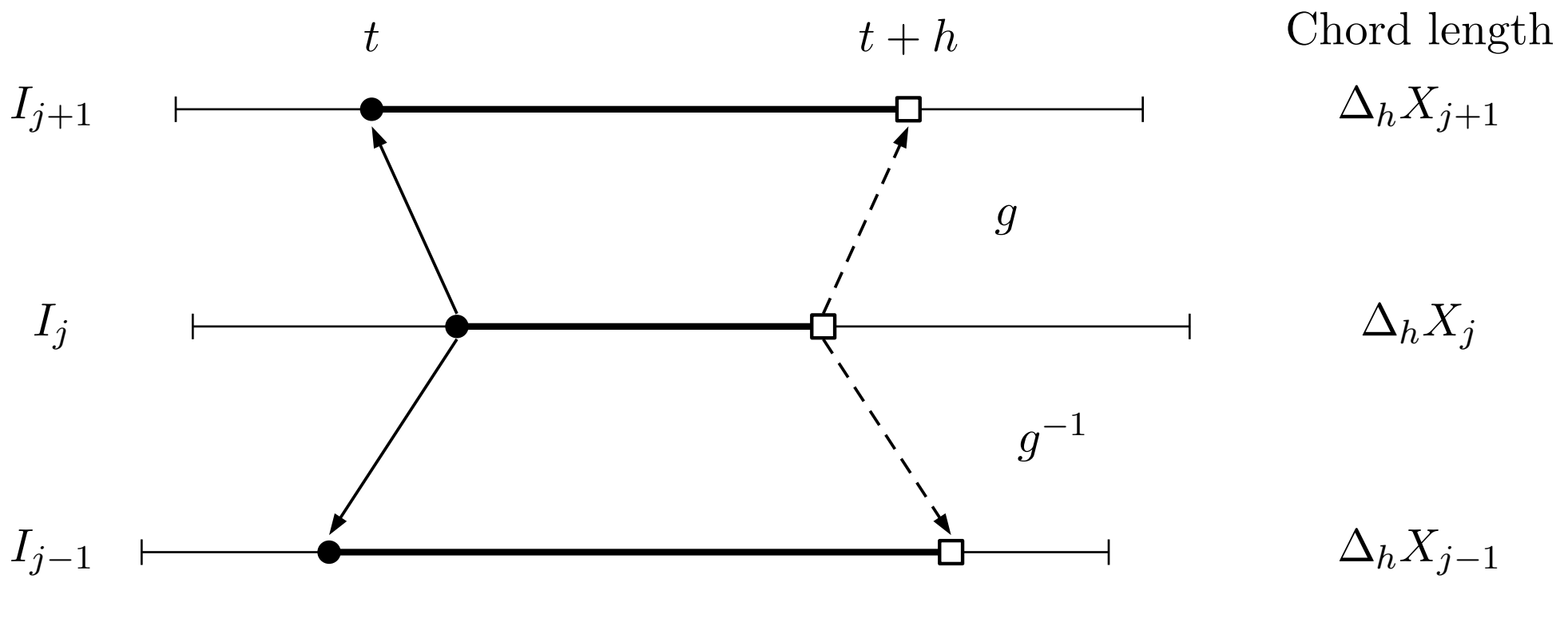

\centering
\begingroup\color{black}\normalsize\SinaiChordsBW\endgroup
\caption{For the same two parameter values $t,t+h$, the three chord lengths are $\Delta_hX_i=X_i(t+h)-X_i(t)$, $i=j-1,j,j+1$. }\label{fig:chords}
\end{figure}

\begin{lemma}[The forward and backward derivatives]\label{gl:forward-backward}
For each $j$, the following derivatives exist at almost every
$t\in(0,1)$:
\begin{equation}\label{gl:chain}
 {\mathrm{D}g}(X_j(t))=\frac{X'_{j+1}(t)}{X'_j(t)},\qquad
 {\mathrm{D}(g^{-1})}(X_j(t))=\frac{X'_{j-1}(t)}{X'_j(t)} .
\end{equation}
At these points,
\begin{equation}\label{gl:density-sum}
 \bigl({\mathrm{D}g}+{\mathrm{D}(g^{-1})}\bigr)(X_j(t))
   =\frac{X'_{j+1}(t)+X'_{j-1}(t)}{X'_j(t)}
   =2+\mathcal Q_j(X_j(t)).
\end{equation}
The same identities hold on the orbit of $1/2$, with $X$ replaced by $\widehat X$.
\end{lemma}
\begin{proof}
Each $X_i$ is the integral of its density, with
$X'_i=\rho_i/S$ almost everywhere. Take $t$ in the intersection
of the full-measure sets where this holds for all $i$, with
$X'_i(t)>0$ there. Figure~\ref{fig:chords} compares the same pair of parameter values on three consecutive intervals. For $s\ne t$, the definitions of the interval maps give
\[
 \frac{g(X_j(s))-g(X_j(t))}{X_j(s)-X_j(t)}
 =\frac{(X_{j+1}(s)-X_{j+1}(t))/(s-t)}
        {(X_j(s)-X_j(t))/(s-t)}.
\]
For the inverse, $g^{-1}(X_j(s))=X_{j-1}(s)$, so
\[
 \frac{g^{-1}(X_j(s))-g^{-1}(X_j(t))}{X_j(s)-X_j(t)}
 =\frac{(X_{j-1}(s)-X_{j-1}(t))/(s-t)}
        {(X_j(s)-X_j(t))/(s-t)}.
\]
Let $s\to t$ from either side. Its denominator tends to
$X'_j(t)>0$, and the numerators tend to $X'_{j+1}(t)$ and
$X'_{j-1}(t)$, respectively. Since $X_j$ is a homeomorphism onto its interval, these limits
give the derivatives at $X_j(t)$. Addition gives \eqref{gl:density-sum}; cancellation of
the common factor $1/S$ gives \eqref{gl:Q}. The argument uses only difference quotients.
In terms of the densities, \eqref{gl:chain} reads
\[
 \mathrm Dg(X_j(t))=\frac{\rho_{j+1}(t)}{\rho_j(t)},\qquad
 \mathrm D(g^{-1})(X_j(t))=\frac{\rho_{j-1}(t)}{\rho_j(t)}.
\]
The inverse derivative is the reciprocal of $\mathrm Dg$ at
$X_{j-1}(t)$, not at $X_j(t)$. These statements hold on the common
full-measure set used above; they make no assertion about individual
derivatives at a parameter in $E$.
\end{proof}

The following identity for increments holds for all parameter
values, including $t\in E$. If $0\le s<t\le1$,
then \eqref{gl:actual}, evaluated at $X_j(t)$ and $X_j(s)$, gives
\begin{equation}\label{gl:all-labels}
 \phi(X_j(t))-\phi(X_j(s))
 =\frac1S\int_s^t\!\bigl(\rho_{j+1}+\rho_{j-1}-2\rho_j\bigr)(v)\,dv
 =\frac1S\int_s^t Q_j(v)\rho_j(v)\,dv .
\end{equation}
Indeed the translation constants in the three coordinates disappear
upon subtraction. Thus, if $\mathcal Q_j$ is continuous at $X_j(t)$,
division by $X_j(t+h)-X_j(t)$ gives a weighted average of $Q_j$
over the interval between $t$ and $t+h$. This tends to $Q_j(t)$ from
either side, since the weights are positive. The derivative of the
sum can therefore exist even where its two individual derivatives
do not.

\begin{lemma}[Exact cancellation near the singular set]\label{gl:local-affine}
For each fixed $j\in\Z$ and each $t_*\in E$, there is an open
interval $V\subset(0,1)$ containing $t_*$ such that
\[
 \phi(x)=a_j^A x+c\qquad(x\in X_j(V))
\]
for a constant $c$ depending on $j$ and $V$. In particular, $\phi$
is smooth on a neighborhood of every point of $X_j(E)$, even though
$g$ and $X_j$ need not be differentiable at the corresponding points.
\end{lemma}
\begin{proof}
Only finitely many positive cutoff radii enter the densities at the
three indices $j-1,j,j+1$. Take $V$ so small that all these cutoffs
are one and $\eta_c=1$ on $V$. For $i=j-1,j,j+1$ we then have
\[
 \rho_i(t)=A_i+\frac{\Theta(t)}4 A_iT_i\qquad(t\in V).
\]
Both $A_i$ and $A_iT_i$ satisfy the recurrence with coefficient
$a_j^A$. Its linearity gives the pointwise identity
\begin{equation}\label{gl:local-density-cancellation}
 \rho_{j+1}(t)+\rho_{j-1}(t)-2\rho_j(t)
       =a_j^A\rho_j(t)\qquad(t\in V).
\end{equation}
For $s<t$ in $V$, the increment identity \eqref{gl:all-labels} yields
\[
 \phi(X_j(t))-\phi(X_j(s))
 =\frac{a_j^A}S\int_s^t\rho_j(v)\,dv
 =a_j^A\bigl(X_j(t)-X_j(s)\bigr).
\]
Since $X_j$ is a homeomorphism, $X_j(V)$ is an open interval, and
the displayed identity proves the asserted affine formula there.
No derivative of $\Theta$, $X_j$, $g$, or $g^{-1}$ at a singular
point has been used.
\end{proof}

Figure~\ref{appfig:regularity} in Appendix~\ref{app:maps} shows
where local affine cancellation applies and where estimates are needed
at the cutoff transitions and near $K$.

In the circle and parameter coordinates, respectively, we have
\[
 \phi|_{X_j(V)}\text{ is affine},\qquad
 (\phi\circ X_j)|_V=a_j^A X_j|_V+c.
\]
The second function need not be differentiable. The neighborhood $V$
may depend on $j$ and on $t_*$. Regularity across $K$ is established below.

\begin{remark}[Relation to the interval constructions]\label{rem:construction-sources}
\leavevmode
\begin{itemize}
\item Herman's construction supplies the circular order, the collapsing
semiconjugacy, and the ordered extension of the interval maps
\cite[Chap.~X, \S\S3.8--3.13]{Her79}. We use its absolute-continuity
argument for the symmetric sum $g+g^{-1}-2\Id$.
\item Arnaud \cite[\S2.3]{Ar11} prescribes distinct one-sided slopes
satisfying the same recurrence, with a modified condition at one index.
Arnaud \cite[\S\S3.1--3.2, equation~(9)]{Ar14} makes the symmetric sum
affine near each corner. Thus her perturbation is already smooth near
every nondifferentiability point of the invariant graph.
\item In our construction, the common coefficient $\Theta(t)/4$ gives
\[
 \rho_j(t)=A_j+\frac{\Theta(t)}4 A_jT_j,
 \qquad
 \rho_{j+1}+\rho_{j-1}-2\rho_j=a_j^A\rho_j
 \quad\hbox{near }E.
\]
The choice of $\Theta$ prescribes the singular parameter set. Cutoffs and
interpolation at selected indices make the lengths summable; the
estimates in the circle coordinate give the extension across $K$.
\end{itemize}
\end{remark}

The functions $\mathcal Q_j,\widehat{\mathcal Q}_j$ are defined
by the density quotients in \eqref{gl:Q}; their interval estimates
come from Sections~\ref{sec:profiles}--\ref{sec:scales}.
We first extend these local functions to a coefficient $Q$ on the
circle and prove its regularity. We then identify $Q$ with
$\mathrm D\phi$, using the absolute continuity of $\phi$ established
in \eqref{gl:actual}.

\vspace{1ex}\paragraph{\bf The reference coefficient.}
The only possible jumps of $q$ are at
$0,\alpha,1/2,1/2+\alpha$. Let
\[
 \mathcal J_0:=\{I_0,I_1,\widehat I_0,\widehat I_1\}.
\]
For $J\in\mathcal J_0$, let $\eta_J$ be its rotation angle and
choose a smooth function $\psi_J$ on $\overline J$ such that
\[
 \psi_J=q(\eta_J^-)\ \hbox{near the left endpoint},
 \qquad
 \psi_J=q(\eta_J^+)\ \hbox{near the right endpoint}.
\]
Define
\[
 q_*(x):=
 \begin{cases}
  \psi_J(x),&x\in\overline J,\ J\in\mathcal J_0,\\
  q(h_g(x)),&x\notin\displaystyle\bigcup_{J\in\mathcal J_0}\overline J.
 \end{cases}
\]
On each complementary component the second expression is constant.
The endpoint conditions therefore give $q_*\in C^\infty(\T)$.
{The four transition functions are constant near their
endpoints, so $q_*$ is locally constant in a neighborhood of each
point of $K$. This step uses the continuity of $h_g$ and the
endpoint values, without differentiating $h_g$.}

Write $Q_J=\mathcal Q_j$ for $J=I_j$ and
$Q_J=\widehat{\mathcal Q}_j$ for $J=\widehat I_j$. The endpoint densities in
Sections~\ref{sec:profiles}--\ref{sec:reciprocal} give
\[
 Q_J=q(\eta_J^-)\ \hbox{near the left endpoint of }J,
 \qquad
 Q_J=q(\eta_J^+)\ \hbox{near the right endpoint of }J,
\]
where $\eta_{I_j}=j\alpha$ and
$\eta_{\widehat I_j}=1/2+j\alpha$. At an interpolated index, this
follows by taking the distance function smaller than the finitely
many cutoff radii involved at that index and its neighbours: the interpolation then agrees with the endpoint densities
before applying the cutoffs.
Consequently
\[
 f_J:=Q_J-q_*|_J\in C^\infty(J),
 \qquad f_J=0\quad\hbox{near both endpoints of }J.
\]
Smoothness inside the exceptional intervals follows from
Lemma~\ref{gl:local-affine} near $X_j(E)$ and from the smooth
coordinates elsewhere. Outside finitely many intervals, $f_J$ is
zero unless $J$ corresponds to a modification index; there it equals
the rescaled error $e_j$ or $\widehat e_j$.

\begingroup
\vspace{1ex}\paragraph{\bf Extension from disjoint intervals.}
We first state the elementary extension result used below. Its
hypotheses concern derivatives in the interval length coordinate.

\begin{lemma}[Extension from disjoint intervals]\label{gl:extension-lemma}
Let $K$ be a nonempty closed subset of $\T$, and enumerate the
components of $\T\setminus K$ as $(J_\ell)_{\ell\ge1}$. Fix $s\ge0$.
For each $\ell$, let $f_\ell\in C^\infty(J_\ell)$ vanish in a
neighborhood of each endpoint, and let $\widetilde f_\ell$ be its
extension by zero to $\T$. Suppose that
\[
 \mathcal M_N:=
 \sup_{\ell>N}\|f_\ell\|_{C^s(\overline J_\ell)}
 \longrightarrow0.
\]
Then
\[
 f:=\sum_{\ell\ge1}\widetilde f_\ell\in C^s(\T),
 \qquad
 \left\|f-\sum_{\ell=1}^N\widetilde f_\ell\right\|_{C^s(\T)}
 \le C_s\mathcal M_N.
\]
Moreover,
\[
 \mathrm D^i f(z)=0
 \qquad(z\in K,\ 0\le i\le\lfloor s\rfloor).
\]
If the hypothesis holds for every integer $s\ge0$ for the same
functions $f_\ell$, then $f\in C^\infty(\T)$.
No measure assumption on $K$ is required.
\end{lemma}

\begin{proof}
Each $\widetilde f_\ell$ is smooth because $f_\ell$ vanishes near
both endpoints. Set
\[
 v_N:=\sum_{\ell=1}^N\widetilde f_\ell,\qquad m:=\lfloor s\rfloor.
\]
At each point, at most one summand or any of its derivatives is
nonzero. Thus, for $M>N$,
\[
 \|\mathrm D^i(v_M-v_N)\|_\infty
 \le\sup_{\ell>N}\|\mathrm D^if_\ell\|_\infty
 \le\mathcal M_N,\qquad 0\le i\le m.
\]
For each $0\le i\le m$, define $V_i$ as the uniform limit of
$\mathrm D^iv_N$. Each $V_i$ is continuous. Since $v_N$ converges
pointwise to the prescribed sum $f$, we have $V_0=f$.
For $i\ge1$, these limits have not yet been identified as derivatives
of $f$. To do so, for $0\le i<m$ and real lifts $x<y$, pass to
the uniform limit in
\[
 \mathrm D^iv_N(y)-\mathrm D^iv_N(x)
 =\int_x^y\mathrm D^{i+1}v_N(t)\,dt
\]
to obtain
\[
 V_i(y)-V_i(x)=\int_x^yV_{i+1}(t)\,dt.
\]
The fundamental theorem of calculus gives $\mathrm D V_i=V_{i+1}$.
Hence $f\in C^m$, with $\mathrm D^if=V_i$ everywhere,
including on $K$. For $m=0$, only uniform convergence and continuity
are needed. Since $\mathrm D^iv_N$ vanishes on $K$, so does $V_i$.
For integer $s=m$, summing the supremum estimates proves the norm
bound with $C_s=m+1$.

Suppose now that $s=m+\beta$, where $0<\beta<1$. Put
$w:=\mathrm D^m(v_M-v_N)$ and take $x\ne y$ in $\T$.
If a shortest joining arc stays in one complementary interval,
that interval's H\"older estimate gives
\[
 |w(x)-w(y)|\le\mathcal M_N\,\dist(x,y)^\beta.
\]
If the arc meets $K$, choose the first points $z_x,z_y\in K$
encountered from $x,y$, respectively; set $z_x=x$ if $x\in K$,
and similarly for $y$. Since $w(z_x)=w(z_y)=0$, the interval
estimates up to the endpoints give
\[
 \begin{aligned}
 |w(x)-w(y)|
 &\le |w(x)-w(z_x)|+|w(y)-w(z_y)|\\
 &\le \mathcal M_N
       \bigl(\dist(x,z_x)^\beta+\dist(y,z_y)^\beta\bigr)\\
 &\le 2\mathcal M_N\,\dist(x,y)^\beta.
 \end{aligned}
\]
This argument also covers points in the same interval when their
shortest joining arc leaves it. Consequently
\[
 \|v_M-v_N\|_{C^s(\T)}\le(m+3)\mathcal M_N.
\]
Letting $M\to\infty$ in the pointwise H\"older inequalities proves
the asserted estimate for $f-v_N$, and hence $f\in C^{m,\beta}$.
Applying the integer-order argument to every fixed $m$ proves
the smooth assertion. The proof uses a supremum of the interval
norms, not their sum.
\end{proof}

\vspace{1ex}\paragraph{\bf Verification of the extension hypothesis.}
We apply Lemma~\ref{gl:extension-lemma} to
$f_\ell=f_{J_\ell}=Q_{J_\ell}-q_*|_{J_\ell}$, enumerating both
families in increasing order of the absolute iterate index.
The finitely many central or exceptional intervals have smooth
errors vanishing near both endpoints. Outside this finite set,
the errors are zero except at the prescribed modification indices,
where they are the affine rescalings of $e_{\pm k_n}$ or
$\widehat e_{\pm q_{2n}}$.

For finite $r\ge1$, the first-family errors tend to zero in
$C^{r-1}$ by \eqref{sc:norm} when $r>1$ and
\eqref{sc:endpoint-error} when $r=1$.
Equation~\eqref{pf:recipbounds} gives the same conclusion for the
second family; for noninteger $r-1$, interpolate between adjacent
integer orders as in the proof of \eqref{sc:norm}.
The fixed coordinate normalization described after \eqref{gl:Q}
preserves these convergences. Each bounded range of cutoff orders
contributes only finitely many nonzero errors. Hence
\begin{equation}\label{gl:error-tail}
 \mathcal M_N:=
 \sup_{\ell>N}\|f_{J_\ell}\|_{C^{r-1}(\overline J_\ell)}
 \longrightarrow0.
\end{equation}
For the smooth construction, Section~\ref{sc:smooth-subsection}
and \eqref{pf:recipbounds} give the same conclusion in every fixed
nonnegative integer order, for the same densities and radii.

\vspace{1ex}\paragraph{\bf The global coefficient.}
Let $\widetilde f_J$ be the extension of $f_J$ by zero. Define
$Q$ pointwise from the local functions by
\begin{equation}\label{gl:global-coefficient}
 Q:=q_*+\sum_{\ell\ge1}\widetilde f_{J_\ell}
 =\begin{cases}
 \mathcal Q_j,&\text{on }I_j,\\
 \widehat{\mathcal Q}_j,&\text{on }\widehat I_j,\\
 q_*,&\text{on }K.
 \end{cases}
\end{equation}
At each point at most one summand is nonzero, so this definition
requires no differentiation. Lemma~\ref{gl:extension-lemma}, applied
using \eqref{gl:error-tail}, gives $Q-q_*\in C^{r-1}(\T)$ and hence
$Q\in C^{r-1}(\T)$. The derivative conclusion of the lemma gives
\[
 \mathrm D^iQ(z)=\mathrm D^iq_*(z)
 \qquad(z\in K,\ 0\le i\le\lfloor r-1\rfloor).
\]
For the smooth construction, the same lemma applies at every
integer order to the same pointwise function $Q$, giving
$Q\in C^\infty(\T)$ and the displayed identity at every order.
This extension uses the interval estimates
and endpoint matching; it does not require $\Leb(K)=0$.

\vspace{1ex}\paragraph{\bf Recovery of the perturbation.}
The function $\phi$ is already Lipschitz and absolutely
continuous by \eqref{gl:actual}. The increment identity
\eqref{gl:all-labels}, and its counterpart on the second family,
identify $\mathrm D\phi$ with the continuous local coefficients on
$\T\setminus K$. By \eqref{gl:global-coefficient} these coefficients
are the restrictions of $Q$. Since $\Leb(K)=0$, we have
$\mathrm D\phi=Q$ almost everywhere on the whole circle.
The fundamental theorem for absolutely continuous functions yields,
for real lifts $x<y$,
\begin{equation}\label{gl:primitive-increment}
 \phi(y)-\phi(x)
 =\int_x^y\mathrm D\phi(t)\,dt
 =\int_x^yQ(t)\,dt.
\end{equation}
{Taking $y=x+1$ in \eqref{gl:primitive-increment}
also gives
\[
 \int_0^1 Q(t)\,dt=0,
\]
by periodicity of the already constructed function $\phi$.
This verifies the compatibility of the extended coefficient with
periodicity of its primitive.}
For every $x$, including $x\in K$, and every sufficiently small
$h\ne0$, it follows that
\[
 \left|\frac{\phi(x+h)-\phi(x)}h-Q(x)\right|
 \le
 \sup_{|t-x|\le|h|}|Q(t)-Q(x)|
 \longrightarrow0.
\]
Continuity of $Q$ therefore gives
\[
 \mathrm D\phi(x)=Q(x)\qquad\text{for every }x\in\T.
\]
For integer $r\ge1$, repeated differentiation of this everywhere
valid identity gives
\[
 \mathrm D^i\phi=\mathrm D^{i-1}Q,\qquad 1\le i\le r,
 \qquad \phi\in C^r.
\]
For noninteger $r=k+\beta$, $k\ge1$ and $0<\beta<1$, we have
$Q\in C^{k-1,\beta}$, so
\[
 \mathrm D^k\phi=\mathrm D^{k-1}Q,\qquad
 [\mathrm D^k\phi]_{\beta,\T}
 =[\mathrm D^{k-1}Q]_{\beta,\T}<\infty.
\]
Hence $\phi\in C^{k,\beta}=C^r$. If $Q\in C^\infty$, the same
argument at every finite order gives $\phi\in C^\infty$.
Absolute continuity is used in~\eqref{gl:primitive-increment} to
recover the increments of $\phi$; it excludes a singular part
concentrated on $K$. No differentiability of $g$ on $K$ is needed
in this argument. This is the zero-measure extension argument in
\cite[Chap.~X, proof of Proposition~3.13]{Her79}, applied here to
the symmetric sum rather than to $g$ itself.
\par\endgroup

\section{The  nondifferentiability set}\label{sec:singular}

Since $u=g-\Id$, we have $N(u)=N(g)$. We prove
\[
 N(g)\cap I_j=X_j(E),\qquad
 N(g)\cap\widehat I_j=\varnothing,\qquad
 g\text{ is strictly differentiable on }K.
\]
The first two assertions concern the interval densities; the third
follows from uniqueness of a normalized positive solution of the
limiting recurrence.

\subsection{Difference quotients inside the wandering intervals}

Fix $j\in\Z$ and $t_*\in E$. Choose a neighbourhood of $t_*$
on which $\eta_c=1$ and all cutoffs used at $j,j+1$ equal one.
The prescribed ratios are affine in $T_i$, so interpolation gives
\[
 \rho_i(t)=A_i\bigl(1+T_i\Theta(t)/4\bigr),
 \qquad i=j,j+1.
\]
Using $g\circ X_j=X_{j+1}$ and integrating the densities, we obtain
\begin{equation}\label{gl:secant}
 \begin{aligned}
 \frac{g(X_j(t_*+h))-g(X_j(t_*))}
      {X_j(t_*+h)-X_j(t_*)}
 &=\frac{A_{j+1}}{A_j}
   \frac{1+T_{j+1}\overline\Theta_h/4}
        {1+T_j\overline\Theta_h/4},\\
 \overline\Theta_h&:=\frac1h\int_{t_*}^{t_*+h}\Theta(t)\,dt .
 \end{aligned}
\end{equation}
By Section~\ref{sec:scales}, including its singleton case, there are
sequences $h_n^\pm\downarrow0$ such that
\[
 \overline\Theta_{h_n^\pm}\longrightarrow\pm1,
 \qquad
 s_j^\pm:=\frac{A_{j+1}}{A_j}
          \frac{1\pm T_{j+1}/4}{1\pm T_j/4}.
\]
Both $s_j^+$ and $s_j^-$ are limits of right difference quotients.
By \eqref{pf:T-increment},
\[
 s_j^+-s_j^-
 =\frac{A_{j+1}}{A_j}
   \frac{(T_{j+1}-T_j)/2}{1-T_j^2/16}
 =\frac1{2WA_j^2(1-T_j^2/16)}>0.
\]
Since $0<X_j(t_*+h_n^\pm)-X_j(t_*)\to0$, the derivative of $g$
at $X_j(t_*)$ does not exist. Outside $X_j(E)$, the adjacent densities
are smooth and positive. The second family has this property throughout
each open interval. Hence
\[
 g|_{I_j\setminus X_j(E)}\in C^\infty,\qquad
 g|_{\widehat I_j}\in C^\infty,
 \qquad
 N(g)\cap I_j=X_j(E).
\]

\subsection{Difference quotients on the minimal set}

Fix $z\in K$ and put $\eta=h_g(z)$. By absolute continuity, it
suffices to find $\ell(z)$ such that
\[
 \omega_z(\delta):=
 \esssup_{|x-z|<\delta}|\mathrm Dg(x)-\ell(z)|
 \longrightarrow0\qquad(\delta\downarrow0).
\]
All intervals and differences in this subsection are taken in local lifts.

\vspace{1ex}\paragraph{\bf The limiting recurrence.}
Let $\Omega_{\mathrm{par}},\widehat\Omega_{\mathrm{par}}\subset(0,1)$
be common full-measure parameter sets on which \eqref{gl:ratios}
and Lemma~\ref{gl:forward-backward} hold for every $j\in\Z$
in the first and second orbit families, respectively. Such sets
are obtained by countable intersection. Define
\[
 \Omega:=\bigcup_{j\in\Z}
 \bigl(X_j(\Omega_{\mathrm{par}})
       \cup\widehat X_j(\widehat\Omega_{\mathrm{par}})\bigr)
 \subset\T\setminus K.
\]
Since the parameterizations are Lipschitz and $\Leb(K)=0$,
\[
 \Leb(\T\setminus\Omega)=0.
\]
First take $x_i\in\Omega$ tending to $z$ through pairwise distinct
intervals. Pass to a subsequence in one family. We write the formulas
for the first family:
\[
 x_i=X_{j_i}(t_i),\qquad |j_i|\longrightarrow\infty,
 \qquad \eta_i=j_i\alpha\pmod1\longrightarrow\eta.
\]
For the second family, use $\widehat\rho$, $\widehat Q$ and
$\eta_i=1/2+j_i\alpha\pmod1$ in the same formulas.
Choose a subsequence approaching $\eta$ from a fixed side
$\epsilon\in\{-,+\}$. Fibre ordering gives
\[
 \begin{aligned}
 z=a(\eta)<b(\eta)&\ \Longrightarrow\ \epsilon=-,\\
 z=b(\eta)>a(\eta)&\ \Longrightarrow\ \epsilon=+.
 \end{aligned}
\]
Indeed, the positive length of the fibre excludes approach through
distinct intervals from its interior side.

Set
\[
 v_{i,k}:=\frac{\rho_{j_i+k}(t_i)}{\rho_{j_i}(t_i)},
 \qquad v_{i,0}=1,
 \qquad C^{-|k|}\le v_{i,k}\le C^{|k|},
 \qquad \mathrm Dg(x_i)=v_{i,1}.
\]
A diagonal subsequence gives $v_{i,k}\to v_k>0$ for every $k\in\Z$,
with $v_0=1$. Before taking the limit, \eqref{gl:Q} gives
\[
 v_{i,k+1}+v_{i,k-1}
   =\bigl(2+Q_{j_i+k}(t_i)\bigr)v_{i,k}.
\]
For each fixed $k$, the indices $j_i+k$ eventually avoid the finitely
many central exceptions. The error estimates and the one-sided
limits of $q$ yield
\[
 \begin{aligned}
 Q_{j_i+k}(t_i)-q(\eta_i+k\alpha)&\longrightarrow0,\\
 q(\eta_i+k\alpha)&\longrightarrow q((\eta+k\alpha)^\epsilon).
 \end{aligned}
\]
Translation by $k\alpha$ preserves the local side of approach, so the
same $\epsilon$ applies to every equation. Consequently
\begin{equation}\label{gl:limit}
 v_{k+1}+v_{k-1}
 =\bigl(2+q((\eta+k\alpha)^\epsilon)\bigr)v_k
 \qquad(k\in\Z).
\end{equation}

\vspace{1ex}\paragraph{\bf Uniqueness of the positive solution.}
Define
\[
 P_0:=1,\qquad
 \frac{P_{k+1}}{P_k}:=\lambda((\eta+k\alpha)^\epsilon).
\]
By the definition of $q$, this positive sequence satisfies
\eqref{gl:limit}. Subtracting the equations for $v$ and $P$ gives
\[
 \begin{aligned}
 P_kv_{k+1}-P_{k+1}v_k
   &=P_{k-1}v_k-P_kv_{k-1}=c,\\
 \frac{v_{k+1}}{P_{k+1}}-\frac{v_k}{P_k}
   &=\frac{c}{P_kP_{k+1}},
 \end{aligned}
\]
where $c$ is independent of $k$. Put
\[
 R_N^+:=\sum_{k=0}^{N-1}\frac1{P_kP_{k+1}},\qquad
 R_N^-:=\sum_{k=-N}^{-1}\frac1{P_kP_{k+1}}.
\]
Summation and positivity imply
\[
 \frac{v_N}{P_N}=1+cR_N^+>0,\qquad
 \frac{v_{-N}}{P_{-N}}=1-cR_N^->0,
 \qquad -\frac1{R_N^+}<c<\frac1{R_N^-}.
\]
At each convergent denominator, Lemma~\ref{ar:grid} gives
\[
 e^{-4}\le P_{\pm q_l}\le e^4,
 \qquad e^{-1}\le\frac{P_{k+1}}{P_k}\le e.
\]
Thus the terms with $k=\pm q_l$ satisfy
$(P_kP_{k+1})^{-1}\ge e^{-9}$.
{The denominators $q_1,\ldots,q_l$ are strictly
increasing. They therefore supply $l$ distinct terms in each sum
below, regardless of whether $q_0=q_1$.} In particular,
\[
 R_{q_l+1}^+\ge le^{-9},\qquad
 R_{q_l+1}^-\ge le^{-9},\qquad
 |c|<\frac{e^9}{l}\longrightarrow0.
\]
Hence $c=0$, $v=P$, and every limit obtained above satisfies
\[
 \lim_{i\to\infty}\mathrm Dg(x_i)=P_1=\lambda(\eta^\epsilon).
\]

\vspace{1ex}\paragraph{\bf Strict differentiability.}
Define the candidate derivative by
\[
 \ell(z):=
 \begin{cases}
 \lambda(\eta),&a(\eta)=b(\eta)=z,\\
 \lambda(\eta^-),&z=a(\eta)<b(\eta),\\
 \lambda(\eta^+),&z=b(\eta)>a(\eta).
 \end{cases}
\]
For a singleton fibre, $\eta$ lies outside both replaced orbits,
so $\lambda(\eta^-)=\lambda(\eta^+)$. At a split endpoint, the
preceding argument determines the limit through distinct intervals.
Inside its own interval, the endpoint densities give, near $z$,
\[
 \mathrm Dg(x)=
 \begin{cases}
 \lambda(\eta^-),&z=a(\eta),\ x>z,\\
 \lambda(\eta^+),&z=b(\eta),\ x<z,
 \end{cases}
 \qquad x\in\Omega.
\]
These values agree with $\ell(z)$.

If convergence to $\ell(z)$ failed, there would be $\varepsilon>0$
and a sequence $x_i\in\Omega$ such that
\[
 x_i\longrightarrow z,\qquad
 |\mathrm Dg(x_i)-\ell(z)|\ge\varepsilon.
\]
A subsequence lies either in an interval adjacent to $z$ or in
pairwise distinct intervals. The endpoint formula in the first
case, and the compactness and uniqueness argument in the second,
give a further subsequence with $\mathrm Dg(x_i)\to\ell(z)$.
This contradiction proves
\[
 \lim_{\substack{x\to z\\x\in\Omega}}\mathrm Dg(x)=\ell(z),
 \qquad \omega_z(\delta)\longrightarrow0.
\]
For $z-\delta<x<y<z+\delta$, absolute continuity now yields
\[
 \begin{aligned}
 \left|\frac{g(y)-g(x)}{y-x}-\ell(z)\right|
 &=\left|\frac1{y-x}\int_x^y
                  \bigl(\mathrm Dg(v)-\ell(z)\bigr)\,dv\right|\\
 &\le\omega_z(\delta)\longrightarrow0.
 \end{aligned}
\]
Thus $g$ is strictly differentiable at $z$, and the strict derivative
of $u=g-\Id$ is $\ell(z)-1$.

\subsection{Proofs of the two theorems}

The preceding subsections give
\begin{equation}\label{gl:exactN}
 N(u)=N(g)=\bigcup_{j\in\Z}X_j(E).
\end{equation}
The definitions of the interval maps and the semiconjugacy give
\[
 \begin{gathered}
 X_j(E)=g^j(X_0(E)),\qquad I_j=g^j(I_0),
 \qquad I_j\cap I_k=\varnothing\quad(j\ne k),\\
 h_g\circ g=R_\alpha\circ h_g,\qquad \rho(g)=\alpha.
 \end{gathered}
\]
Moreover, $u=g-\Id$ is periodic and Lipschitz. By
Proposition~\ref{in:herman},
\[
 \int_0^1\phi(x)\,dx=0,\qquad
 F_\phi\circ\iota_u=\iota_u\circ g.
\]
The primitive following that proposition proves exactness.
The graph is strictly differentiable on $K$.
Figure~\ref{appfig:regularity} summarizes the location of the singular
set and the regularity of the perturbation.

\begin{proof}[Proof of Theorem~\ref{thm:finite}]
Choose $E$ and the radii as in Section~\ref{sec:scales}, and put
$d=2/(r+1)$. Section~\ref{sec:realization} gives
\[
 Q\in C^{r-1},\qquad \mathrm D\phi=Q,
 \qquad \phi\in C^r,
\]
where $C^{r-1}=C^0$ when $r=1$. Each $X_j$ is bi-Lipschitz, so
\[
 \dim_H X_j(E)=\dim_H E=d,\qquad
 \mathcal H^d(X_j(E))=0.
\]
In particular, $N(u)\cap I_0=X_0(E)$ is a compact Cantor set.
By \eqref{gl:exactN}, countable stability and subadditivity give
\[
 \begin{aligned}
 \dim_H N(u)&=\sup_{j\in\Z}\dim_H X_j(E)=d,\\
 \mathcal H^d(N(u))&\le\sum_{j\in\Z}\mathcal H^d(X_j(E))=0.
 \end{aligned}
\]
Together with $X_j(E)=g^j(X_0(E))$, this proves all assertions.
\end{proof}

\begin{proof}[Proof of Theorem~\ref{thm:smooth}]
Take $E=\{1/2\}$ and the radii \eqref{sc:smooth-radii}.
The smooth case of the extension following
\eqref{gl:global-coefficient} and the recovery argument following
\eqref{gl:primitive-increment} give
\[
 Q\in C^\infty,\qquad \mathrm D\phi=Q,\qquad \phi\in C^\infty.
\]
For $x_*:=X_0(1/2)\in I_0$, \eqref{gl:exactN} gives
\[
 N(u)=\bigcup_{j\in\Z}\{X_j(1/2)\}
     =\{g^j(x_*):j\in\Z\}.
\]
Thus $u$ is not $C^1$. The construction uses one map for all derivative
orders; Theorem~\ref{thm:finite} allows the map to depend on $r$.
\end{proof}


\appendix
\begingroup
\renewcommand{\thefigure}{A.\arabic{figure}}
\renewcommand{\theHfigure}{A.\arabic{figure}}
\setcounter{figure}{0}
\section{Maps and coordinates in the construction}\label{app:maps}

The three figures below describe the circles, the interval coordinates, and the maps that define the invariant graph and the perturbation. They use the notation of
Sections~\ref{sec:profiles}--\ref{sec:singular}. All circles are copies
of $\T=\R/\Z$; formulas involving real differences use compatible
degree-one lifts.

\subsection{The circles and their fibres}\label{app:circles}

Figure~\ref{appfig:circles} distinguishes the rotation coordinate
$\theta$ (or $\eta$) from the coordinate $x$ on the circle of $g$.
The rotation $R_\alpha(\theta)=\theta+\alpha$ orders the two families
$\{j\alpha\}$ and $\{1/2+j\alpha\}$, with $j\in\Z$ and angles
taken modulo one. These are exactly the points with nondegenerate fibres.

\begin{SinaiAppendixFigure}
\centering
\SinaiOverviewCirclesBW
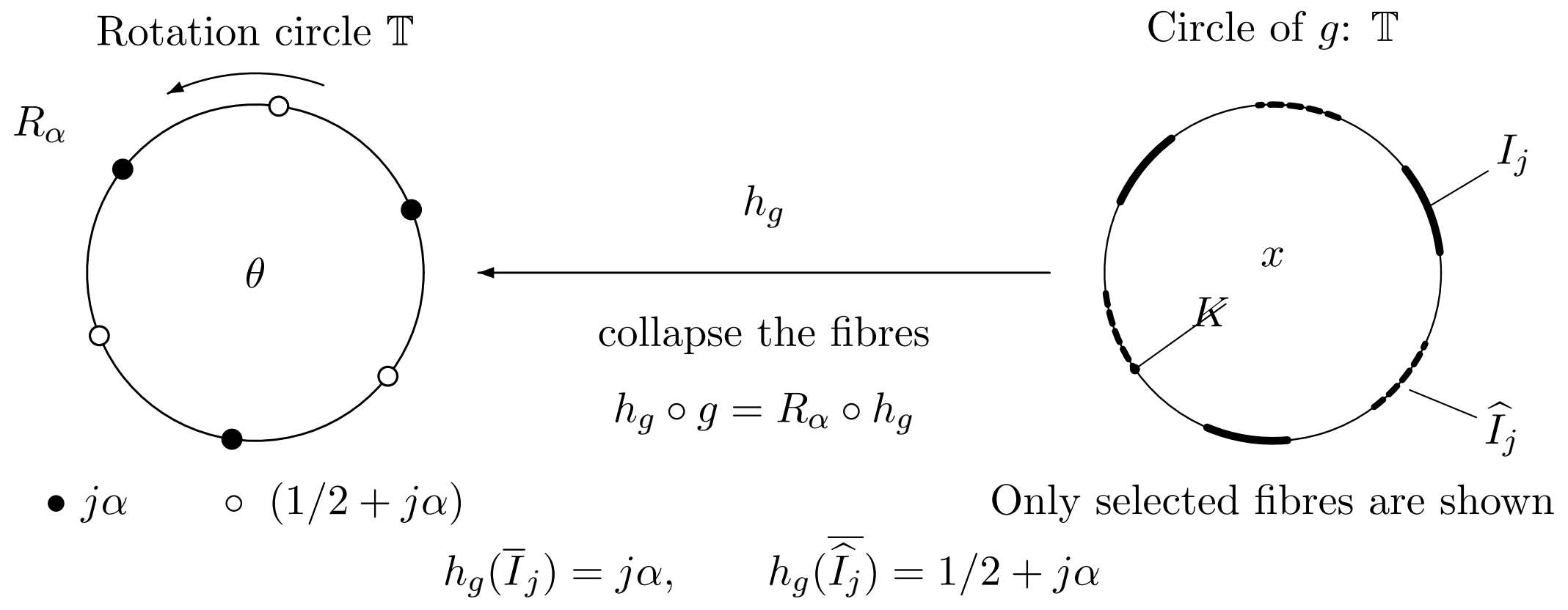
\caption{The two circles and the collapsing semiconjugacy.
}
\label{appfig:circles}
\end{SinaiAppendixFigure}

The endpoint maps, defined from the probability $\nu$ in
\eqref{gl:placement-measure}, are
\[
 a(\eta)=\nu([\eta_*,\eta)),\qquad
 b(\eta)=a(\eta)+\nu(\{\eta\}),\qquad \eta_*=1/4,
\]
for $\eta_*\le\eta<\eta_*+1$, extended by translation by one.
The closed fibre $[a(\eta),b(\eta)]$ has interior $I_j$ over $j\alpha$
and $\widehat I_j$ over $1/2+j\alpha$; otherwise it is a singleton.
The maps $a,b$ need not be continuous.

The continuous nondecreasing degree-one map $h_g$ collapses each
closed fibre to its angle. The homeomorphism $g$ sends the fibre over
$\eta$ to the fibre over $\eta+\alpha$, with
$h_g\circ g=R_\alpha\circ h_g$. The complement $K$ of the open
fibres consists of singleton fibres and interval endpoints. It is the
minimal Cantor set and has Lebesgue measure zero.

\subsection{Interval coordinates and densities}\label{app:densities}

Figure~\ref{appfig:densities} shows how integration of the densities,
followed by normalization and translation, places the parameter
intervals on the circle.

\begin{SinaiAppendixFigure}
\centering
\SinaiOverviewDensitiesBW
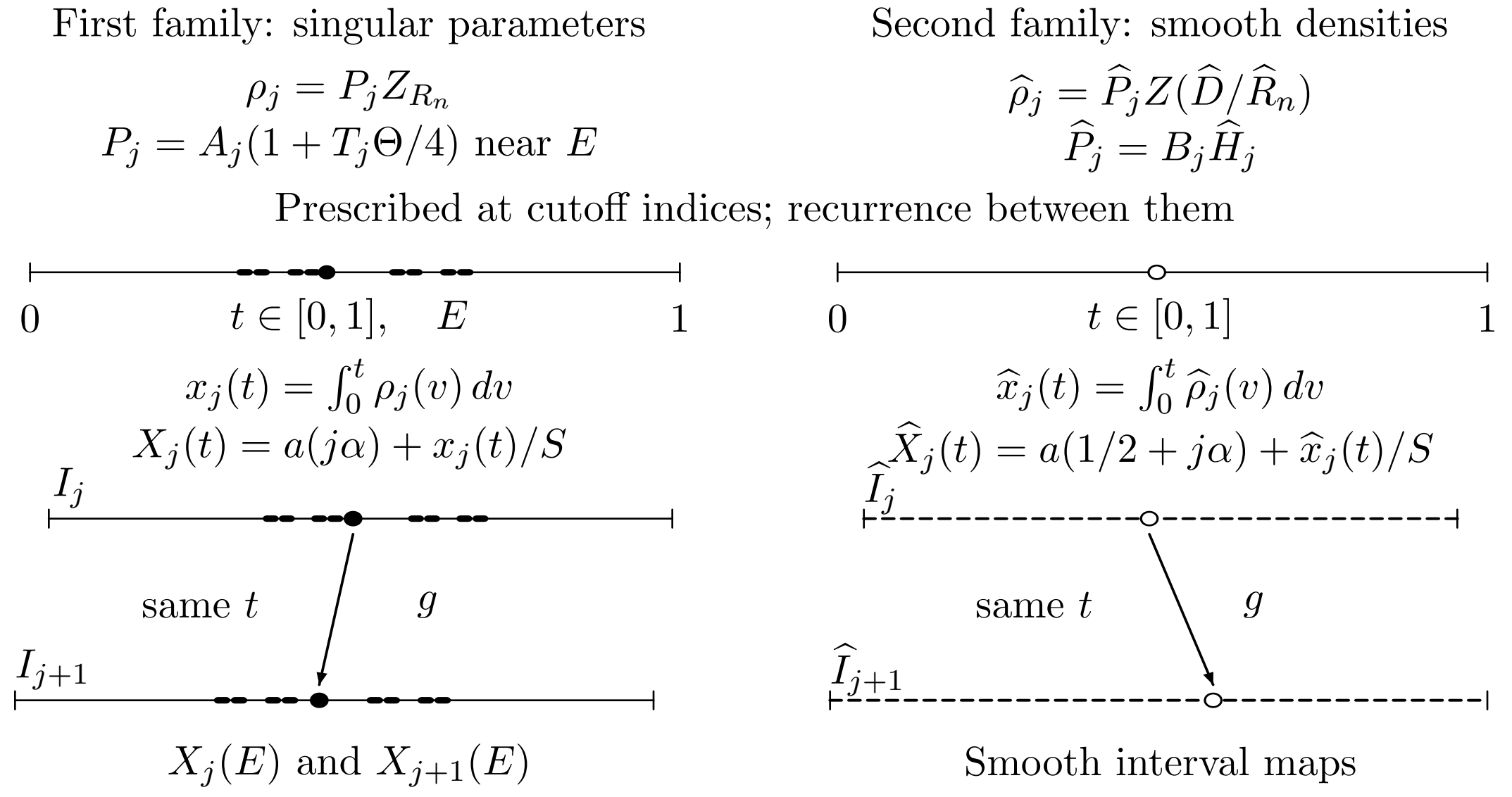
\caption{Densities, coordinates and interval maps. }
\label{appfig:densities}
\end{SinaiAppendixFigure}

The integrals $x_j,\widehat x_j$ measure unnormalized length;
$w_j=x_j(1)$, $\widehat w_j=\widehat x_j(1)$ and
$S=\sum_j(w_j+\widehat w_j)$. The maps $X_j,\widehat X_j$ place
these coordinates in $\overline I_j,\overline{\widehat I_j}$ by
\eqref{gl:charts}. Each $X_j$ is bi-Lipschitz, with constants depending
on $j$; each $\widehat X_j$ is a smooth increasing parameterization.

The restrictions of $g$ are
$g_j=X_{j+1}\circ X_j^{-1}$ and
$\widehat g_j=\widehat X_{j+1}\circ\widehat X_j^{-1}$;
see \eqref{gl:interval-maps}. Iteration changes the index $j$ and leaves $t$ fixed. The densities $\rho_j/S,\widehat\rho_j/S$ describe
 length, whereas the invariant probability $\mu=dh_g$ is
supported on $K$.

The uncut densities are $P_j,\widehat P_j$.
Near $E$, $P_j=A_j(1+T_j\Theta/4)$ combines the positive recurrence
solutions $A_j,A_jT_j$ of \eqref{pf:linear-recurrence} with the common
oscillating function $\Theta$. On the second orbit $B_j$ is the
positive reference sequence and $\widehat H_j=\widehat P_j/B_j$
is its smooth density quotient. The fixed positive cutoff $Z$ equals
one near zero. We write $Z_{R_n}=Z(D/R_n)$, where $D$ is comparable to the distance from
$E\cup\{0,1\}$; on the second orbit $\widehat D(t)=t(1-t)$.
The radii $R_n,\widehat R_n$ are parameter scales, not circle lengths.
Sections~\ref{sec:profiles}--\ref{sec:reciprocal} specify the cutoff
indices and the interpolation.

\subsection{The singular set and the perturbation}\label{app:regularity}

Figure~\ref{appfig:regularity} shows how the densities determine the singular set and the regularity of the perturbation. The set $E\subset(0,1)$ belongs to the parameter
interval; its image $X_j(E)$ belongs to the wandering interval
$I_j$ on the circle.

\begin{SinaiAppendixFigure}
\centering
\SinaiOverviewRegularityBW
\caption{Recovery of the graph and the perturbation.
}
\label{appfig:regularity}
\end{SinaiAppendixFigure}

The identity map is denoted by $\Id$. The periodic function
$u=g-\Id$ defines the graph, and the Herman--Mather formula gives
$\phi=g+g^{-1}-2\Id$. The quotient $Q_j$ in \eqref{pf:Q} is the
derivative of $\phi$ after passage to the circle coordinate:
$\mathrm D\phi(X_j(t))=Q_j(t)$ for $0<t<1$. The increment identity \eqref{gl:all-labels} proves this equality even when $g$ or $X_j$ is not differentiable at the corresponding point.
The corresponding identity holds on the second family of intervals.

Near each point of $X_j(E)$ the three neighbouring densities
solve the same second-order recurrence. Their quotient $Q_j$ is
therefore constant, so $\phi$ is affine on a neighborhood of that
point (Lemma~\ref{gl:local-affine}). The sizes of these neighborhoods
may depend on $j$ and on the point. Global $C^r$ regularity requires the additional
estimates across $K$, proved in Section~\ref{sec:realization}.
Section~\ref{sec:singular} proves
$N(u)=N(g)=\bigcup_jX_j(E)$ and strict differentiability on $K$.
There are no singular points in the second family.

Finally, the standard twist map $F_\phi$ and the graph embedding
$\iota_u$ are
\[
 F_\phi(x,y)=(x+y,\,y+\phi(x+y)),\qquad
 \iota_u(x)=(x,u(x)).
\]
They satisfy $F_\phi\circ\iota_u=\iota_u\circ g$.
The projection $\pi_1(x,y)=x$ recovers the circle dynamics as
$g=\pi_1\circ F_\phi\circ\iota_u$; the domains are those in
\eqref{in:map} and~\eqref{pre:maps}. These identities relate the circle map to the dynamics on the invariant graph.

\endgroup

\section{Admissible rotation numbers}
\label{app:frequencies}

\subsection{Rotation estimates}
Keep $S_j=\log {U}_j$ as in \eqref{ar:fields}. For fixed $B,k_0$, put
$Q_k=q_{2k^2}$ and $M_k^*=\max_{0\le j<Q_{k+1}}S_j$. We prove that
\begin{align}
 c k^4\le M_k^*-M_{k-1}^*&\le C_B k^4,
 &M_k^*&\asymp_B k^5,\label{freq:peak}\\
 \log Q_{k+1}&=O_B(k^2),
 &\min_{|j|\le Q_{k+1}}S_j&\ge-C_B k^2.\label{freq:loss}
\end{align}
The estimates below hold for all sufficiently large $k$; the constants may also depend on the finite initial segment. Lemma~\ref{ar:grid} applies to every
irrational rotation. {For incomplete blocks we use a separate grid comparison from
\eqref{ar:convergents}:} among the first $i\le q_n$ signs, a translation by
at most $C/q_n$ changes only $O_C(1)$ signs. Indeed, only points
within $C/q_n$ of a discontinuity can change sign; each of the
$q_n$ points is within $1/q_n$ of a different grid point.

At a designated position $m=2k^2$, the recurrence for $q_m$ and
coprimality of consecutive denominators show that $h=q_m$ is odd.
Write $a=a_{m+1}$, $H=q_{m+1}=ah+q_{m-1}$ and
$\delta=h{\alpha}-{p}_m\in(0,1/H)$. At time $lh+i$, $0\le i<h$,
the position differs from $i{p}_m/h$ by $(l+i/h)\delta$ to the right.
For $0\le l<a$ this is less than $1/h$. No grid point crosses zero;
the only one that can cross $1/2$ is $(h-1)/(2h)$. For $l<a/2$
the displacement is less than $1/(2h)$, so the first $a/2$
complete blocks each have sum $+1$. The remaining complete blocks have
sums $+1$ or $-1$.
{Since $a$ is even, the sum of the first $l$ complete
blocks is bounded below by
\[
 \begin{cases}
 l,&0\le l\le a/2,\\
 a-l,&a/2\le l\le a.
 \end{cases}
\]
In particular, the cumulative sum at every block start is
nonnegative.} The last partial
block, like any earlier one, differs from the corresponding prefix at
zero by $O(1)$, since $l\delta<1/h$ for $l\le a$.

It follows that passing from $[0,h)$ to $[0,H)$ reduces the minimum
by at most a constant and increases the maximum by at most $a+O(1)$.
For the lower bound, use in block $\lfloor a/8\rfloor$ a prefix on which the maximum over $[0,h)$ is attained: all the signs agree up to that time
in the two blocks, and the intervening complete blocks add
$\lfloor a/8\rfloor$. Extending to $q_{m+2}=2H+h$ changes either extremum by at most $O(1)$: there are at most two complete $H$-blocks,
each bounded by Lemma~\ref{ar:grid}, and the last partial block starts
within $2/H$ of zero. Whenever both succeeding digits are bounded by
$B$, the same reasoning gives a change of $O_B(1)$ on extending
$[0,q_i)$ to $[0,q_{i+2})$: there are at most $(B+1)^2$ complete
$q_i$-blocks, and the last starting displacement is $O_B(1/q_i)$.

Between the designated large digits at successive stages there are
$O(k)$ such bounded pairs. Maxima cannot decrease when the prefix
is extended. Summation gives \eqref{freq:peak} and the lower bound
in \eqref{freq:loss} at positive times; the additional endpoint changes the sum by at most one. For negative times use $S_{-j}=S_{j+1}-1$. Finally
$q_i\le\prod_{v=1}^i(a_v+1)$ gives the denominator bound.

\subsection{Passage to the circle construction}
Let $t_k$ attain $M_k^*$. For large $k$ one has
$Q_k<t_k<Q_{k+1}<t_{k+1}$: the previous maximum is smaller and
Lemma~\ref{ar:grid} bounds $|S_{Q_k}|$. The same lemma and $|S_{j+1}-S_j|=1$ show that the distance from $t_k$ to either endpoint is at least $M_k^*-C$. Use $\pm t_k$ on the first orbit and $\pm Q_k$ on the
second. Thus Sections~\ref{sec:profiles}--\ref{sec:scales} are
used with the substitutions
\[
 (q_{2n},k_n,M_n)\longmapsto(Q_k,t_k,M_k^*).
\]
The uniform bound ${U}_j\ge1$ used for \eqref{in:frequency}
is no longer assumed. On the blocks of orders $k-1,k,k+1$,
\eqref{freq:loss} gives the replacements
\[
 |j|+{U}_j^{-1}+1+|L_j|\le\exp(O_B(k^2)),
 \qquad H_j\ge c/(1+|L_j|).
\]
For the first inequality sum the increments $({U}_i{U}_{i+1})^{-1}$ in
\eqref{pf:AT}; for the second use \eqref{pf:raw} and the fact that
${U}_j^-/{U}_j$ is still either $1$ or $e^{-2}$.

At $j=\pm t_k$, $\log A_j=M_k^*+O_B(k^2)$, using
$|S_{-t_k}-S_{t_k}|\le2$; at an intervening return $\pm Q_k$ both
$A_j,A_{j+1}\le\exp(O_B(k^2))$. Thus in \eqref{pf:conductance}
and \eqref{pf:mass}, $O_B(k^2)$ replaces $O(k\log(k+2))$.
The weighted derivatives of $1/H_j$ of any fixed order are at most
$\exp(O_{B,m}(k^2))$ by the reciprocal differentiation formula
following \eqref{pf:nodeidentity}. Consequently these replacements
add only $O_{B,m}(k^2)$ to the logarithms in the first-orbit
differentiation estimates. On the second orbit we have $B_j={U}_j^{-1}\le\exp(O_B(k^2))$, whereas $B_{\pm Q_k}\asymp1$, and
$B_j,B_{j+1}\le Ce^{-M_k^*}$ at $j=\pm t_k$. Thus the bound for
the errors at returns still holds, and the logarithm of the mass
bound \eqref{pf:recipmass} is at most
\[
 -\sqrt{M_k^*}+O_B(k^2)+O(\log M_k^*).
\]
Since $\sqrt{M_k^*}\asymp_B k^{5/2}$, the masses are summable.
For $1<r<\infty$ and for the smooth example the choices
\eqref{sc:choice} and \eqref{sc:smooth-radii} therefore give the
same conclusions: the successive increments are $O_B(k^4)$ and the
correction $k^{9/2}$ dominates both those increments and $O_B(k^2)$.
At $r=1$ only the supremum norm of the first-orbit error is needed;
the endpoint estimates in Section~\ref{sec:scales} remain valid after adding $O_B(k^2)$ to the logarithms of the error and the mass. Finally the uniqueness of the positive
solution on the minimal set uses Lemma~\ref{ar:grid} at arbitrary
starting points and thus holds for these frequencies as well.

\subsection{Dimension and irrationality exponent}
The inverse branch of a continued-fraction word $w$ has derivative
of absolute value $(q+q't)^{-2}$, $0<t<1$, where $q,q'$ are its
last two denominators. Its distortion is at most four. Prescribing
the next digit $A$ retains at most $4/[A(A+1)]$ times the length
of each parent interval. At $A=2(k+2)^4$ this bound tends to zero.
Since the parents have disjoint interiors, each fixed $(B,k_0)$
family, and hence their countable union, has Lebesgue measure zero.

To estimate the dimension, let $D_B$ consist of the numbers whose digits are all at most $B$. Let $\mathcal T_B:D_B\to\mathcal S$ be the map
that inserts the prescribed triples into their digit sequences.
If two outputs have a longest common prefix of length $n$, their next
digits differ at a free position, since the prescribed digits depend
only on the past. Let $Q$ and $q$ be the denominators of the full
common prefix and of its free digits, respectively. For finite words
$u,v$, multiplication of continued-fraction matrices gives
$q(u)q(v)\le q(uv)\le2q(u)q(v)$.
There are $O(\sqrt n)$ inserted triples, of size polynomial in $n$.
Applying the inequality to the successive words and triples gives
\[
 Q\le q\exp(O_B(\sqrt n\log(n+2))),\qquad q\ge e^{c n}
\]
for large $n$; the latter bound follows from the Fibonacci recursion
on the $n-O(\sqrt n)$ free digits. If the first unequal digits are
adjacent, the two tails have the form $1/(a+u)$, $1/(a+1+v)$.
Here $a+1\le B$, $0<u,v<1$, and $v\ge c_B/(n+2)^2$ because the
following digit is $O_B((n+2)^2)$. Thus their separation is at
least $v/(B+1)^2$. Nonadjacent digits are separated by at least
$1/(B+1)^2$.
{If $Q'$ is the preceding denominator of the common
output prefix $w$, its inverse branch $f_w$ satisfies
\[
 |f_w'(t)|=(Q+Q't)^{-2}\ge\frac1{4Q^2},
 \qquad 0<t<1,
\]
because $0\le Q'\le Q$. Hence the two outputs are separated by
at least $c_B(n+2)^{-2}Q^{-2}$, whereas the inputs are at distance
at most $q^{-2}$. For every $0<\gamma<1$,
\[
 \log\frac{q^{-2}}{\bigl((n+2)^{-2}Q^{-2}\bigr)^\gamma}
 \le -2(1-\gamma)\log q
      +O_B(\sqrt n\log(n+2))+2\gamma\log(n+2).
\]
The right-hand side tends to $-\infty$ since $\log q\ge cn$.
Absorbing the finitely many smaller prefix lengths into the
constant proves}
\[
 |x-y|\le C_{B,\gamma}
       |\mathcal T_B(x)-\mathcal T_B(y)|^\gamma
 \qquad(x,y\in D_B).
\]
Applying this inverse H\"older bound to interval covers gives
\[
 \dim_H\mathcal T_B(D_B)\ge\gamma\dim_H D_B.
\]
Letting $\gamma\uparrow1$ proves
$\dim_H\mathcal T_B(D_B)\ge\dim_H D_B$.

{By Hensley's dimension estimate for bounded partial quotients
\cite{Hensley}, $\lim_{B\to\infty}\dim_H D_B=1$.}

For an arbitrary open interval $J$, choose a finite prefix whose cylinder lies in $J$ and start the prescribed positions after this prefix. Its inverse branch
is bi-Lipschitz; the preceding dimension estimates on the remaining
free digits are unchanged. Increasing $B$ proves $\dim_H(\mathcal S\cap J)=1$.

The irrationality exponent of an irrational number $\alpha$ is
\[
 \mu_{\mathrm{irr}}(\alpha)=\sup\bigl\{\nu>0:
 |\alpha-p/q|<q^{-\nu}\text{ for infinitely many reduced fractions }p/q,
 \ q\ge1\bigr\}.
\]
We say that $\alpha$ is of \emph{Roth type} if, for every
$\varepsilon>0$, there is $c_{\alpha,\varepsilon}>0$ such that
\begin{equation}\label{freq:roth-condition}
 \left|\alpha-\frac pq\right|
 \ge \frac{c_{\alpha,\varepsilon}}{q^{2+\varepsilon}}
 \qquad(p\in\Z,\ q\in\mathbb N,\ q\ge1).
\end{equation}
For irrational numbers this is equivalent to
$\mu_{\mathrm{irr}}(\alpha)=2$; the constant is allowed to depend on
$\varepsilon$. This is the Diophantine condition used in
\cite[Chap.~V, \S\S7.8.1--7.8.2]{Her79}.

For $\alpha\in\mathcal S$, the prescribed digits give
$a_{n+1}=O_\alpha((n+2)^2)$. Since the denominators $q_n$ grow at
least exponentially, for each $\varepsilon>0$ there is
$C_{\alpha,\varepsilon}>0$ such that
\[
 a_{n+1}+2\le C_{\alpha,\varepsilon}q_n^\varepsilon
 \qquad(n\ge0).
\]
Given $q\ge1$, choose $n$ with $q_n\le q<q_{n+1}$.
The best-approximation property of convergents and their standard
error bounds imply, for every $p\in\Z$,
\[
 |q\alpha-p|\ge |q_n\alpha-p_n|
 >\frac1{q_{n+1}+q_n}
 \ge\frac1{(a_{n+1}+2)q_n}
 \ge\frac1{C_{\alpha,\varepsilon}q_n^{1+\varepsilon}}.
\]
Dividing by $q$ and using $q_n\le q$ proves
\eqref{freq:roth-condition}. On the other hand, the convergents satisfy
$|\alpha-p_n/q_n|<q_n^{-2}$ for infinitely many $n$.
Consequently $\mu_{\mathrm{irr}}(\alpha)=2$.


\subsection{Liouville frequencies}
For a different family keep $a_{2n+1}=2(n+2)^4$ but allow arbitrary
even-position digits subject to $\log q_{2n+2}=O(n^2)$.
{Lemma~\ref{ar:lower} applies because every odd-indexed digit is even.}
Again $h=q_{2n}$ is odd and $H=q_{2n+1}$ even. The preceding
large-digit computation gives an increase of order $n^4$ up to $H$.
For $0<j<q_{2n+2}$ the displacement from $j{p}_{2n+1}/H$ is strictly
between $-1/H$ and zero by \eqref{ar:convergents}. Both discontinuities
belong to this even grid, so every sign agrees with the left limit
at its grid point. At the initial time $j=0$ we use the right value instead.
Each complete period of the left grid signs has sum zero. Extending
from $H$ to $q_{2n+2}$ therefore changes the maximum by at most a
constant, independently of $a_{2n+2}$. This proves the estimates
above with $Q_n=q_{2n}$.

For an example take $N_1\ge2$, $N_{k+1}=N_k^3$, and put
$a_{2N_k+2}=\lfloor e^{N_k^2}\rfloor$, taking all other
even-position digits equal to two. Then $\log q_{2n+2}=O(n^2)$, but
just before the exceptional digit $\log q_{2N_k+1}
=O(N_k\log N_k+\sum_{i<k}N_i^2)=O(N_k\log N_k)$.
For large $k$, the next digit exceeds every fixed power of this denominator. The corresponding approximation errors are therefore smaller than every fixed negative power of the denominator along a subsequence. Thus $\alpha$ is a Liouville number.

\section{Sharpness under local affine cancellation}\label{app:sharp}

This appendix gives an upper bound for a class containing the examples
of Theorem~\ref{thm:finite}. The additional condition concerns the
perturbation near the singular orbits. It does not prescribe the
dimension of the singular set. We use the convention for regularity
spaces stated in the introduction.

\begin{proposition}[Upper bound on a wandering interval]
\label{sharp:upper}
Let $1\le r<\infty$, let $\phi:\R\to\R$ be a $1$-periodic function
of class $C^r$, and suppose that the standard twist map $F_\phi$
preserves the graph of a Lipschitz function $u:\T\to\R$.
Assume that the induced circle homeomorphism $g$ preserves orientation,
and let $I\subset\T$ be an open wandering interval. Suppose that
\begin{equation}\label{sharp:affine-condition}
 \begin{gathered}
 \text{for every }x\in N(u)\cap I\text{ and every }j\in\Z,
 \text{ there is a neighborhood }V\text{ of }g^j(x)\\
 \text{on which }\phi\text{ is affine in a real coordinate.}
 \end{gathered}
\end{equation}
The neighborhood may depend on $x$ and $j$; no uniform radius is
assumed. Then
\begin{equation}\label{sharp:upper-conclusion}
 \mathcal H^{2/(r+1)}\bigl(N(u)\cap I\bigr)=0,
 \qquad
 \dim_H\bigl(N(u)\cap I\bigr)\le\frac{2}{r+1}.
\end{equation}
No assumption on the rotation number is needed for this upper bound.
\end{proposition}

\subsection{The linearized equation and a slope criterion}

We first prove the properties of the linearized equation needed below. All the notation introduced in this subsection is
local to Appendix~\ref{app:sharp}; in particular, the operator
$\mathcal J_x$ below is distinct from the scalar density functions
$H_j$ in Section~\ref{sec:profiles}.

Choose increasing real lifts. Invariance of the graph gives
$g(x)=x+u(x)$ and
\[
 \phi(g(x))=u(g(x))-u(x)=g^2(x)-2g(x)+x.
\]
Substituting $g^{-1}(x)$ for $x$ gives
$g+g^{-1}-2\Id=\phi$. The first identity shows that $g$ is Lipschitz, and the second shows that $g^{-1}=2\Id+\phi-g$ is Lipschitz. Choose $L\ge1$ such that both have Lipschitz constant
at most $L$. The derivatives of every fixed iterate exist almost
everywhere, and
\[
 L^{-1}\le \mathrm Dg(x)\le L
 \quad\text{at every point where }\mathrm Dg(x)\text{ exists}.
\]
The lower bound follows by applying the Lipschitz inequality for
$g^{-1}$ to $g(x+h)$ and $g(x)$. By deleting the countably many images
and inverse images of the exceptional null sets, we obtain an invariant
set $X_*$ of full Lebesgue measure on which the derivatives and chain
rules for all iterates hold. This deletion preserves full measure
because every fixed iterate is bi-Lipschitz.

For $x\in\R$ put $\xi_j:=g^j(x)$. The orbit equation is
\begin{equation}\label{sharp:orbit-equation}
 \xi_{j+1}+\xi_{j-1}-2\xi_j=\phi(\xi_j),\qquad j\in\Z.
\end{equation}
For a real square-summable sequence $v=(v_j)_{j\in\Z}$ define
\begin{equation}\label{sharp:jacobi}
 (\mathcal J_xv)_j
 :=\bigl(2+\mathrm D\phi(\xi_j)\bigr)v_j-v_{j-1}-v_{j+1}.
\end{equation}
Here $\ell^2(\Z)$ is the space of sequences satisfying
$\sum_jv_j^2<\infty$, with inner product
$\langle v,w\rangle:=\sum_jv_jw_j$. Since $\mathrm D\phi$ is bounded, $\mathcal J_x$ is a bounded symmetric operator on this space.
Its quadratic form is
\begin{equation}\label{sharp:quadratic-form}
 \mathcal E_x(v):=\langle\mathcal J_xv,v\rangle
 =\sum_j\bigl((v_{j+1}-v_j)^2+\mathrm D\phi(\xi_j)v_j^2\bigr).
\end{equation}

\begin{lemma}[Positivity and capacities]\label{sharp:capacity}
For every $x$, the quadratic form $\mathcal E_x$ is nonnegative.
Define
\begin{equation}\label{sharp:delta}
 \Delta(x):=
 \inf\{\mathcal E_x(v):v_0=1, v\text{ has finite support}\}.
\end{equation}
Then $\Delta\ge0$ is upper semicontinuous,
\begin{equation}\label{sharp:capacity-bound}
 \mathcal E_x(v)\ge\Delta(x)v_0^2
 \qquad(v\in\ell^2(\Z)),
\end{equation}
and
\begin{equation}\label{sharp:singular-inclusion}
 N(u)\subset\{x:\Delta(x)>0\}.
\end{equation}
Moreover, $\Delta(x)=0$ implies strict differentiability of $u$ at $x$.
\end{lemma}

\begin{proof}
For $x\in X_*$, the positive sequence $p_j:=\mathrm D(g^j)(x)$ satisfies
$p_0=1$ and $\mathcal J_xp=0$, by differentiating
\eqref{sharp:orbit-equation}. For a finitely supported $v$, expansion
gives
\begin{equation}\label{sharp:positive-identity}
 \mathcal E_x(v)=
 \sum_jp_jp_{j+1}
 \left(\frac{v_{j+1}}{p_{j+1}}-\frac{v_j}{p_j}\right)^2.
\end{equation}
Indeed, the coefficient of $v_j^2$ on the right is
$(p_{j+1}+p_{j-1})/p_j=2+\mathrm D\phi(\xi_j)$, and the cross terms are
$-2v_jv_{j+1}$. Thus the form is nonnegative on $X_*$. For fixed
finite support its value is continuous in $x$; since $X_*$ is dense,
nonnegativity holds at every $x$. Approximation by finite truncations
then proves it on $\ell^2(\Z)$. The same approximation, with $v_0$
unchanged, proves \eqref{sharp:capacity-bound} from
\eqref{sharp:delta}. If $v_0=0$, this is just nonnegativity.

We next relate the capacity to slopes by considering finite intervals. Every nonempty finite principal matrix of
$\mathcal J_x$ is positive definite. Otherwise, a nonzero finitely
supported $v$ would have $\mathcal E_x(v)=0$. The inequality
$\mathcal E_x(v+tw)\ge0$ for all real $t$ and all finitely supported
$w$ would give $\mathcal J_xv=0$ at every index. If $j_0$ is the first
nonzero index of $v$, its equation at $j_0-1$ reads $-v_{j_0}=0$, a
contradiction.

These finite matrices also have the following maximum principle:
if $Mz=f\ge0$ componentwise, then $z\ge0$. To see this, put
$(z_+)_j:=\max(z_j,0)$ and $(z_-)_j:=\max(-z_j,0)$.
The off-diagonal entries of $M$ are nonpositive, so
$\langle Mz_+,z_-\rangle\le0$. If $z_-\ne0$, positive definiteness
would imply
\[
 0\le\langle f,z_-\rangle
 =\langle Mz_+,z_-\rangle-\langle Mz_-,z_-\rangle<0.
\]
This proves the assertion.

For $n\ge1$, minimize $\mathcal E_x(v)$ under
$v_0=1$ and $v_j=0$ for $|j|\ge n$. Denote the unique minimizer by
$V^{(n)}(x)$ and its value by $c_n(x)$. Existence and uniqueness follow
from the positive definiteness just proved. Finite matrix inversion
also shows that every component of $V^{(n)}$ and $c_n$ is continuous
in $x$. The minimizer satisfies
\[
 (\mathcal J_xV^{(n)})_j=0\quad(0<|j|<n),
 \qquad (\mathcal J_xV^{(n)})_0=c_n(x).
\]
The second equality follows by multiplying the equations by the
components of $V^{(n)}$ and summing.

Define the continuous functions
\begin{equation}\label{sharp:slope-bounds}
 \underline m_n(x):=V^{(n)}_1(x)-1,
 \qquad
 \overline m_n(x):=1+\mathrm D\phi(x)-V^{(n)}_{-1}(x).
\end{equation}
For $n=1$ the two components at indices $\pm1$ mean zero.
The equation at zero gives
\begin{equation}\label{sharp:finite-gap}
 \overline m_n(x)-\underline m_n(x)=c_n(x).
\end{equation}
For $x\in X_*$, compare $p$ and $V^{(n)}$ separately on the positive
and negative half-intervals. Their difference vanishes at index zero,
has positive values $p_n,p_{-n}$ at the outer endpoints, and satisfies
the homogeneous equation at all intervening indices. The finite
maximum principle gives
$p_{\pm1}\ge V^{(n)}_{\pm1}$. Since
$p_1+p_{-1}=2+\mathrm D\phi(x)$ and $\mathrm Du(x)=p_1-1$, it follows that
\begin{equation}\label{sharp:ae-squeeze}
 \underline m_n(x)\le \mathrm Du(x)\le \overline m_n(x)
 \quad\text{for almost every }x.
\end{equation}
Extending an admissible vector by zero proves $c_{n+1}\le c_n$, and
every finite support is contained in one such interval. Hence
$\Delta=\inf_n c_n$, which is upper semicontinuous as an infimum of
continuous functions.

Suppose now that $\Delta(x)=0$. Absolute continuity of $u$ and
\eqref{sharp:ae-squeeze} imply, for $y<z$ in a local coordinate,
\begin{equation}\label{sharp:secant-squeeze}
 \inf_{[y,z]}\underline m_n
 \le\frac{u(z)-u(y)}{z-y}
 \le\sup_{[y,z]}\overline m_n.
\end{equation}
The difference quotients based at $x$ are bounded and therefore have cluster values,
and every such value belongs to $[\underline m_n(x),\overline m_n(x)]$ for every $n$.
By \eqref{sharp:finite-gap} these intervals have lengths tending to
zero, so all cluster values are one number $m$. Given $\epsilon>0$,
choose $n$ with $c_n(x)<\epsilon$. Continuity of the two endpoints
and \eqref{sharp:secant-squeeze} then show that every secant slope with distinct endpoints sufficiently close to $x$ differs
from $m$ by at most $2\epsilon$. Thus $\mathrm Du(x)=m$ and all such secant slopes converge to this value. This is strict differentiability and
proves \eqref{sharp:singular-inclusion}.
\end{proof}

\subsection{The orbit measure and the covering estimate}

\begin{proof}[Proof of Proposition~\ref{sharp:upper}]
For $r=1$, the conclusion follows from Rademacher's theorem: a
Lipschitz function is differentiable almost everywhere, and
$\mathcal H^1$ on the circle agrees with length. Assume now that $r>1$ and put $d:=2/(r+1)$.

For each Borel set $B\subset I$, define
\begin{equation}\label{sharp:orbit-measure}
 \nu_I(B):=\sum_{j\in\Z}|g^j(B)|,
\end{equation}
where $|\cdot|$ denotes Lebesgue measure on the circle. Each summand
is a Borel measure on $I$, because $g^j$ is a homeomorphism.
The wandering property gives $\nu_I(I)\le1$. Moreover,
$\nu_I\ll\operatorname{Leb}$: every fixed Lipschitz iterate takes
Lebesgue-null sets to null sets. The nondifferentiability set is
Borel, by the Cauchy criterion for difference quotients, and is null
by Rademacher's theorem. Consequently
\begin{equation}\label{sharp:nu-null}
 \nu_I(N(u)\cap I)=0.
\end{equation}

Choose $x<z$ in a lift of $I$ and set
\[
 \xi_j:=g^j(x),\qquad \zeta_j:=g^j(z),\qquad h_j:=\zeta_j-\xi_j>0.
\]
The interval $g^j([x,z])$ lies in a lift of the wandering interval
$g^j(I)$, so its length is $h_j$. Since $\nu_I$ has no atoms,
\begin{equation}\label{sharp:l1-chords}
 \sum_jh_j=\nu_I([x,z])<\infty.
\end{equation}
In particular $h\in\ell^1(\Z)\subset\ell^2(\Z)$. Subtracting the
orbit equations gives
\[
 h_{j+1}+h_{j-1}-2h_j=\phi(\zeta_j)-\phi(\xi_j).
\]
The identity
$\sum_j(h_{j+1}-h_j)^2
=-\sum_jh_j(h_{j+1}+h_{j-1}-2h_j)$ is justified by absolute
convergence. Substituting into \eqref{sharp:quadratic-form}, we obtain
\begin{equation}\label{sharp:endpoint-identity}
 \frac{\mathcal E_x(h)+\mathcal E_z(h)}2
 =\sum_jh_j^2
 \left(
 \frac{\mathrm D\phi(\xi_j)+\mathrm D\phi(\zeta_j)}2
 -\frac{\phi(\zeta_j)-\phi(\xi_j)}{h_j}
 \right).
\end{equation}
All terms are absolutely summable: the expression in parentheses is
bounded by a constant depending on $\|\mathrm D\phi\|_\infty$, and
$h\in\ell^2$.

Let $a$ be a point satisfying the local affine condition \eqref{sharp:affine-condition}. If $r=k+\beta$ with an integer $k\ge1$ and
$0<\beta<1$, then $f:=\mathrm D\phi$ belongs to $C^{k-1,\beta}$. Local
affineness gives $f^{(m)}(a)=0$ for $1\le m\le k-1$.
Taylor's formula with a H\"older remainder therefore yields
\[
 |f(a+t)-f(a)|\le C_r|t|^{k-1+\beta}=C_r|t|^{r-1}.
\]
For $k=1$ this is simply the H\"older inequality for $f$.
If $r=k\ge2$ is an integer, Taylor's integral remainder, using
boundedness of $f^{(k-1)}=\mathrm D^k\phi$ and the vanishing lower
coefficients, gives the same bound with exponent $k-1$.
In both cases,
\begin{equation}\label{sharp:flat-taylor}
 |\mathrm D\phi(a+t)-\mathrm D\phi(a)|\le C_r|t|^{r-1}
 \qquad(|t|\le1).
\end{equation}
The constant depends on the global $C^r$ norm of the periodic
function, not on $a$ or on the radius of its affine neighborhood.
The neighborhood is used only to identify the derivatives at $a$.
Since
\[
 \frac{\phi(a+h)-\phi(a)}h
 =\int_0^1\mathrm D\phi(a+th)\,dt,
\]
\eqref{sharp:flat-taylor} implies the uniform estimate
\begin{equation}\label{sharp:trapezoidal}
 \left|
 \frac{\mathrm D\phi(a)+\mathrm D\phi(a+h)}2
 -\frac{\phi(a+h)-\phi(a)}h
 \right|
 \le C_r|h|^{r-1},\qquad 0<|h|\le1.
\end{equation}

Fix $x\in N(u)\cap I$. Each $a=\xi_j$ satisfies
\eqref{sharp:affine-condition}, and $0<h_j<1$. Combining
\eqref{sharp:capacity-bound}, \eqref{sharp:endpoint-identity}, and
\eqref{sharp:trapezoidal}, we obtain
\begin{align}
 \frac{\Delta(x)+\Delta(z)}2|z-x|^2
 &\le C_r\sum_jh_j^{r+1}\notag\\
 &\le C_r\left(\sum_jh_j\right)^{r+1}
 =C_r\nu_I([x,z])^{r+1}.
 \label{sharp:key-estimate}
\end{align}
The elementary inequality used in the second line follows first for
finite sums of nonnegative numbers and then by increasing limits.
In particular, if $\Delta(x)\ge\varepsilon>0$, then
\begin{equation}\label{sharp:measure-lower}
 \begin{gathered}
 \nu_I([x,z])\ge c_{\varepsilon,r}|z-x|^d
 \\{\bigl(x\in N(u)\cap I,\quad
       \Delta(x)\ge\varepsilon,\quad x<z\in I\bigr).}
 \end{gathered}
\end{equation}

We prove the zero Hausdorff-measure assertion by a covering argument. Fix a compact subinterval
$J\Subset I$ and set
\[
 E_\varepsilon:=J\cap N(u)\cap\{\Delta\ge\varepsilon\}.
\]
By \eqref{sharp:nu-null}, this set has zero $\nu_I$-measure. Given
$\eta>0$, regularity of the finite Borel measure $\nu_I$ on $I$
provides an open set $O\subset I$ containing $E_\varepsilon$ with
$\nu_I(O)<\eta$. For every center $x\in E_\varepsilon$, take
arbitrarily small intervals $B(x,h)\subset O$, with their radii also
less than any prescribed $\delta>0$. Apply
\eqref{sharp:measure-lower} to $z=x+h/2$. This gives
\[
 \nu_I(B(x,h))\ge c_{\varepsilon,r}'h^d.
\]
Select a countable disjoint subfamily $B_i=B(x_i,h_i)$ such that the
fivefold intervals $5B_i$ cover $E_\varepsilon$. To construct such a subfamily, consider the radii in successive dyadic ranges,
choosing a maximal disjoint family in each range among intervals
not intersecting those already chosen. Every rejected interval meets
a chosen interval of radius at least half its own radius, and is
contained in the fivefold enlargement of that chosen interval.
This proves the required covering property.

Since all the selected intervals lie in $O$,
\[
 \sum_i\bigl(\operatorname{diam}(5B_i)\bigr)^d
 \le10^d\sum_i h_i^d
 \le C_{\varepsilon,r}\sum_i\nu_I(B_i)
 \le C_{\varepsilon,r}\eta.
\]
Their diameters are at most $10\delta$. Let $\eta$ and $\delta$
tend to zero. This proves $\mathcal H^d(E_\varepsilon)=0$.
Finally exhaust $I$ by compact subintervals and take
$\varepsilon=1/m$, $m\ge1$. Inclusion
\eqref{sharp:singular-inclusion} and countable subadditivity prove
\eqref{sharp:upper-conclusion}.
\end{proof}

\subsection{Attainment and the regularity of the perturbation}

Fix the frequency $\alpha$ in \eqref{in:frequency}.
For $1\le r<\infty$, let $\mathcal C_r(\alpha)$ be the class of
triples $(\phi,u,I)$ satisfying the hypotheses of
Proposition~\ref{sharp:upper}, with induced rotation number $\alpha$.
In particular, condition \eqref{sharp:affine-condition} is part of
the definition of this class.

\begin{proposition}[Sharpness and exact regularity]
\label{sharp:optimal}
For every $1\le r<\infty$,
\begin{equation}\label{sharp:maximum}
 \max_{(\phi,u,I)\in\mathcal C_r(\alpha)}
 \dim_H\bigl(N(u)\cap I\bigr)=\frac2{r+1}.
\end{equation}
Moreover, the perturbation $\phi_r$ constructed in
Theorem~\ref{thm:finite} for the exponent $r$ satisfies
\begin{equation}\label{sharp:exact-regularity}
 \phi_r\in C^r,
 \qquad \phi_r\notin C^s\quad\text{for every }s>r.
\end{equation}
Hence $r$ is the exact regularity exponent of this perturbation.
\end{proposition}

\begin{proof}
Proposition~\ref{sharp:upper} gives the upper bound in
\eqref{sharp:maximum}. The construction in the proof of
Theorem~\ref{thm:finite} provides a wandering
interval $I$ on which the nondifferentiability set has dimension
$2/(r+1)$. By \eqref{gl:exactN}, its singular points and their iterates
are the sets $X_j(E)$. By Lemma~\ref{gl:local-affine}, $\phi_r$ is affine in a neighborhood of every such point.
The example therefore belongs to $\mathcal C_r(\alpha)$, proving
that the maximum is attained.

Suppose that the same perturbation belonged to $C^s$ for some
$s>r$. The graph, its wandering interval, its rotation number, and
its local affine property would be unchanged. Applying
Proposition~\ref{sharp:upper} with exponent $s$ would give
\[
 \frac2{r+1}
 =\dim_H\bigl(N(u)\cap I\bigr)
 \le\frac2{s+1}<\frac2{r+1},
\]
a contradiction. This proves \eqref{sharp:exact-regularity}.
\end{proof}


\par\medskip\noindent\textbf{Data availability.}
No datasets were generated or analyzed during the current study.
\par\medskip\noindent\textbf{Conflict of interest.}
The author declares no conflict of interest.
\par\medskip\noindent\textbf{AI declaration.}
The author used OpenAI's ChatGPT as an auxiliary tool for proofreading and improving the presentation of some arguments.
All mathematical arguments were independently verified by the author, who takes full
responsibility for the content.


\begin{thebibliography}{99}
\setlength{\itemsep}{1pt}

\bibitem{Ar09}
M.-C.~Arnaud, Three results on the regularity of the curves that are invariant by an exact symplectic twist map,
\emph{Publ. Math. Inst. Hautes \'Etudes Sci.} \textbf{109} (2009), 1--17.

\bibitem{Ar11}
M.-C.~Arnaud, A non-differentiable essential irrational invariant curve for a $C^1$ symplectic twist map,
\emph{J. Mod. Dyn.} \textbf{5} (2011), 583--591.

\bibitem{ArAnn11}
M.-C.~Arnaud, The link between the shape of the irrational
Aubry--Mather sets and their Lyapunov exponents,
\emph{Ann. of Math.} \textbf{174} (2011), 1571--1601.

\bibitem{Ar14}
M.-C.~Arnaud, Boundaries of instability zones for symplectic twist maps,
\emph{J. Inst. Math. Jussieu} \textbf{13} (2014), 19--41.

\bibitem{AB16}
M.-C.~Arnaud and P.~Berger, The non-hyperbolicity of irrational
invariant curves for twist maps and all that follows,
\emph{Rev. Mat. Iberoam.} \textbf{32} (2016), no.~4, 1295--1310.

\bibitem{ALC17}
M.-C.~Arnaud and P.~Le~Calvez, A notion of Denjoy sub-system,
\emph{C. R. Math. Acad. Sci. Paris} \textbf{355} (2017), 914--919.

\bibitem{Arn63}
V.~I.~Arnol'd, Proof of a theorem of A.~N.~Kolmogorov on the invariance
of quasi-periodic motions under small perturbations of the Hamiltonian,
\emph{Russian Math. Surveys} \textbf{18} (1963), no.~5, 9--36.

\bibitem{ALD83}
S.~Aubry and P.~Y.~Le~Daeron, The discrete Frenkel--Kontorova model
and its extensions. I. Exact results for the ground-states,
\emph{Physica D} \textbf{8} (1983), 381--422.

\bibitem{AF}
A.~Avila and B.~Fayad, Non-differentiable irrational curves for $C^1$ twist map,
\emph{Ergodic Theory Dynam. Systems} \textbf{42} (2022), 491--499.

\bibitem{Ban88}
V.~Bangert, Mather sets for twist maps and geodesics on tori,
in \emph{Dynamics Reported}, Vol.~1,
U.~Kirchgraber and H.-O.~Walther (eds.),
Wiley, Chichester, and Teubner, Stuttgart, 1988, 1--56.

\bibitem{BR}
M.~Boshernitzan and D.~Ralston, Continued fractions and heavy sequences,
\emph{Proc. Amer. Math. Soc.} \textbf{137} (2009), 3177--3185.

\bibitem{FK}
B.~Fayad and R.~Krikorian, Some questions around quasi-periodic dynamics,
in \emph{Proceedings of the International Congress of Mathematicians---Rio de Janeiro 2018}, Vol.~III, World Scientific, 2018, 1909--1932.

\bibitem{GLW}
J.~Guo, Y.~Liu and L.~Wang, On non-KAM invariant circles for area-preserving twist maps, preprint (2025).
{arXiv:2510.21098}.

\bibitem{Hall83}
G.~R.~Hall, Bifurcation of an attracting invariant circle:
a Denjoy attractor,
\emph{Ergodic Theory Dynam. Systems} \textbf{3} (1983), 87--118.

\bibitem{Hensley}
{D.~Hensley, Continued fraction Cantor sets, Hausdorff dimension, and
functional analysis, \emph{J. Number Theory} \textbf{40} (1992), 336--358.}

\bibitem{Her79}
M.~R.~Herman, Sur la conjugaison diff\'erentiable des diff\'eomorphismes
du cercle \`a des rotations, \emph{Publ. Math. Inst. Hautes \'Etudes Sci.}
\textbf{49} (1979), 5--233.

\bibitem{Her}
M.~R.~Herman, \emph{Sur les courbes invariantes par les diff\'eomorphismes de l'anneau}, Vol.~I,
Ast\'erisque \textbf{103--104}, Soci\'et\'e Math\'ematique de France, 1983, 1--221.

\bibitem{Her86}
M.~R.~Herman, \emph{Sur les courbes invariantes par les
diff\'eomorphismes de l'anneau}, Vol.~II,
Ast\'erisque \textbf{144}, Soci\'et\'e Math\'ematique de France, 1986.

\bibitem{HLW26}
Y.~Huang, J.~Li and L.~Wang,
On the sharpness of the $C^1$-norm threshold for perturbations in the
normally hyperbolic invariant manifold theorem---a toy model perspective,
preprint, arXiv:2608.04862 (2026).

\bibitem{Kol54}
A.~N.~Kolmogorov, On conservation of conditionally periodic motions
for a small change in Hamilton's function,
\emph{Dokl. Akad. Nauk SSSR} \textbf{98} (1954), 527--530.

\bibitem{KMW26}
D.~Kong, J.~E.~Massetti and L.~Wang,
Quantitative persistence and regularity of invariant graphs for
dissipative twist maps,
preprint, arXiv:2610.35068 (2026).

\bibitem{LeC88}
P.~Le~Calvez, Propri\'et\'es des attracteurs de Birkhoff,
\emph{Ergodic Theory Dynam. Systems} \textbf{8} (1988), 241--310.

\bibitem{LW26}
Q.~Li and L.~Wang, On the dynamics of invariant graphs for
dissipative twist maps,
preprint, arXiv:2506.04938v3 (2026).

\bibitem{MY}
S.~Marmi and J.-C.~Yoccoz, Some open problems related to small divisors,
in \emph{Dynamical Systems and Small Divisors} (Cetraro, 13--20 June 1998),
Lecture Notes in Math., vol.~1784, Springer, Berlin, 2002, 175--191.

\bibitem{Mat82}
J.~N.~Mather, Existence of quasi-periodic orbits for twist
homeomorphisms of the annulus,
\emph{Topology} \textbf{21} (1982), 457--467.

\bibitem{Mat84}
J.~N.~Mather, Non-existence of invariant circles,
\emph{Ergodic Theory Dynam. Systems} \textbf{4} (1984), 301--309.

\bibitem{Mos62}
J.~Moser, On invariant curves of area-preserving mappings of an annulus,
\emph{Nachr. Akad. Wiss. G\"ottingen Math.-Phys. Kl.~II}
(1962), 1--20.

\bibitem{Sinai}
Ya.~G.~Sinai, Some problems in the theory of dynamical systems and mathematical physics,
\emph{Geom. Funct. Anal.}, GAFA 2000 Special Volume, Part~I (2000), 425--433.
Lecture at \emph{Visions in Mathematics}, Tel Aviv, 25 August--3 September 1999.

\bibitem{SW26}
A.~Sorrentino and L.~Wang,
On the destruction of invariant Lagrangian graphs for conformal
symplectic twist maps,
\emph{Calc. Var. Partial Differential Equations}
\textbf{65} (2026), article~153.

\end{thebibliography}
\end{document}